\documentclass[12pt]{article}
\usepackage{amsmath,amssymb,amsthm,bbm}
\usepackage[colorlinks]{hyperref}
\usepackage[shortlabels]{enumitem}
\numberwithin{equation}{section}

\newcommand{\dto}{\stackrel{d}{\longrightarrow}}
\newcommand{\Pto}{\stackrel{\pr}{\longrightarrow}}
\newcommand{\asto}{\stackrel{\mathrm{a.s.}}{\longrightarrow}}
\newcommand{\vto}{\stackrel{v}{\longrightarrow}}

\newcommand{\wto}{\stackrel{w}{\longrightarrow}}
\newcommand{\fidi}{\stackrel{\mathrm{fi.di.}}{\longrightarrow}}
\newcommand{\toi}{\to\infty}
\newcommand{\eind}{\stackrel{d}{=}}

\newcommand{\dvague}{\stackrel{vd}{\longrightarrow}}

\newcommand{\pr}{\mathbb{P}}
\newcommand{\ex}{\mathbb{E}}
\newcommand{\var}{\text{Var}}
\DeclareMathOperator{\Var}{Var}
\newcommand{\law}{\mathcal{L}}
\newcommand{\diag}{\operatorname{diag}}

\newcommand{\sgn}{\operatorname{sgn}}

\newcommand{\dd}{\mathrm{d}}
\newcommand{\ee}{\mathrm{e}}
\newcommand{\ii}{\mathrm{i}}

\newcommand{\Rset}{\mathbb{R}}
\newcommand{\Nset}{\mathbb{N}}
\newcommand{\Zset}{\mathbb{Z}}

\newcommand{\C}{\mathbb{C}}

\newcommand{\N}{\mathbb{N}}
\newcommand{\R}{\mathbb{R}}

\newcommand{\Z}{\mathbb{Z}}

\newcommand{\NN}{\mathbb{N}}
\newcommand{\RR}{\mathbb{R}}
\newcommand{\ZZ}{\mathbb{Z}}
\newcommand{\QQ}{\mathbb{Q}}

\newcommand{\cA}{\mathcal{A}}
\newcommand{\cB}{\mathcal{B}}
\newcommand{\cF}{\mathcal{F}}
\newcommand{\cL}{\mathcal{L}}
\newcommand{\cN}{\mathcal{N}}

\newcommand{\hC}{\widehat{C}}

\newcommand{\tX}{\widetilde{X}}
\newcommand{\tY}{\widetilde{Y}}
\newcommand{\tZ}{\widetilde{Z}}
\newcommand{\tmu}{\widetilde{\mu}}
\newcommand{\ttheta}{\widetilde{\theta}}
\newcommand{\tkappa}{\widetilde{\kappa}}

\newcommand{\vare}{\epsilon}
\newcommand{\bzero}{{\boldsymbol{0}}}

\newcommand{\1}[1]{\mathbbm{1}_{#1}}

\newcommand{\bC}{{\boldsymbol{C}}}

\newcommand{\bX}{{\boldsymbol{X}}}

\newcommand{\bY}{\boldsymbol{Y}}
\newcommand{\bU}{\boldsymbol{U}}
\newcommand{\bV}{\boldsymbol{V}}

\newcommand{\bc}{{\boldsymbol{c}}}
\newcommand{\bu}{\boldsymbol{u}}
\newcommand{\bv}{\boldsymbol{v}}
\newcommand{\bx}{\boldsymbol{x}}
\newcommand{\by}{\boldsymbol{y}}
\newcommand{\bz}{\boldsymbol{z}}

\newcommand{\bone}{\boldsymbol{1}}

\newcommand{\hmuA}{\widehat{\mu_A}}
\newcommand{\hphi}{\widehat{\phi}}

\newcommand{\dsum} {\displaystyle\sum}

\newcommand{\proofend}{\hfill\mbox{$\Box$}}
\newcommand{\stoch}{\stackrel{\pr}{\longrightarrow}}
\newcommand{\qmean}{\stackrel{L_2}{\longrightarrow}}

\newtheorem{theorem}{Theorem}[section]
\newtheorem{lemma}[theorem]{Lemma}

\newtheorem{proposition}[theorem]{Proposition}

\theoremstyle{definition}
\newtheorem{remark}[theorem]{Remark}

\begin{document}

\begin{center}
 {\bfseries\large Statistical inference of subcritical strongly stationary  \\[2mm]
                   Galton--Watson processes with regularly varying immigration} \\[5mm]
 {\sc M\'aty\'as $\text{Barczy}^{*,\diamond}$, \ Bojan $\text{Basrak}^{**}$, \ P\'eter $\text{Kevei}^{***}$,} \\[2mm]
  {\sc \framebox[1.1\width]{Gyula $\text{Pap}^{***,****}$}, \ Hrvoje $\text{Planini\'c}^{**}$}
\end{center}

\vskip0.2cm

\noindent
 * MTA-SZTE Analysis and Stochastics Research Group,
   Bolyai Institute, University of Szeged,
   Aradi v\'ertan\'uk tere 1, 6720 Szeged, Hungary.

\noindent
 ** Department of Mathematics,
    Faculty of Science, University of Zagreb,
    Bijeni\v{c}ka 30, 10000 Zagreb, Croatia

\noindent
 *** Bolyai Institute, University of Szeged,
     Aradi v\'ertan\'uk tere 1, 6720 Szeged, Hungary.

\noindent
 **** Alfr\'ed R\'enyi Institute of Mathematics,
      Re\'altanoda u. 13-15, 1053 Budapest, Hungary.

\noindent e-mails: barczy@math.u-szeged.hu (M. Barczy),
                   bbasrak@math.hr (B. Basrak),
                   kevei@math.u-szeged.hu (P. Kevei),
                   hrvoje.planinic@math.hr (H. Planini\'c).

\noindent $\diamond$ Corresponding author.

\begin{abstract}
We describe the asymptotic behavior of the conditional least squares estimator of the offspring mean
 for subcritical strongly stationary Galton--Watson processes with regularly varying immigration with tail index $\alpha \in (1, 2)$.
The limit law is the ratio of two dependent stable random variables with indices $\alpha/2$ and $2\alpha/3$, respectively, and it has a continuously differentiable density function.
We use point process technique in the proofs.
\end{abstract}

\renewcommand{\thefootnote}{}
\footnote{Supported by the Hungarian Croatian Intergovernmental S\&T Cooperation Programme
 for 2017-2018 under Grant No.\ 16-1-2016-0027.
M\'aty\'as Barczy and P\'eter Kevei are supported by the J\'anos Bolyai Research Scholarship of the Hungarian
 Academy of Sciences.
Bojan Basrak, P\'eter Kevei and Hrvoje Planini\'c are financed within the Croatian-Swiss Research Program
 of the Croatian Science Foundation and the Swiss National Science Foundation - grant CSRP 018-01-180549.
P\'eter Kevei and Gyula Pap are supported by the Ministry for Innovation and Technology, Hungary grant TUDFO/47138-1/2019-ITM,
 and by the EU-funded Hungarian grant EFOP-3.6.1-16-2016-00008.
P\'eter Kevei is supported by the NKFIH grant FK124141.}

\vspace*{-10mm}

\section{Introduction}

The theory and estimation of branching processes, especially Galton--Watson processes without or with immigration has a long history, see, e.g.,
 the survey paper by Winnicki \cite{winnicki:1987}.
In this paper we will consider a Galton--Watson process with regularly varying immigration distribution, and we will study the
 Conditional Least Squares (CLS) estimation of the mean of the offspring distribution.
Heavy-tailed Galton--Watson processes with immigration, especially with regularly varying immigration distribution,
 have been in the focus of research for a long time.
We name only two papers here.
Seneta \cite{seneta:1973} derived conditions under which there exists a sequence of positive constants such that the logarithm of a Galton--Watson process
 having offspring distribution with infinite mean and normalized by the sequence in question converges in distribution to a non-degenerate distribution.
Schuh and Barbour \cite{schuh:barbour:1977} derived necessary and sufficient conditions for the almost sure convergence
 of some slowly varying function of a Galton-Watson process having offspring distribution with infinite mean.
Heavy-tailed branching processes are important not only from the theoretical point of view, but they are also used in modeling of biological phenomena,
 e.g., for modeling the development of multi-focal tumors, see Ernst et al. \cite{ErnKimKurZho}.
Concerning the estimation theory of heavy-tailed Galton--Watson processes with immigration we are not aware of any results.
Let us recall now a result on the estimation of the offspring mean under finite third order moment assumptions.

Let $\Zset$ and $\Nset$ denote the set of integers and positive integers, respectively.
Every random variable will be defined on a fixed probability space
 $(\Omega, \cA, \pr)$.
For each $i, j \in \Nset$, the number of individuals in the
 $i^\mathrm{th}$ generation will be denoted by $X_i$, the number of
 offspring produced by the $j^\mathrm{th}$ individual belonging to the
 $(i-1)^\mathrm{th}$ generation will be denoted by $A_{j}^{(i)}$, and
 the number of immigrants in the \ $i^\mathrm{th}$ \ generation
 will be denoted by \ $B_i$.
Further, $X_0$ \ denotes the size of the initial population.
Then we have
 \begin{align*}
  X_i = \sum_{j=1}^{X_{i-1}} A_{j}^{(i)} + B_i, \qquad i\in\Nset,
 \end{align*}
where $\{A, A_{j}^{(i)} : i, j\in\Nset\}$ are independent, identically distributed (i.i.d.) nonnegative integer-valued
 random variables independent of another i.i.d.\ sequence $\{B, B_i : i\in\Nset\}$ of nonnegative integer-valued random variables.
Assuming that $X_0$ is independent of $\{A_{j}^{(i)}, B_i : i, j\in\Nset\}$, $\ex[X_0] < \infty$, $\mu_A := \ex[A] < \infty$, $\mu_B := \ex[B] \in (0, \infty)$ and that $\mu_B$ is known, the CLS estimator of $\mu_A$ based on the observations $X_0, X_1, \dots, X_n$
 has the form
 \begin{equation}\label{CLSE}
  \hmuA^{(n)} := \frac{\sum_{i=1}^n X_{i-1} (X_i - \mu_B)}{\sum_{i=1}^n X_{i-1}^2}
 \end{equation}
 on the set $\left\{\sum_{i=1}^n X_{i-1}^2 > 0\right\}$, see Klimko and Nelson \cite{klimko:nelson:1978}.
We have $\pr(\sum_{i=1}^n X_{i-1}^2 > 0) \to 1$ as $n \toi$, since $\pr(\sum_{i=1}^n X_{i-1}^2 = 0) = \pr(X_0 = 0, B_1 = 0, \ldots, B_{n-1} = 0) \leq \pr(B = 0)^{n-1} \to 0$ as $n \toi$ due to $\pr(B = 0) \in [0, 1)$.
If, in addition, $\mu_A \in (0, 1)$, then the Markov chain $(X_i)_{i\geq0}$ admits a unique
 stationary distribution, see, e.g., Quine \cite{quine:1970a}.
If, in addition, $\ex[X_0^3] < \infty$, $\ex[A^3] < \infty$ and $\ex[B^3] < \infty$, then
 \begin{equation}\label{m_xi_sub1}
   n^{1/2} (\hmuA^{(n)} - \mu_A)
   \dto
   \cN\biggl(0, \frac{\sigma_A^2\ex[\tX^3]+\sigma_B^2\ex[\tX^2]}{\bigl(\ex[\tX^2]\bigr)^2}\biggr)
   \qquad \text{as \ $n \to \infty$,}
 \end{equation}
 where $\sigma_A^2 := \Var(A)$ and $\sigma_B^2 := \Var(B)$, the distribution of the random variable $\tX$ is the unique stationary distribution of $(X_i)_{i\geq0}$, and $\dto$ denotes convergence in distribution.
The paper by Klimko and Nelson \cite[Section 5]{klimko:nelson:1978} contains a similar result for the CLS estimator
 of $(\mu_A, \mu_B)$, and \eqref{m_xi_sub1} can be derived by the method of
 that paper.
Note that $\ex[\tX^2]$ and $\ex[\tX^3]$ can be expressed in terms of the first three
 moments of $A$ and $B$,
 see, e.g., Quine \cite[formula (26) and page 422]{quine:1970a} and Barczy et al.~\cite[formulae (14), (16) and (20)]{barczy:nedenyi:pap:2018}.

In contrast with those earlier results, we explore the case  where the distribution of $B$ is regularly varying with tail index in $(1, 2)$, thus having infinite variance.
In the sequel, we will always assume the following conditions:
 \begin{itemize}
  \item[(i)] $\mu_A \in (0, 1)$,
  \item[(ii)] $\sigma_A^2\in(0,\infty)$,
  \item[(iii)] $B$ is regularly varying with tail index $\alpha \in (1, 2)$, i.e.,
               \begin{align*}
                  \lim_{x\to\infty}
                  \frac{\pr(B>qx)}{\pr(B>x)} = q^{-\alpha} \qquad \text{for all $q \in (0, \infty)$.}
               \end{align*}
 \end{itemize}
In particular, $\mu_B \in (0, \infty)$, and there exists a strongly stationary process $(X_i)_{i\in\Zset}$ satisfying
\begin{equation}\label{GWI}
X_i=\sum_{j=1}^{X_{i-1}} A_{j}^{(i)} + B_i, \qquad i \in \Zset,
\end{equation}
where $\{A,A_{j}^{(i)}:j\in\Nset\,,i\in\Zset\}$ are i.i.d.\ nonnegative integer-valued random variables independent of another i.i.d.\ sequence
 $\{B,B_i:i\in\Zset\}$ of nonnegative integer-valued random variables.
Note that, for simplicity, in \eqref{GWI} we keep the same letter $X$ to denote the strongly stationary extension of the process indexed by $\ZZ$.
In this case, the distribution of $X_0$ is also regularly varying with the same tail index $\alpha$ having infinite variance, or more precisely,
\begin{align}\label{eq:tail_of_X}
\pr(X_0>x)\sim \frac{1}{1-\mu_A^{\alpha}} \pr(B>x) \qquad \text{as $x\toi$,}
\end{align} see Basrak et al.~\cite[Theorem 2.1.1]{basrak:kulik:palmowski:2013}.

Our aim is to study the limiting behavior of $\hmuA^{(n)}$ as $n \to\infty$ for the strongly stationary process $(X_i)_{i\in\ZZ}$ given in \eqref{GWI}.
For each $n \in \Nset$, by \eqref{CLSE}, we have
 \[
   \hmuA^{(n)} - \mu_A
   = \frac{\sum_{i=1}^n X_{i-1} (X_i - \mu_B)}{\sum_{i=1}^n X_{i-1}^2} - \mu_A
   = \frac{\sum_{i=1}^n X_{i-1} M_i}{\sum_{i=1}^n X_{i-1}^2}
 \]
 on the set $\left\{\sum_{i=1}^n X_{i-1}^2 > 0\right\}$, where, by \eqref{GWI},
 \begin{align}\label{def:M}
   M_i
   := X_i - \mu_A X_{i-1} -\mu_B
   = \sum_{j=1}^{X_{i-1}} (A_j^{(i)} - \mu_A) + (B_i-\mu_B)
   =: \sum_{j=1}^{X_{i-1}} \tilde{A}_j^{(i)} + \tilde{B}_i, \qquad i \in \Zset .
 \end{align}
Intuitively, by the central limit theorem, for large $X_{i-1}$, the distribution of $M_i/\sqrt{X_{i-1}}$ is approximately normal,
 thus, in the spirit of Breiman's lemma, $M_i X_{i-1} = X_{i-1}^{3/2}\, M_i/\sqrt{X_{i-1}}$ is regularly varying with tail index $2\alpha/3$.
This argument is made precise in Proposition \ref{prop:pseudoTail}.

Our analysis relies on the fact that one can determine the weak limit of the point processes
 \begin{equation}\label{conv_point_process}
  \sideset{}{^{*}} \sum_{j=1}^n \delta_{\bigl(\frac{X_j}{a_n},\,\frac{M_{j+1}}{\sqrt{ X_j}}\bigr)}
  := \sum_{\{ j\in\{1,\ldots,n\}\,:\,X_j>0\}} \delta_{\bigl(\frac{X_j}{a_n},\,\frac{M_{j+1}}{\sqrt{X_j}}\bigr)}
 \end{equation}
as $n\to\infty$ on a suitable space of point measures on $(0,\infty)\times\Rset$ and with a scaling sequence $(a_n)_{n\in\Nset}$ satisfying $n \pr(X_0 > a_n) \to 1$ as $n \toi$
 (see \eqref{a_n_quantile}),
 where $\delta_{(x,y)}$ denotes the Dirac measure concentrated on $(x,y)\in (0,\infty)\times \R$,
 see Theorem \ref{th:PPconv}.
For a possible choice of a suitable sequence $(a_n)_{n\in\Nset}$ and its asymptotic behavior as $n\to\infty$, see the beginning of Section \ref{sec:Ppc}.
The proof of Theorem \ref{th:PPconv} is based on general results of Kallenberg \cite[Theorems 4.11 and 4.22]{kallenberg:2017}
 for convergence in distribution of random measures with respect to the vague topology.
Point processes have been  often applied to analyze regularly varying observations, see, for instance, Resnick \cite{resnick:1987,resnick:2007} and
 Kulik and Soulier~\cite{kulik:soulier:2020}, however, our approach here is not standard since we topologize the space of point measures on $(0,\infty)\times\Rset$ using vague topology with "bounded" Borel sets being those which are bounded away from the vertical line $\{0\} \times \mathbb{R}$ instead of being bounded away from the point $(0,0)$.
For a detailed discussion on our setup, see the beginning of Section \ref{sec:Ppc} and Appendix \ref{app:topology}.
In the course of the proof of Theorem \ref{th:PPconv}, the joint tail behavior of $(X_i)_{i\in\ZZ}$ and $(M_{i+1}/\sqrt{\max(X_i,1)})_{i\in\ZZ_+}$,
 especially the so-called tail process of \ $(X_i)_{i\in\ZZ}$, \ also plays a crucial role, see Proposition \ref{prop:pseudoTail} and \eqref{tail_process}.

Based on convergence of the point processes in \eqref{conv_point_process}, we obtain
 \[
   \sqrt{a_n}(\hmuA^{(n)} - \mu_A) \dto \frac{V^{(2)}}{V^{(1)}} \qquad \text{as $n \toi$,}
 \]
 where $V^{(1)}$ is an $\alpha/2$-stable positive random variable, $V^{(2)}$ is a symmetric $2\alpha/3$-stable random variable, and $V^{(1)}$ and $V^{(2)}$ are dependent with an explicitly given joint characteristic function, see Theorem \ref{th:main}.
Concerning the asymptotic behavior of $(a_n)_{n\in\NN}$, note that if $x^\alpha \pr(B > x) \to 1$ as $x \toi$, i.e.,
 the distribution of $B$ is asymptotically equivalent to a Pareto distribution with parameter $\alpha$,
 then $n^{-1/\alpha} a_n \to (1-\mu_A^\alpha)^{-1/\alpha}$ as $n \toi$.
Indeed, using \eqref{eq:tail_of_X} and that $n\pr(X_0>a_n)\to 1$ as $n\toi$, we have
 \begin{align*}
  n^{-1/\alpha} a_n  = (n^{-1}a_n^\alpha)^{1/\alpha}
    \sim (a_n^\alpha \pr(X_0>a_n) )^{1/\alpha}
    \sim \big( a_n^\alpha (1-\mu_A^\alpha)^{-1} \pr(B>a_n) \big)^{1/\alpha}
    \sim (1-\mu_A^\alpha)^{-1/\alpha}
 \end{align*}
 as $n\to\infty$, as desired.

In Section \ref{section:properties}, we collect several properties of $(V^{(1)}, V^{(2)})$ and $V^{(2)}/V^{(1)}$, including that the distribution of $(V^{(1)}, V^{(2)})$ is operator stable and $V^{(2)}/V^{(1)}$ has a continuously differentiable density function.
In Appendix \ref{app:topology}, we collect some topological properties of $(0,\infty)\times \R$.
Appendix \ref{app:vague_convergence_of_point_measures} contains the proof of Lemma \ref{lem:vagueConv_pointMeas} describing vague convergence of point measures.
In Appendix \ref{app:laplace} we show that the process $\bigl(X_i\1{\{X_i>0\}},\frac{M_{i+1}}{\sqrt{X_i}}\1{\{X_i>0\}}\bigr)_{i\geq0}$ satisfies a certain mixing condition.
Appendix \ref{app:slutsky+CMT} is devoted to a conditional Slutsky's lemma and a conditional continuous mapping theorem.
In Appendix \ref{app:regvarrelproc}, we show that the process $(X_i^{3/2}, X_i M_{i+1})_{i\in\Zset}$ is regularly varying with tail index $2\alpha/3$ with an explicitly given forward tail process.

We note that our proof technique does not work for the case $\alpha=2$.
Formally, one can see that the fact $\alpha<2$ is used many times in the proofs,
 and, for example, in case $\alpha=2$ the series $\sum_{i=1}^\infty P_i^2$ appearing as the weak limit of $(1-\mu_A^2)a_n^{-2}\sum_{j=1}^n X_j^2$ as $n\to\infty$
 in Theorem \ref{thm:partial_sums} is not convergent almost surely by Campbell's theorem (see, e.g., Kingman \cite[Section 3.2]{kingman:1993})
 due to $\int_0^\infty (y^2\wedge 1)\,\dd (-y^{-2}) = \infty$.
In case of $\alpha=2$, $X_j^2$, $j\in\NN$, is positive almost surely and regularly varying with tail index $1$,
 so in order to get some weak limit of $\sum_{j=1}^n X_j^2$ as $n\to\infty$ one may need to introduce an appropriate centering as well.

Our results can be compared with the results on AR(1) processes
 \[
   \xi_i = \phi \xi_{i-1} + \vare_i , \qquad i \in \Nset ,
 \]
 where $\phi\in\RR$ and $\{\vare, \vare_i : i \in \Nset\}$ are i.i.d.\ random variables.
Assuming that $\xi_0$ is independent of $\{\vare_i : i \in\Nset\}$, $\ex[\vare] = 0$, the CLS estimator of $\phi$ based on the observations
 $\xi_0, \xi_1, \dots, \xi_n$
 has the form
 \begin{equation*}
  \hphi_n := \frac{\sum_{i=1}^n \xi_{i-1} \xi_i}{\sum_{i=1}^n \xi_{i-1}^2}
 \end{equation*}
 on the set $\left\{\sum_{i=1}^n \xi_{i-1}^2 > 0\right\}$.
If, in addition, $\phi \in (-1, 1)$, then the Markov chain $(\xi_i)_{i\geq0}$ admits a unique
 stationary distribution, and there exists a strongly stationary process $(\xi_i)_{i\in\Zset}$ satisfying
\begin{equation*}
\xi_i = \phi \xi_{i-1} + \vare_i , \qquad i \in \Zset ,
\end{equation*}
where $\{\vare, \vare_i : i \in \Zset\}$ are i.i.d.\ random variables.
If, in addition, $\vare$ is symmetric and regularly varying with tail index $\alpha \in (0, 2)$,
 then
 \[
   \hphi_n \to \phi \qquad \text{as $n \toi$ almost surely,}
 \]
 see Hannan and Kanter \cite{hannan:kanter:1977}, and
 \[
   b_n (\hphi_n - \phi) \dto \frac{U^{(2)}}{U^{(1)}} \qquad \text{as $n \toi$}
 \]
 under the additional assumption $\lim_{x\to\infty}\frac{\pr(|\vare_0\vare_1|>x)}{\pr(|\vare_0|>x)}=2\ex[|\vare|^\alpha]$ in case of
 $\ex[|\vare|^\alpha]<\infty$,
 where $(b_n)_{n\in\Nset}$ is a suitable scaling sequence, $U^{(1)}$ is an $\alpha/2$-stable positive random variable and $U^{(2)}$ is a symmetric $\alpha$-stable random variable, see Davis and Resnick \cite[Theorem 3.6]{davis:resnick:1985b} for the case of $\ex[|\vare|^\alpha] < \infty$ and Davis and Resnick \cite[Theorem 4.4]{davis:resnick:1986} for the case of $\ex[|\vare|^\alpha] = \infty$.
This representation of the limit distribution is derived in Davis and Resnick \cite[Example 5]{davis:resnick:1986}.
Further, $U^{(1)}$ and $U^{(2)}$ are claimed to be dependent in case of $\ex[|\vare|^\alpha] < \infty$,
while they are independent in case of $\ex[|\vare|^\alpha] = \infty$.
In fact, if $\ex[|\vare|^\alpha] < \infty$, then $b_n$, $n\geq 2$, is the $1-\frac{1}{n}$ lower quantile of $|\vare|$, and one can
 also write $b_n = n^{\frac{1}{\alpha}}L_1(n)$, $n\in\NN$, with some slowly varying function $L_1:(0,\infty)\to(0,\infty)$.
If $\ex[|\vare|^\alpha] = \infty$, then $b_n = \widetilde c_n^{-1}c_n^2$, $n\geq 2$, where $c_n$ is the $1-\frac{1}{n}$ lower quantile of $|\vare_1|$,
 and $\widetilde c_n$ is the $1-\frac{1}{n}$ lower quantile of $|\vare_0\vare_1|$, and in this case one can also write
 $b_n = n^{\frac{1}{\alpha}}L_2(n)$, $n\in\NN$, with some slowly varying function $L_2:(0,\infty)\to(0,\infty)$.
Moreover, if $x^\alpha \pr(|\vare| > x) \to 1$ as $x \toi$, i.e., the distribution of $|\vare|$ is asymptotically equivalent to a Pareto distribution with
 parameter $\alpha$, then $\ex[|\vare|^\alpha] = \infty$ and $n^{-1/\alpha} (\log(n))^{1/\alpha} b_n \to 1$ as $n \toi$, see Resnick \cite[Problem 9.13]{resnick:2007}.
We also emphasize that Galton--Watson processes with immigration are quite different from AR(1) process due to the branching property of the process.
As a consequence, for example, the unique stationary distribution of the process can be represented in a more complicated way compared to that of an AR(1) process,
 see Lemma \ref{lem:reprX}.
So our point process technique for the proof is not a simple modification of the known one for AR(1) processes.

Finally, we recall two results on the CLS estimator of some parameters for related heavy-tailed continuous time processes.

Hu and Long \cite{hu:long:2009} studied the asymptotic behavior of the least squares estimator of the drift parameter
 for a generalized Ornstein-Uhlenbeck process driven by a symmetric $\alpha$-stable L\'evy motion with $\alpha\in(0,2)$ in the ergodic case
 based on discrete time infill-increasing (high frequency) observations.
Using some results of Davis and Resnick \cite{davis:resnick:1986}, Hu and Long \cite{hu:long:2009} proved strong consistency of the least squares estimator in question,
 and they also described its asymptotic behavior with a limit distribution being the fraction of two independent stable random variables.

Li and Ma \cite{li:ma:2015}, using a similar point process technique, described somewhat similar asymptotic behavior of the weighted and non-weighted CLS estimators of the drift parameters for a stable Cox--Ingersoll--Ross model based on low frequency observations.
This process can be viewed as a special subcritical continuous state and continuous time branching process with immigration.

Convergence in probability and in $L_2$ will be denoted by $\stoch$ and $\qmean$, respectively.
Weak convergence of finite measures will be denoted by $\wto$, and convergence of finite dimensional distributions is denoted by $\fidi$.
The diagonal matrix with diagonal entries $a_1,\ldots,a_n\in\R$ is denoted by $\diag_n(a_1,\dots,a_n)$.
We write $\delta_h$ for the Dirac measure on a set $H$ concentrated at $h \in H$.
For a random vector $\bX$ and event $A \in \cA$ such that $\pr(A)>0$, $\law(\bX)$ and $\law (\bX \,|\, A)$ denote the law of $\bX$ and law of $\bX$ conditionally on $A$, respectively.

\section{Tail behavior of the process}

In what follows, let us consider the strongly stationary process $(X_i)_{i\in\ZZ}$ given in \eqref{GWI}.
For determining the weak limit of the point processes given in \eqref{conv_point_process}, in the proof of Theorem \ref{th:PPconv}
 the joint tail behavior of $(X_i)_{i\in\ZZ}$ and $(M_{i+1}/\sqrt{\max(X_i,1)})_{i\in\ZZ_+}$,
 especially the so-called tail process of \ $(X_i)_{i\in\ZZ}$, \ plays a crucial role.
First, we present the tail process of $(X_i)_{i\in\ZZ}$, and then we formulate a result on the tail behavior of the above mentioned two stochastic processes.

The process $(X_i)_{i\in \Z}$ is jointly regularly varying with tail index $\alpha$ and admits a tail process $(Y_i)_{i\in \Z}$, i.e.\
 all the finite dimensional distributions of  $(X_i)_{i\in \Z}$ are regularly varying with tail index $\alpha$, and for all $m\in \N$,
\begin{align*}
\law (x^{-1}X_{-m},\dots, x^{-1} X_m \mid X_0 >x) \wto \law(Y_{-m},\dots,Y_m) \, ,
\end{align*}
as $x\toi$.
Indeed, Basrak et al.~\cite[Lemma 3.1]{basrak:kulik:palmowski:2013} showed the existence of the forward tail process $(Y_i)_{i\geq 0}$
 of the sequence $(X_i)_{i\in\Zset}$, and, by Theorem 2.1 in Basrak and Segers~\cite{basrak:segers:2009}, the existence of the forward tail process
 is equivalent to the existence of the (whole) tail process and to the joint regular variation of  $(X_i)_{i\in\Zset}$ as well.
We claim that
 \[
   Y_{-K+i} = \left\{
  \begin{array}{ll}
   \mu_A^{-K+i} Y_0,   & i \geq 0, \\[2mm]
   0, &  i < 0 ,
  \end{array}
  \right.
 \]
 yielding
  \begin{align}\label{tail_process}
  Y_i = \left\{
  \begin{array}{ll}
   \mu_A^ i Y_0,   & i \geq 0, \\[2mm]
   \mu_A^ i \1{\{K\geq |i|\}} Y_0, &  i < 0 ,
  \end{array}
  \right.
 \end{align}
 where $Y_0$ is a Pareto distributed random variable such that $\pr(Y_0\geq y)=y^{-\alpha}$ for $y\geq 1$, and $K$ is a geometrically distributed random variable independent of $Y_0$ such that
 $$
   \pr (K=k) = \mu_A^{\alpha k}  (1-\mu_A^\alpha)\,, \quad k =0,1,2,\ldots
 $$
Indeed,
 as shown in Basrak et al.~\cite[Lemma 3.1]{basrak:kulik:palmowski:2013},  $(Y_i)_{i\geq 0}$ is the forward tail process of the sequence $(X_i)_{i\in\Zset}$.
On the other hand, by Janssen and Segers \cite[Example 6.2]{janssen:segers:2014}, $(Y_i)_{i\in\Zset}$ is the tail process of the stationary solution
 $(\tilde\xi_i)_{i\in \Z}$ to the stochastic recurrence equation $\tilde\xi_i=\mu_A \tilde\xi_{i-1} + B_i$, $i \in \Zset$.
Since the distribution of the forward tail process determines the distribution of the (whole) tail process (see Basrak and Segers~\cite[Theorem 2.1]{basrak:segers:2009}), it follows that $(Y_i)_{i\in\Zset}$ represents the tail process of $(X_i)_{i\in\Zset}$.

For the ease of notation, denote
$$
  W_i := \frac{M_{i+1}}{\sqrt{\mu_A^iX_0}} \,, \qquad
  W'_i := \frac{M_{i+1}}{\sqrt{X_i\lor1}} \,, \qquad i \geq 0 \,,
$$
on the set $\{X_0 > 0\}$, where $a \vee b := \max\{a, b\}$, $a, b \in \R$.
Note that $\pr(X_0 = 0) > 0$ might occur.
For example, if $\pr(B \geq k) = c k^{-\alpha}$ for $k \in \Nset$ and $\pr(B = 0) = 1 - c$ with some \ $c \in (0, 1)$, then $\mu_B = c \sum_{k=1}^\infty k^{-\alpha}$, hence in case of $c < \frac{1-\mu_A}{\sum_{k=1}^\infty k^{-\alpha}}$, we have $\ex(X_0) = \frac{\mu_B}{1-\mu_A} < 1$, yielding $\pr(X_0 = 0) > 0$, since $X_0$ is a nonnegative integer valued random variable.

\begin{proposition}\label{prop:pseudoTail}
As $x\toi$,
\begin{equation}\label{eq:pseudoTail}
\law (x^{-1}X_{-m},\dots,x^{-1}X_{m},W_0',\dots, W_m' \mid X_0>x)
 \wto \law(Y_{-m},\dots,Y_{m},Z_0,\dots,Z_m)
\end{equation}
for all $m\in\N$, where $(Z_i)_{i\geq 0}$ is an i.i.d.\ sequence of $\cN(0,\sigma_A^2)$--distributed random variables
 being independent of $Y_0$ and $K$ with $\sigma^2_A =\Var(A) \in (0, \infty)$.
\end{proposition}

Although one can prove that the two-dimensional process $(X_i^{3/2}, X_i M_{i+1})_{i\in\ZZ}$ admits a tail process in the sense of Basrak and Segers \cite{basrak:segers:2009}, see Appendix \ref{app:regvarrelproc}, and then use standard point process convergence results for describing the asymptotic behavior of $\hmuA^{(n)}$ as $n\toi$, such an approach turns out to be rather complicated.
However, for the purpose of our analysis, the statement of Proposition \ref{prop:pseudoTail} turns out to be sufficient.
Note that this approach is similar to the so--called conditional extreme value approach, see Kulik and Soulier~\cite{kulik:soulier:2015} and references therein.

\noindent{\bf Proof of Proposition \ref{prop:pseudoTail}.}
Let $m\in\N$ be fixed.
First, note that $Y_0^{-1} Y_i= \mu_A^{i}$, $i\geq 0$, is the forward spectral process of $(X_k)_{k\in\Z}$, so by part (ii) of Corollary 3.2 in Basrak and Segers~\cite{basrak:segers:2009},
\begin{align}\label{eq:pseudoTail_inter1}
\cL(X_0^{-1}X_0, \ldots, X_0^{-1}X_m \mid X_0 > x) \wto \delta_{(1,\mu_A,\ldots,\mu_A^m)} \qquad \text{as $x \to \infty$.}
\end{align}
In particular, for every $i \geq 0$ and $\epsilon>0$, since convergence in distribution to a constant implies convergence in probability
 (formally applying part (ii) of Lemma \ref{lemma:cmt} with the Borel measurable function $h:\RR\to\RR$, $h(x):=\bone_{\{\vert \mu_A^{-i} x - 1\vert> \vare\}}$, $x\in\RR$,
 satisfying $D_h = \{\mu_A^i(1\pm \vare)\}$ and $\pr(\mu_A^i \in D_h)=0$ with $D_h$ being the set of discontinuities of $h$), we have
\begin{align}\label{eq:pseudoTail_inter0}
\lim_{x\toi} \pr(|(\mu_A^{i}X_0)^{-1}X_i - 1|>\epsilon  \mid X_0>x)=0 \, .
\end{align}
Moreover, it is enough to show \eqref{eq:pseudoTail} for $W_0$, \ldots, $W_m$ instead of $W'_0$, \ldots, $W'_m$.
Indeed, one can easily check that
 \[
   \cL(X_0^{-1} \mid X_0 > x) \wto \delta_0 \qquad \text{as $x \to \infty$,}
 \]
 hence, by \eqref{eq:pseudoTail_inter1} and Lemma \ref{lemma:Slutsky4}, we obtain
 \[
   \cL(X_0^{-1}, X_0^{-1}X_0, \ldots, X_0^{-1}X_m \mid X_0 > x) \wto \delta_{(0,1,\mu_A,\ldots,\mu_A^m)} \qquad \text{as $x \to \infty$.}
 \]
Then, identifying $\R^{(3m+2)\times(3m+2)}$ with $\R^{(3m+2)^2}$ in a natural way, we can use a conditional version of the continuous mapping theorem (see part (i) of Lemma \ref{lemma:cmt}), and we get
 \begin{align*}
  \cL\biggl(\diag_{3m+2}\biggl(1, \ldots, 1, \sqrt{\frac{X_0}{X_0\lor1}}, \ldots, \sqrt{\frac{\mu_A^m X_0}{X_m\lor1}}\biggr) \,\bigg|\, X_0 > x\biggr)
  \wto \delta_{\diag_{3m+2}(1,\ldots,1)}
 \end{align*}
 as $x \to \infty$.
Consequently, \eqref{eq:pseudoTail} with $W_i'$ replaced by $W_i$ and Lemma \ref{lemma:Slutsky3} imply \eqref{eq:pseudoTail}.

Define now for $i\geq 0$ and $n\geq 1$,
 \begin{align*}
 W''_i(n)  :=  (\mu_A^{i}n)^{-1/2} \left( \sum_{j=1}^{\lfloor \mu_A^{i}n\rfloor} \tilde{A}_j^{(i+1)} + \tilde{B}_{i+1}\right) ,
 \end{align*}
 where $\tilde{A}_j^{(i+1)}$ and $\tilde{B}_{i+1}$ are introduced in \eqref{def:M}.
Take arbitrary $y\geq 1 , u_{-m},\dots,u_m \in \R$ and introduce the events
\begin{align*}
 C(x) &:= \{X_{-m}\leq x u_{-m},\dots,X_{-1} \leq x u_{-1}\} , \qquad x \in (0, \infty), \\
 D(n) &:= \{W''_0(n) \leq u_0,\dots, W_m''(n) \leq u_m\} , \qquad n \in \Nset .
\end{align*}
Observe that the definition of the tail process $(Y_j)_{j\in\Zset}$ of $(X_j)_{j\in\Zset}$ implies
\begin{align}\label{help12}
\lim_{x\toi} \pr(X_0>xy , C(x) \mid X_0>x) = \pr(Y_0 > y, Y_{-m} \leq u_{-m},\dots , Y_{-1} \leq u_{-1}) \, .
\end{align}
Further,
 the central limit theorem, Slutsky's lemma and independence of $W''_0(n),\ldots,W''_m(n)$ imply that
 \begin{align}\label{help13}
 \begin{split}
  \lim_{n\toi} \pr(D(n))
  &=\lim_{n\toi} (\pr(W''_0(n) \leq u_0) \cdots \pr(W''_m(n) \leq u_m)) \\
  &= \pr(Z_0 \leq u_0) \cdots \pr(Z_m \leq u_m)
   = \pr(Z_0 \leq u_0,\dots, Z_m \leq u_m),
 \end{split}
 \end{align}
 where, recall, $Z_i$, $i\geq 0$, are independent and $\cN(0,\sigma_A^2)$--distributed random variables.
Since $(\tilde{A}_j^{(i+1)},\tilde{B}_{i+1} : i\geq 0 , j\geq 1 )$ are independent of $X_{-m},\dots, X_0$, \eqref{help12} and \eqref{help13}
 imply that for every $\epsilon>0$,
 \begin{align*}
 &\limsup_{x\toi} \pr(X_0>xy, C(x) , D(X_0) \mid X_0>x) \\
 &= \limsup_{x\toi}\pr(X_0>x)^{-1} \sum_{n>xy,\,n\in\NN} \pr(X_0 = n, C(x)) \pr(D(n))\\
 &\leq \limsup_{x\toi}\pr(X_0>x)^{-1} \sum_{n>xy,\,n\in\NN} \pr(X_0 = n, C(x)) (\pr(Z_0 \leq u_0,\dots, Z_m \leq u_m) + \epsilon)  \\
 &= (\pr(Z_0 \leq u_0,\dots, Z_m \leq u_m) + \epsilon)  \limsup_{x\toi} \pr(X_0 > x y, C(x) \mid X_0 > x) \\
 &= (\pr(Z_0 \leq u_0,\dots, Z_m \leq u_m) + \epsilon) \pr(Y_0 > y, Y_{-m} \leq u_{-m},\dots, Y_{-1} \leq u_{-1}) \, .
\end{align*}
An analogous claim holds for the limit inferior, and now letting $\epsilon\to 0$ yields that
\begin{multline*}
\lim_{x\toi} \pr(X_0>xy, C(x) , D(X_0) \mid X_0>x)\\
= \pr(Y_0 > y, Y_{-m} \leq u_{-m},\dots, Y_{-1} \leq u_{-1}, Z_0 \leq u_0,\dots, Z_m \leq u_m) \, , \,
\end{multline*}
since $(Z_j)_{j\geq0}$ is assumed to be independent of $(Y_j)_{j\in\Zset}$, hence
\begin{multline*}
 \law (x^{-1}X_{-m},\dots,x^{-1}X_{0},W''_0(X_0),\dots, W''_m(X_0) \mid X_0>x)\\
 \wto \law(Y_{-m},\dots,Y_{0}, Z_0,\dots, Z_m) \qquad \text{as $x \to \infty$.}
\end{multline*}
Using \eqref{eq:pseudoTail_inter1} and part (ii) of Lemma \ref{lemma:Slutsky4}, we get
\begin{multline*}
 \law (x^{-1}X_{-m},\dots,x^{-1}X_{0},X_0^{-1}X_{1},\dots,X_0^{-1}X_{m},W''_0(X_0),\dots, W''_m(X_0) \mid X_0>x)\\
 \wto \law(Y_{-m},\dots,Y_{0},\mu_A,\dots,\mu_A^m, Z_0,\dots, Z_m) \qquad \text{as $x \to \infty$.}
\end{multline*}
Hence, using $x^{-1}X_{i} = (X_0^{-1}X_{i}) (x^{-1}X_{0})$ and $Y_i = \mu_A^i Y_0$ for $i = 1, \dots, m$, together with Lemma \ref{lemma:cmt}, we obtain
\begin{multline}
\law (x^{-1}X_{-m},\dots,x^{-1}X_{m},W''_0(X_0),\dots, W''_m(X_0) \mid X_0>x)\\
 \wto \law(Y_{-m},\dots,Y_{m},Z_0,\dots,Z_m) \qquad \text{as $x \to \infty$.}
 \label{eq:pseudoTail_inter2}
\end{multline}
Next, we show that the above convergence holds with $W''_i(X_0)$ replaced by $W_i$, $i = 1, \dots, m$, which yields \eqref{eq:pseudoTail} (as explained at the beginning of the proof).
For this purpose, we prove that
\begin{align}\label{Delta_i}
 \cL(\Delta_0,\dots,\Delta_m\mid X_0>x) \wto \delta_{(0,\ldots,0)} \qquad \text{as $x\toi$,}
\end{align}
where, for $i\geq 0$,
\begin{align*}
\Delta_i:=W_i -  W''_i(X_0)  = \frac{1}{\sqrt{\mu_A^i X_0}}
\left(
\sum_{j=1}^{X_{i}} \tilde{A}_j^{(i+1)} -\sum_{j=1}^{\lfloor \mu_A^i X_0\rfloor} \tilde{A}_j^{(i+1)}  \right)
\end{align*}
 on the set $\{X_0>0\}$.
Note that $\Delta_0=0$.
It is enough to check that for all $\epsilon\in(0,1)$,
\begin{align}\label{Delta_i_proof}
\limsup_{x\toi} \pr(|\Delta_i|>\epsilon \mid X_0>x) \leq 2\epsilon \sigma_A^2, \qquad i = 0, \dots, m ,
\end{align}
 since then, by part (i) of Lemma \ref{lemma:Slutsky4}, \eqref{Delta_i} follows.
We will adapt the proof of R\'enyi's version of the Anscombe's theorem in Gut \cite[page 347]{gut:2005}.
For proving \eqref{Delta_i_proof}, let $i\in \{0,\dots,m\}$ be fixed.
For all $\vare\in(0,1)$ and $\ell\in\N$, let us introduce the notations
 \[
   n_0(\ell):=\lfloor\mu_A^i\ell\rfloor, \qquad
   n_1(\ell):=\lfloor(1-\vare^3)\mu_A^i\ell\rfloor + 1 , \qquad
   n_2(\ell):=\lfloor(1+\vare^3)\mu_A^i\ell\rfloor .
 \]
For all $\vare \in (0, 1)$ and $k\in\Nset$, with the notation $S_k:=\sum_{j=1}^k \tilde{A}_j^{(i+1)}$, we have
 \begin{align*}
   \pr(\vert \Delta_i\vert>\vare \mid X_0>x)
     & = \pr\big( \vert S_{X_i} - S_{n_0(X_0)}\vert > \vare \sqrt{\mu_A^i X_0} \mid X_0>x\big)\\
     & = \pr\big( \vert S_{X_i} - S_{n_0(X_0)}\vert > \vare \sqrt{\mu_A^i X_0}, X_i\in[n_1(X_0), n_2(X_0)] \mid X_0>x\big) \\
     &\phantom{=} + \pr\big( \vert S_{X_i} - S_{n_0(X_0)}\vert > \vare \sqrt{\mu_A^i X_0}, X_i\notin[n_1(X_0), n_2(X_0)] \mid X_0>x\big) \\
     & \leq \pr\Bigl(\max_{n_1(X_0)\leq k \leq n_2(X_0)}
     \vert S_k - S_{n_0(X_0)}\vert > \vare \sqrt{\mu_A^i X_0} \mid X_0>x\Bigr) \\
     &\phantom{\leq} + \pr(X_i\notin [n_1(X_0), n_2(X_0)] \mid X_0>x ).
 \end{align*}
Here, by Kolmogorov's theorem (see, e.g., Gut \cite[Theorem 3.1.6]{gut:2005}) and the independence of $X_0$ and $\tilde A_j^{(i+1)}$, $i \in \{0, \dots, m\}$, $j\in\NN$,
 we have for $x \geq 2/(\mu_A^i\vare^3)$,
 \begin{align}\label{help_Kolmogorov_ineq}
  \begin{split}
   &\pr\big( \max_{n_1(X_0)\leq k \leq n_2(X_0)}\vert S_k - S_{n_0(X_0)}\vert > \vare \sqrt{\mu_A^i X_0} \mid X_0>x\big) \\
     &= \sum_{\ell>x, \,\ell\in\Nset}\frac{1}{\pr(X_0>x)}
         \pr\Bigl(\max_{n_1(\ell)\leq k \leq n_2(\ell)}\vert S_k - S_{n_0(\ell)}\vert > \vare \sqrt{\mu_A^i\ell} , X_0=\ell\Bigr)\\
     &\leq \sum_{\ell>x,\, \ell\in\Nset}\frac{1}{\pr(X_0>x)}
           \frac{(n_2(\ell) - n_1(\ell))\sigma_A^2}{\vare^2\mu_A^i\ell} \pr(X_0=\ell)
      \leq 2 \vare\sigma_A^2.
  \end{split}
 \end{align}
Indeed, for $x \geq 2/(\mu_A^i\vare^3)$, we have $n_1(\ell)\leq n_0(\ell)\leq n_2(\ell)$ for $\ell\geq x$, $\ell\in\NN$,
 and
 \begin{align*}
   \vert S_k - S_{n_0(\ell)}\vert
      = \begin{cases}
          \left\vert \sum_{j=n_0(\ell)+1}^k \tilde A_j^{(i+1)} \right\vert & \text{if $n_0(\ell) \leq k\leq n_2(\ell)$,}\\
          \left\vert \sum_{j=k+1}^{n_0(\ell)} \tilde A_j^{(i+1)} \right\vert & \text{if $n_1(\ell) \leq k\leq n_0(\ell)-1$,}\\
        \end{cases}
 \end{align*}
 so using the independence of $\tilde A_j^{(i+1)}$, $j\in\NN$, by Kolmogorov's theorem, we have
 \begin{align*}
  & \pr\Bigl(\max_{n_1(\ell)\leq k \leq n_2(\ell)}\vert S_k - S_{n_0(\ell)}\vert > \vare \sqrt{\mu_A^i\ell} \Bigr) \\
  &  = \pr\left(  \max\left(  \max_{n_0(\ell)\leq k \leq n_2(\ell)}  \left\vert \sum_{j=n_0(\ell)+1}^k \tilde A_j^{(i+1)} \right\vert,
         \max_{n_1(\ell)\leq k \leq n_0(\ell)-1}  \left\vert \sum_{j=k+1}^{n_0(\ell)} \tilde A_j^{(i+1)} \right\vert  \right)  > \vare \sqrt{\mu_A^i\ell} \right)\\
  & =1-\left( 1 -  \pr\Bigl(\max_{n_0(\ell)\leq k \leq n_2(\ell)}  \left\vert \sum_{j=n_0(\ell)+1}^k \tilde A_j^{(i+1)} \right\vert > \vare \sqrt{\mu_A^i\ell}\Bigr)\right)\\
  &\phantom{=1-\;\;}
        \times \left( 1 -  \pr\Bigl(\max_{n_1(\ell)\leq k \leq n_0(\ell)-1}  \left\vert \sum_{j=k+1}^{n_0(\ell)} \tilde A_j^{(i+1)} \right\vert > \vare \sqrt{\mu_A^i\ell}\Bigr)\right)\\
  & \leq 1-\left( 1 - \frac{(n_2(\ell) - n_0(\ell))\sigma_A^2}{\vare^2\mu_A^i\ell} \right)
            \left( 1 - \frac{(n_0(\ell) - n_1(\ell))\sigma_A^2}{\vare^2\mu_A^i\ell} \right)\\
  & =  \frac{(n_2(\ell) - n_1(\ell))\sigma_A^2}{\vare^2\mu_A^i\ell}
        - \frac{(n_2(\ell) - n_0(\ell))(n_0(\ell) - n_1(\ell)) \sigma_A^2}{\vare^4\mu_A^{2i}\ell^2} \\
  & \leq \frac{(n_2(\ell) - n_1(\ell))\sigma_A^2}{\vare^2\mu_A^i\ell} ,
 \end{align*}
 yielding the first inequality in \eqref{help_Kolmogorov_ineq}.
The second inequality in \eqref{help_Kolmogorov_ineq} follows by
 \begin{align*}
  \frac{n_2(\ell) - n_1(\ell)}{\ell}
  \leq \frac{(1+\vare^3)\mu_A^i\ell - (1-\vare^3)\mu_A^i\ell}{\ell}
  = 2\mu_A^i \vare^3 .
 \end{align*}
Consequently, for all $\vare\in(0,1)$,
 \begin{align*}
   \pr(\vert \Delta_i\vert>\vare \mid X_0>x)
     \leq 2\vare\sigma_A^2 + \pr(X_i\notin [n_1(X_0), n_2(X_0)] \mid X_0>x ).
 \end{align*}
Here
 \begin{align*}
  &\pr(X_i\notin [n_1(X_0), n_2(X_0)] \mid X_0>x) \\
  &\leq \pr(X_i\notin [(1-\vare^3)\mu_A^iX_0+1, (1+\vare^3)\mu_A^iX_0-1] \mid X_0>x) \\
  &= \pr((\mu_A^iX_0)^{-1}X_i-1\notin [-\vare^3+(\mu_A^iX_0)^{-1}, \vare^3-(\mu_A^iX_0)^{-1}] \mid X_0>x) \\
  &\leq \pr((\mu_A^iX_0)^{-1}X_i-1\notin [-\vare^3+(\mu_A^ix)^{-1}, \vare^3-(\mu_A^ix)^{-1}] \mid X_0>x) \\
  &\leq \pr((\mu_A^iX_0)^{-1}X_i-1\notin [-\vare^3/2, \vare^3/2] \mid X_0>x)
 \end{align*}
 if $(\mu_A^ix)^{-1} \leq \vare^3/2$, i.e., for $x \geq 2/(\mu_A^i\vare^3)$, hence
 \begin{align*}
  \pr(X_i\notin [n_1(X_0), n_2(X_0)] \mid X_0>x)
   \to 0 \qquad \text{as $x\to\infty$,}
 \end{align*}
 due to \eqref{eq:pseudoTail_inter0}.
Consequently, we have \eqref{Delta_i_proof}.
By \eqref{Delta_i} and part (i) of Lemma \ref{lemma:Slutsky4}, we obtain
 \[
   \law(0,\dots,0,\Delta_0,\dots,\Delta_m \mid X_0>x)\\
 \wto \delta_{(0,\dots,0,0,\dots,0)} \qquad \text{as $x \to \infty$.}
 \]
Hence, using \eqref{eq:pseudoTail_inter2} and Lemma \ref{lemma:Slutsky3}, we obtain \eqref{eq:pseudoTail} with $W'_i$ replaced by $W_i$, $i \in \{0,\ldots,m\}$, as desired.
\proofend

\section{Point process convergence}\label{sec:Ppc}

Let $(a_n)_{n\in\NN}$ be a sequence of positive real numbers satisfying $\lim_{n\to\infty} a_n=\infty$ and
 \begin{align}\label{a_n_quantile}
   \lim_{n\to\infty} n \pr( X_0 > a_n ) = 1 .
 \end{align}
Note that for $n\geq 2$, one can choose $a_n$ to be the maximum of 1 and the $1-\frac{1}{n}$ lower quantile of $X_0$.
In fact, $a_n = n^{1/\alpha} L(n)$, $n \in \NN$, for some slowly varying continuous function $L: (0, \infty) \to (0, \infty)$, see, e.g., Araujo and Gin\'e \cite[Exercise 6 on page 90]{araujo:gine:1980}.
In what follows we fix such a sequence $(a_n)_{n\in\NN}$ satisfying \eqref{a_n_quantile}.
By the proof of Lemma 3.2 in Basrak et al.~\cite{basrak:kulik:palmowski:2013}, for any sequence of positive integers $(r_n)_{n\geq 1}$
 such that $\lim_{n\toi}r_n=\infty$ and $\lim_{n\toi}n/r_n=\infty$, we have
 \begin{align}\label{eq:AC}
 \lim_{m\toi}\limsup_{n\toi}\pr\left(\max_{m< |i| \leq r_n} X_{i}>a_nu\; \big| \; X_0>a_nu\right)=0 \qquad \text{for all $u>0$.}
 \end{align}
Condition \eqref{eq:AC} is sometimes called the anticlustering condition or finite mean cluster size condition for the strongly stationary process $(X_i)_{i\in\ZZ}$.
Later, in the proof of Theorem \ref{th:PPconv} we will choose a particular sequence $(r_n)_{n\in\NN}$ satisfying a kind of mixing condition \eqref{eq:A(a_n)} as well.
For a possible choice of such a sequence $(r_n)_{n\in\NN}$, see the proof of Lemma \ref{lemma:mixing}.

Consider the space $S:=(0,\infty)\times \R$ with the usual topology and call a Borel set $B \subset S$ bounded if it is separated from the vertical line $\{(0, y) : y \in \RR\}$, i.e., there exists $\epsilon>0$ such that $B \subset \{(x,y)\in(0,\infty)\times \R : x>\epsilon\}$.
Using the terminology of Basrak and Planini\'c~\cite{basrak:planinic:2019}, the collection of bounded sets is a boundedness (Borel subsets of bounded sets are bounded, and finite union of bounded sets are bounded), which properly localizes $S$.
Denote by $\cB(S)$ the Borel $\sigma$-algebra on $S$ and by $\hC_S$ the class of bounded, continuous functions $f : S \to [0, \infty)$ with
bounded support.
Hence, if $f \in \hC_S$, then there exists an $\vare>0$ such that $f(x,y)=0$ for all $(x,y)\in S$ with $x\leq \vare$.

Let $\mathcal{M}_p (S)$ be the space of integer-valued measures (or point measures) on $S$ which are finite on bounded Borel sets
 (called locally finite measures), topologized by the so--called vague topology.
The associated notion of vague convergence of $\mu_n\in\mathcal{M}_p(S)$ towards $\mu\in\mathcal{M}_p(S)$ as $n\to\infty$, denoted by $\mu_n \vto \mu$ as $n\to\infty$, is defined by the condition $\mu_n(f) \to \mu(f)$ as $n\to\infty$ for all $f \in \hC_S$, where $\kappa(f) := \int_S f(x,y) \, \kappa(\dd x, \dd y)$ for each $\kappa \in \mathcal{M}_p(S)$.
Alternatively, we could have used the framework of so--called $\mathbb{M}_{\mathbb{O}}$--convergence from Lindskog et al.~\cite{lindskog:resnick:roy:2014}, see also Basrak and Planini\'c~\cite{basrak:planinic:2019} for the discussion on vague convergence.
Convergence in distribution of random measures with respect to the vague topology will be denoted by $\dvague$.

In what follows, we use the theory of Kallenberg \cite[Chapter 4]{kallenberg:2017}.
To do so, one can equip $S$ with the metric $d:S\times S\to [0,\infty)$,
 \begin{align}\label{S_metric}
   d( (x,y),(x',y') ):= \min\big(\sqrt{(x-x')^2 + (y-y')^2},1\big)+ \left\vert \frac{1}{x} - \frac{1}{x'}\right\vert
 \end{align}
 for $(x,y),(x',y')\in S$.
This metric is complete and induces the original topology of $S$, and the family of $d$--bounded sets is precisely the family of bounded sets in $S$ as defined above, see Lemma \ref{Lem_S_top}.
In Kallenberg \cite[page 125]{kallenberg:2017}, one can find a similar metric.
Note that $\mathcal{M}_p(S)$ is a complete separable metric space with a metric inducing the vague topology, see Kallenberg \cite[Lemma 4.6]{kallenberg:2017}.

The following well--known result (stated in our setting) describes vague convergence of point measures and is crucial for the use of continuous mapping arguments (for a proof, based on Lemma 3.13 in Resnick \cite{resnick:1987}, see Appendix \ref{app:vague_convergence_of_point_measures}), we will use it
 in the proof of Theorem \ref{thm:partial_sums}.

\begin{lemma}\label{lem:vagueConv_pointMeas}
Let $\mu,\mu_n\in \mathcal{M}_p (S)$, $n\in\NN$.
Then $\mu_n\vto \mu$ as $n\to\infty$ if and only if for each $\epsilon>0$ satisfying $\mu(\{\epsilon\}\times \R)=0$ there exist integers $n_0,M\geq 0$ and a labeling of the points of $\mu$ and $\mu_n$, $n \geq n_0$, in $(\epsilon,\infty)\times \R$ such that
\begin{align*}
\mu_n|_{(\epsilon,\infty)\times \R}=\sum_{i=1}^M \delta_{(x_i^{(n)}, \, y_i^{(n)})} \, , \qquad  \mu|_{(\epsilon,\infty)\times \R}=\sum_{i=1}^M \delta_{(x_i, y_i)} \, ,
\end{align*}
and $x_i^{(n)}\to x_i$ and $y_i^{(n)}\to y_i$ as $n \toi$ for all $i=1,\dots,M$, where $\mu|_{B}$ denotes the restriction of $\mu$ to the set $B\subset S$.
\end{lemma}

Recall \eqref{tail_process} and define
 \begin{align}\label{theta}
  \theta := \pr (\sup_{j<0} Y_j \leq 1) = \pr(K = 0) = 1-\mu_A^\alpha.
 \end{align}
Indeed, $\{\sup_{j<0} Y_j \leq 1\} = \{K = 0 \}$, since if $K=0$, then, by \eqref{tail_process}, $Y_j=0$ for all $j<0$
 yielding $\sup_{j<0} Y_j =0 \leq 1$; and if $K\geq 1$, then $\sup_{j<0} Y_j \geq Y_{-1} = \mu_A^{-1} Y_0 >1$ almost surely
  due to the facts that $\pr(Y_0\geq 1)=1$ and $\mu_A\in(0,1)$.

\begin{theorem} \label{th:PPconv}
In $\mathcal{M}_p(S)$,
\begin{align} \label{eq:PPconv}
 N_n := \sideset{}{^{*}} \sum_{j=1}^n   \delta_{\bigl( \frac{X_{j}}{a_n}, \frac{M_{j+1}}{\sqrt{X_j}} \bigr)}
      = \sum_{ \{ j\in\{1,\ldots,n\} \,:\, X_j>0 \} } \delta_{\bigl(\frac{X_j}{a_n},\,\frac{M_{j+1}}{\sqrt{X_j}}\bigr)}
 \dvague
 N := \dsum_{i=1}^\infty \dsum_{j=0}^\infty  \delta_{(P_i\mu_A^j,Z_{i,j})}
\end{align}
as $n \toi$, where $\sum_{i=1}^\infty \delta_{P_i}$ denotes a Poisson point process on $(0,\infty)$
 with intensity $\theta \, \dd(-y^{-\alpha})$ which is independent of  an i.i.d.\ array
 of $\cN(0,\sigma_A^2)$-distributed random variables $\{Z_{i,j}: i\in\NN, j\geq 0\}$.
\end{theorem}

\noindent{\bf Proof.}
First, we check that $N$ is a locally finite measure almost surely, i.e., $N(B)$ is finite almost surely for each bounded Borel measurable subset $B$ of $S$.
Assume that $B\in\cB(S)$ is bounded.
Then we have $B \subset (\vare,\infty)\times\RR$ for some $\vare>0$, and
 \[
   N((\vare, \infty) \times \R)
   = \sum_{i=1}^\infty \sum_{j=0}^\infty \1{\{(P_i\mu_A^j,Z_{i,j})\in(\vare,\infty)\times\R\}}
   = \sum_{j=0}^\infty \sum_{i=1}^\infty \1{\{P_i > \mu_A^{-j} \vare\}}
   < \infty
 \]
 almost surely, since, for each $j \geq 0$, there are only finitely many $P_i$'s greater than $\mu_A^{-j} \vare$ almost surely, and for sufficiently large $j \geq 0$, the set $\{i \in \N : P_i > \mu_A^{-j} \vare\}$ is empty, since $\mu_A^{-j} \vare \geq \max_{i\in\N} P_i$ holds for sufficiently large $j \geq 0$ due to $\mu_A^{-j} \to \infty$ as $j \to \infty$.

By Remark 3.1 in Basrak et al.~\cite{basrak:kulik:palmowski:2013} or Lemma F.1 in Barczy et al.~\cite{barczy:nedenyi:pap:2019}, the strongly stationary process $(X_j)_{j\in \Z}$ is strongly mixing, and since $ M_{j+1}/\sqrt{X_j} \in \sigma(X_j,X_{j+1})$ for every $j\geq 0$ (understanding this ratio 0 when $X_j=0$), a modification of Lemma 2.3.9 in Basrak \cite{basrakPhD} or the proof of Proposition 1.34 in Krizmani\'c \cite{Kri}, shows that there exists a sequence of positive integers $(r_n)_{n\in \N}$ satisfying $r_n\toi$ and $r_n/n\to 0$ as $n\to\infty$, such that
\begin{align}\label{eq:A(a_n)}
   \ex\biggl[\exp\biggl\{- \sideset{}{^{*}} \sum_{j=1}^n f\biggl(\frac{X_j}{a_n},  \frac{M_{j+1}}{\sqrt{X_j}} \biggr)\biggr\}\biggr]
   - \biggl(
     \ex\biggl[\exp\biggl\{- \sideset{}{^{*}} \sum_{j=1}^{r_n} f\biggl(\frac{X_j}{a_n}, \frac{M_{j+1}}{\sqrt{X_j}}  \biggr)\biggr\}\biggr]\biggr)^{k_n}
    \to 0
\end{align}
as $n\toi$ for all $f \in \hC_S$, where $k_n:=\lfloor n/r_n\rfloor$, see Lemma \ref{lemma:mixing}.
This is similar to condition (2.1) in Davis and Hsing \cite{davis:hsing:1995}.
Note that for this sequence $(r_n)_{n\in\NN}$, the anticlustering condition \eqref{eq:AC} holds automatically, as explained before \eqref{eq:AC}.

By Theorem 4.11 in Kallenberg \cite{kallenberg:2017} and \eqref{eq:A(a_n)}, $N_n\dvague N$ as $n\to\infty$ if and only if
 $\widetilde N_n:=\sum_{i=1}^{k_n} \widetilde N_{n,i}\dvague N$ as $n\to\infty$, where for each $n\in \N$,
 $\widetilde N_{n,i}$, $i=1,\dots,k_n$, are i.i.d.\ point processes on  $S$ with common distribution
 equal to the distribution of
 \[
  N_{r_n} = \sideset{}{^{*}} \sum_{j=1}^{r_n} \delta_{\left(\frac{X_j}{a_n}, \frac{M_{j+1}}{\sqrt{X_j}}\right)} \, ,
 \]
 since \eqref{eq:A(a_n)} takes the form
 \[
 \ex[\ee^{-N_n(f)}] - \ex[\ee^{-\widetilde N_n(f)}]  \to 0
 \]
 as $n\toi$ for all $f \in \hC_S$.

We will apply Theorem 4.22 in Kallenberg \cite{kallenberg:2017}.
The array of point processes $\{ \widetilde N_{n,i}, n\in \N, i=1,\dots, k_n\}$ forms a null array,
 since if $B \subset S$ is a bounded Borel set, then
 there exists an $\vare>0$ such that $B \subset(\vare,\infty)\times\R$, and then
 \begin{align*}
  &\sup_{ i\in\{1,\ldots, k_n\} }
      \pr( \widetilde N_{n,i}(B)>0) = \pr( \widetilde N_{n,1}(B)>0) = \pr\left(\sideset{}{^{*}} \sum_{j=1}^{r_n} \delta_{(X_j/a_n, M_{j+1}/\sqrt{X_j} )}(B)>0\right)\\
  &\leq \pr\left( \bigcup_{j=1}^{r_n} \left\{ \left(\frac{X_j}{a_n},  \frac{M_{j+1}}{\sqrt{X_j}}  \right)\in(\vare,\infty)\times \R \right\} \right)
        \leq \sum_{j=1}^{r_n} \pr\left( \frac{X_j}{a_n} \in (\vare,\infty)\right)\\
  &= r_n \pr\left( \frac{X_1}{a_n} \in (\vare,\infty)\right)
   = \frac{r_n}{n}\cdot n \pr(X_1>a_n) \cdot \frac{\pr(X_1> \vare a_n)}{\pr(X_1 > a_n)}
   \to 0\cdot 1 \cdot \vare^{-\alpha} = 0 \qquad \text{as \ $n\to\infty$,}
 \end{align*}
 where the last step follows by \eqref{a_n_quantile} and the fact that $X_1$ is regularly varying with tail index $\alpha$.
By Kallenberg \cite[Theorem 4.22 and page 89]{kallenberg:2017}, $\widetilde N_n \dvague N$ as $n\to\infty$ if
 there exists a measure $\nu$ on $\mathcal{M}_p(S)\setminus\{0\}$ (furnished with
 the smallest $\sigma$-algebra making all the evaluation maps $\mathcal{M}_p(S)\ni m \mapsto m(F)$ measurable, where $F\in \cB(S)$) such that
 \begin{itemize}
 \item[(i)]
  $\int_{\mathcal{M}_p(S)\setminus\{0\}} \min(\kappa(B),1)\, \nu(\dd\kappa)<\infty$ for all bounded Borel subsets of $S$,
 \item[(ii)]
  $k_n\ex[1-\ee^{-N_{r_n}(f)}]\to \int_{\mathcal{M}_p(S)\setminus\{0\}} (1-\ee^{-\kappa(f)})\, \nu(\dd\kappa)$ as \ $n\to\infty$
        for all $f \in \hC_S$,
 \item[(iii)]
  $-\log\ex[\ee^{-N(f)}]= \int_{\mathcal{M}_p(S)\setminus\{0\}} (1-\ee^{-\kappa(f)})\, \nu(\dd\kappa)$
   for all $f \in \hC_S$.
 \end{itemize}
Here we note that, with the notations of Kallenberg \cite{kallenberg:2017}, the non-random locally finite measure $\alpha$ on $S$
 appearing in Theorem 4.22 in Kallenberg \cite{kallenberg:2017} is the null measure, by Theorem 3.20 in Kallenberg \cite{kallenberg:2017}.

Let $\nu$ be a measure on $\mathcal{M}_p(S)\setminus\{0\}$ given by
\begin{align*}
\nu(\,\cdot\,):= \theta \int_0^\infty \pr\left(
\dsum_{j=0}^{\infty} \delta_{(y\mu_A^j,Z_j)}
 \in \cdot\right) \alpha y^{-\alpha - 1}\, \dd y \, ,
\end{align*}
 where $(Z_j)_{j\geq 0}$ is an i.i.d.\ sequence of $\cN(0,\sigma_A^2)$--distributed random variables
 being independent of $Y_0$ and $K$ with $\sigma^2_A =\Var(A)$.
Note that we have
 \begin{equation}\label{integral}
   \int_{\mathcal{M}_p(S)\setminus\{0\}} h(\kappa(g)) \, \nu(\dd\kappa)
   = \theta \int_0^\infty \ex\left[h\left(\sum_{j=0}^{\infty} g(y\mu_A^j, Z_j) \right)\right] \alpha y^{-\alpha - 1} \, \dd y
 \end{equation}
 for any measurable functions $h : \R \to [0, \infty)$ and $g : S \to [0, \infty)$ such that $g$ has bounded support,
 since this trivially holds for linear combinations of indicator functions, which extends by monotone convergence to arbitrary nonnegative measurable functions.
Consequently, $\nu$ satisfies property (i) above, since if $B$ is a bounded Borel subset of $S$,
 then there exists $\vare>0$ such that $x>\vare$ for all $(x,z)\in B$, and hence
 \begin{align*}
  &\int_{\mathcal{M}_p(S)\setminus\{0\}} \min(\kappa(B),1)\, \nu(\dd\kappa)
       = \theta \int_0^\infty \ex\left[ \left(\sum_{j=0}^{\infty} \1{B}(y\mu_A^j, Z_j) \right) \wedge 1  \right] \alpha y^{-\alpha - 1} \, \dd y \\
      & = \theta \int_{y>\vare}^\infty \ex\left[ \left(\sum_{j=0}^{\infty} \1B(y\mu_A^j, Z_j) \right) \wedge 1  \right] \alpha y^{-\alpha - 1} \, \dd y
       \leq \theta \int_{y>\vare}^\infty \alpha y^{-\alpha - 1} \, \dd y
       = \theta \vare^{-\alpha}
       <\infty,
 \end{align*}
 where we used that if $y\in(0,\vare]$, then $y\mu_A^j\in(0,\vare]$, $j\geq 0$.

Next we turn to prove (iii).
We check that $N$ given in (\ref{eq:PPconv}) satisfies
\begin{equation}\label{N2_check}
  -\log\ex[\ee^{-N(f)}]= \theta \int_0^\infty \ex[1-\ee^{ -\sum_{j=0}^\infty f({y\mu_A^j, Z_j) } }]
                          \alpha y^{-\alpha - 1} \,\dd y
                       = \int_{\mathcal{M}_p(S)\setminus\{0\} } (1-\ee^{-\kappa(f)}) \, \nu(\dd\kappa)
\end{equation}
 for all $f \in \hC_S$.
The second equality follows from \eqref{integral}.
For the first equality, we provide two alternative proofs.
The first equality in \eqref{N2_check} can be derived using the representation
\[
N(f)=\sum_{i=1}^\infty\sum_{j=0}^{\infty} f(P_i \mu_A^{j}, Z_{i,j})=\sum_{i=1}^\infty f'(P_i, (Z_{i,j})_{j\geq 0})=N'(f') ,
\]
 where $f'(p, (z_j)_{j\geq 0}):=\sum_{j=0}^{\infty}f(p \mu_A^{j}, z_j)$ for $p \in (0, \infty)$, $z_j \in \R$, $j \geq 0$, and $N':=\sum_{i=1}^{\infty} \delta_{(P_i, \, (Z_{i,j})_{j\geq 0})}$ is a Poisson point process on $(0,\infty)\times \R^{\NN \cup \{0\}}$ with intensity $\theta \, \dd(-y^{-\alpha}) \times \cL(Z_{1,1})^{\NN \cup \{0\}}$ (see, e.g., Kingman \cite[Section 5.2]{kingman:1993}).
Indeed, $f \in \hC_S$ implies that the sum $\sum_{j=0}^{\infty}f(p \mu_A^{j}, z_j)$ has only finitely many non-zero terms, and
 one can use the expression for the Laplace functional of $N'$ together with Fubini's theorem.
The second proof of the first equality in \eqref{N2_check} is based on the representation
 \begin{align*}
  \ex[\ee^{-N(f)}] &= \ex\Bigl[\ex\bigl[\ee^{-\sum_{i=1}^\infty\sum_{j=0}^\infty f(P_i\mu_A^j, Z_{i,j})} \,\big|\, (P_i)_{i\in\NN}\bigr]\Bigr] \\
  &= \ex\Biggl[\prod_{i=1}^\infty \prod_{j=0}^\infty
               \ex\bigl[\ee^{-f(P_i\mu_A^j,Z)} \,\big|\, (P_i)_{i\in\NN}\bigr]\Biggr]
   = \ex\Biggl[\prod_{i=1}^\infty h(P_i)\Biggr] ,
 \end{align*}
 where $Z$ is an $\cN(0, \sigma_A^2)$-distributed random variable, independent of $(P_i)_{i\in\NN}$, and $h : (0, \infty) \to \R$ is defined by
 \[
   h(x) := \prod_{j=0}^\infty \ex\bigl[\ee^{-f(x\mu_A^j,Z)}\bigr] , \qquad x \in (0, \infty) .
 \]
Consequently, using the Laplace functional of the Poisson point process $\sum_{i=1}^\infty \delta_{P_i}$, we obtain
 \begin{align*}
  \ex[\ee^{-N(f)}] &= \ex\bigl[\ee^{-\sum_{i=1}^\infty (-\log(h(P_i)))}\bigr]
   = \exp\biggl\{- \int_0^\infty \bigl(1 - \ee^{-(-\log(h(x)))}\bigr) \theta \alpha x^{-\alpha-1} \, \dd x\biggr\} \\
  &= \exp\biggl\{- \theta \int_0^\infty (1 - h(x) ) \alpha x^{-\alpha-1} \, \dd x\biggr\} \\
  &= \exp\biggl\{- \theta \int_0^\infty \biggl(1 - \prod_{j=0}^\infty \ex\bigl[\ee^{-f(x\mu_A^j,Z)}\bigr]\biggr) \alpha x^{-\alpha-1} \, \dd x\biggr\} \\
  &= \exp\biggl\{- \theta \int_0^\infty \biggl(1 - \ex\bigl[\ee^{-\sum_{j=0}^\infty f(x\mu_A^j,Z_j)}\bigr]\biggr) \alpha x^{-\alpha-1} \, \dd x\biggr\} ,
 \end{align*}
 as desired.

Now we prove (ii).
By \eqref{N2_check}, we need to show that
\begin{align}\label{eq:intensity}
k_n\ex[1-\ee^{-N_{r_n}(f)}]\to \theta \int_0^\infty \ex[1-\ee^{ -\sum_{j=0}^{\infty} f(y\mu_A^j, Z_j)}]
  \alpha y^{-\alpha - 1} \dd y
\end{align}
 as $n\to\infty$ for all $f \in \hC_S$.
Take an arbitrary $f \in \hC_S$ and let $\epsilon>0$ be such that $x\leq \epsilon$ implies that $f(x,y)=0$ for all $y\in \R$.
To show (\ref{eq:intensity}) we rely on the condition (\ref{eq:AC}) and use similar arguments as in the proof of Basrak and Segers \cite[Theorem 4.3]{basrak:segers:2009}.
Note that, even though condition (\ref{eq:AC}) concerns only the process $(X_j)_{j\in \Z}$, it will be sufficient for (\ref{eq:intensity}), since the bounded sets in $S$ are separated from the vertical line (see Lemma \ref{Lem_S_top}) by our choice of bounded sets in $S$.

For $n \in \N$ and $k, \ell \in \Z$ with $k \leq \ell$, write $M_{k,\ell} := \max_{k\leq j \leq \ell} X_j$ (if $\ell< k$ set $M_{k,\ell} := 0$),
 and
 \[
  c_n(k,\ell) := 1-\exp\left\{- \sideset{}{^{*}} \sum_{j=k}^\ell  f\left( \frac{X_j}{a_n}, \frac{M_{j+1}}{\sqrt{X_j}} \right)\right\}.
 \]
In particular, $c_n(1,r_n)=1-\exp\{- N_{r_n}(f)\}$, and since  $f(x,y)=0$ if $(x,y)\in (0,\epsilon]\times\RR$, we have
 $\{M_{1,r_n}\leq a_n\vare\} \subset \{c_n(1,r_n) =0\}$, and then
\begin{align}\label{help14}
\ex[c_n(1,r_n)] = \ex[c_n(1,r_n)\1{\{M_{1,r_n}>a_n\epsilon\}}]
= \sum_{i=1}^{r_n} \ex[c_n(1,r_n)\1{\{M_{1,i-1}\leq a_n\epsilon < X_i\}}] \, .
\end{align}
Since $k_n\sim (r_n \pr(X_0>a_n))^{-1}$ as $n\to\infty$, (\ref{eq:intensity}) will follow if we show that
\begin{align*}
\frac{\ex[c_n(1,r_n)]}{r_n \pr(X_0>a_n)} \to \theta \int_0^\infty \ex[1-\ee^{ -\sum_{j=0}^{\infty} f(y \mu_A^j, Z_j)}]
  \alpha y^{-\alpha - 1} \dd y \qquad \text{as $n\toi$.}
\end{align*}
Fix now $m\in \N$ and assume that $n$ is large enough so that $r_n\geq 2m +1$. The key observation now is that, since  $f(x,y)=0$ if $x\leq \epsilon$,  for all $m+1\leq i \leq r_n - m$, we have $c_n(1,r_n)=c_n(i-m,i+m)$ if $M_{1, i - m -1} \vee M_{i+m+1, r_n} \leq a_n \epsilon$.
Using the strong stationarity of $(X_i)_{i\in\ZZ}$ and that $c_n(k,l)\in[0,1]$, $n\in\NN$, $k,l\in\ZZ$,
 this implies that for all $i\in\{m+1,\ldots,r_n-m\}$,
\begin{align*}
 &\Big|\ex[c_n(1,r_n)\1{\{M_{1,i-1}\leq a_n\epsilon < X_i\}}] - \ex[c_n(i-m,i+m)\1{\{M_{i-m,i-1}\leq a_n\epsilon < X_i\}}] \Big|  \\
 & \leq \Big|\ex\Big[ \Big( c_n(1,r_n)\1{\{M_{1,i-1}\leq a_n\epsilon < X_i\}} - c_n(i-m,i+m)\1{\{M_{i-m,i-1}\leq a_n\epsilon < X_i\} } \Big)
              \1{\{ M_{1, i-m-1} \vee M_{i+m+1,r_n}  \leq a_n\epsilon \} }  \Big] \Big| \\
 & +   \Big|\ex\Big[ \Big( c_n(1,r_n)\1{\{M_{1,i-1}\leq a_n\epsilon < X_i\}} - c_n(i-m,i+m)\1{\{M_{i-m,i-1}\leq a_n\epsilon < X_i\} } \Big)
              \1{\{ M_{1, i-m-1} \vee M_{i+m+1,r_n}  > a_n\epsilon \}}  \Big] \Big| \\
 &=  \Big|\ex \Big[ {\big( c_n(1,r_n) - c_n(i-m,i+m)\big) \1{\{M_{1,i-1}\leq a_n\epsilon < X_i\}} }
              \1{\{ M_{1, i-m-1} \vee M_{i+m+1,r_n}  \leq a_n\epsilon \} }  \Big] \Big|\\
 & + \Big |\ex\Big[ \Big( c_n(1,r_n)\1{\{M_{1,i-1}\leq a_n\epsilon \}} - c_n(i-m,i+m)\1{\{M_{i-m,i-1}\leq a_n\epsilon \} } \Big)
              \1{\{X_i>a_n\vare\}}\1{\{ M_{1, i-m-1} \vee M_{i+m+1,r_n}  > a_n\epsilon \} }  \Big] \Big|\\
 &\leq \pr (M_{1, i-m-1} \vee M_{i+m+1,r_n}>a_n \epsilon , X_i>a_n\epsilon)\\
 &= \pr (M_{1-i, -m-1} \vee M_{m+1,r_n-i}>a_n \epsilon , X_0>a_n\epsilon)\\
 &\leq \pr (M_{-r_n, -m-1} \vee M_{m+1,r_n}>a_n \epsilon , X_0>a_n\epsilon) ,
\end{align*}
 where the last step follows by $-r_n\leq 1-i$ and $r_n-i\leq r_n$.
For $i\in\{1,\ldots,m\}\cup\{r_n-m+1,\ldots,r_n\}$,
 use the trivial bound
 \begin{align*}
  &\big|\ex[c_n(1,r_n)\1{\{M_{1,i-1}\leq a_n\epsilon < X_i\}}] - \ex[c_n(i-m,i+m)\1{\{M_{i-m,i-1}\leq a_n\epsilon < X_i\}}]\big| \\
  &\leq \ex[|c_n(1,r_n)\1{\{M_{1,i-1}\leq a_n\epsilon\}}
             - c_n(i-m,i+m)\1{\{M_{i-m,i-1}\leq a_n\epsilon\}}| \1{\{ a_n\epsilon < X_i\}} ]
   \leq \pr(X_0>a_n \epsilon) .
 \end{align*}
Consequently, using again the strong stationarity of $(X_i)_{i\in\ZZ}$ and \eqref{help14}, we have
 \begin{align*}
 &\Delta_{n,m}:= \left|\frac{\ex[c_n(1,r_n)]}{r_n \pr(X_0>a_n)}
   - \frac{\ex[c_n(-m,m)\1{\{M_{-m,-1}\leq a_n\epsilon < X_0\}}]}{\pr(X_0>a_n)} \right| \\
 &=\frac{1}{r_n\pr(X_0>a_n)} \left\vert
    \sum_{i=1}^{r_n} \ex[c_n(1,r_n) \1{\{M_{1,i-1}\leq a_n\vare < X_i\}}]
                      - r_n \ex[c_n(-m,m) \1{\{M_{-m,-1}\leq a_n\vare < X_0\}}]\right\vert\\
 &=\frac{1}{r_n\pr(X_0>a_n)} \left\vert
    \sum_{i=1}^{r_n} \ex[c_n(1,r_n) \1{\{M_{1,i-1}\leq a_n\vare < X_i\}}]
                      - \sum_{i=1}^{r_n} \ex[c_n(i-m,i+m) \1{\{M_{i-m,i-1}\leq a_n\vare < X_i\}}]\right\vert
 \end{align*}
 \begin{align*}
 &\leq  \frac{1}{r_n\pr(X_0>a_n)}
  \Bigg( \sum_{i\in\{1,\ldots,m\}\cup \{r_n-m+1,\ldots,r_n\}} \pr(X_0>a_n\vare) \\
 &\phantom{\leq  \frac{1}{r_n\pr(X_0>a_n)} \Big(}
           + \sum_{i=m+1}^{r_n-m} \pr(M_{-r_n, -m-1} \vee M_{m+1,r_n}>a_n \epsilon, X_0 > a_n\vare) \Bigg)\\
 & = \frac{1}{r_n\pr(X_0>a_n)}
     \Big( 2m \pr(X_0>a_n\vare)  +  (r_n - 2 m) \pr(M_{-r_n, -m-1} \vee M_{m+1,r_n}>a_n \epsilon, X_0 > a_n\vare) \Big)\\
 & = \frac{\pr(X_0>a_n\vare)}{\pr(X_0>a_n)}
     \left( \frac{2m}{r_n}   +  \left( 1 -\frac{2m}{r_n}\right) \pr(M_{-r_n, -m-1} \vee M_{m+1,r_n}>a_n \epsilon \mid  X_0 > a_n\vare) \right).
 \end{align*}
Note $\tfrac{\pr(X_0>a_n\epsilon)}{\pr(X_0>a_n)} \sim \epsilon^{-\alpha}$  and $r_n\toi$ as $n\toi$, hence, (\ref{eq:AC}) implies that
 \[
  \lim_{m\toi} \limsup_{n\toi} \Delta_{n,m}=0.
 \]
Observe that on the event $M_{-m,-1}\leq a_n \epsilon$, we have $c_n(-m,m)=c_n(0,m)$, $n\in\NN$.
Hence,
\begin{equation}\label{help11_point}
 \frac{\ex[c_n(-m,m)\1{\{M_{-m,-1}\leq a_n\epsilon < X_0\}}]}{\pr(X_0>a_n)}
 =\frac{\pr(X_0>a_n\epsilon)}{\pr(X_0>a_n)}
  \,\ex[c_n(0,m) \1{\{M_{-m,-1}\leq a_n\epsilon\}} \mid X_0>a_n \epsilon].
\end{equation}
The first term on the right hand side of \eqref{help11_point} tends to $\epsilon^{-\alpha}$ as $n\toi$.
For each $m \in \NN$, we have
 \[
   c_n(0,m)=1-\exp\biggl\{- \sideset{}{^{*}} \sum_{j=0}^m f\biggl(\frac{X_j}{a_n}, \frac{M_{j+1}}{\sqrt{X_j}} \biggr)\biggr\}
   = 1-\exp\biggl\{- \sideset{}{^{*}} \sum_{j=0}^m f\biggl(\frac{X_j}{a_n}, W_j'\biggr)\biggr\} ,
 \]
 hence, by Proposition \ref{prop:pseudoTail} and a consequence of the conditional continuous mapping theorem
 (see part (ii) of Lemma \ref{lemma:cmt}) with a bounded Borel measurable function $\widetilde h : \RR^{3m+2} \to \RR$ satisfying
 \[
   \widetilde h(x_{-m},\ldots,x_m,w_0',\ldots,w_m')
   = \biggl(1- \exp\biggl\{-\sum_{j=0}^mf(\epsilon x_j, w_j')\1{(0,\infty)}(x_j)\biggr\}\biggr)\1{\{\max\{x_{-m},\ldots,x_{-1}\}\leq1\}}
 \]
 for $x_{-m},\ldots,x_m \in[0,\infty)$, $w_0',\ldots,w_m' \in \RR$, we obtain
 \begin{align*}
  &\ex(c_n(0,m) \1{\{M_{-m,-1}\leq a_n\epsilon\}} \mid X_0>a_n\epsilon) \\
  &\to \ex\biggl(\biggl(1- \exp\biggl\{-\sum_{j=0}^mf(\epsilon \mu_{A}^j Y_0, Z_j)\biggr\}\biggr)\1{\{\max\{Y_{-m},\ldots,Y_{-1}\}\leq1\}}\biggr) \qquad \text{as $n \toi$,}
 \end{align*}
 since the absolute continuity of $Y_0$ and the independence of $Y_0$ and $K$ imply
 \begin{align*}
  \pr((Y_{-m},\ldots,Y_m,Z_0,\ldots,Z_m)\in D_{\widetilde h})
  &=\pr(\max\{Y_{-m},\ldots,Y_{-1}\}=1) \\
  &= \pr(Y_0\max\{\mu_A^{-m}\1{\{K\geq m\}},\ldots,\mu_A^{-1}\1{\{K\geq1\}}\}=1)
   = 0 ,
 \end{align*}
 where $D_{\widetilde h}$ denotes the set of discontinuities of $\widetilde h$.
Consequently, by the dominated convergence theorem, the second term on the right hand side of \eqref{help11_point} as $n\toi$ and then as $m\toi$, converges to
\begin{align*}
 &\ex\Big[\Big(1- \exp\Big\{-\sum_{j=0}^{\infty}f(\epsilon \mu_{A}^j Y_0, Z_j)\Big\}\Big)\1{\{\sup_{j<0}Y_j \leq 1\}}\Big]\\
 &=\ex\Big[\Big(1- \exp\Big\{-\sum_{j=0}^{\infty}f(\epsilon \mu_{A}^j Y_0, Z_j)\Big\}\Big)\1{\{K=0\}}\Big]
  =\theta \ex\Big[1- \exp\Big\{-\sum_{j=0}^{\infty}f(\epsilon \mu_{A}^j Y_0, Z_j)\Big\}\Big] \, ,
\end{align*}
 where the last step follows by the fact that $K$ is independent of $Y_0$ and $(Z_j)_{j\geq 0}$ and since $\pr(K=0)=\theta$ (see \eqref{theta}).
Altogether, using the fact that $Y_0$ is Pareto distributed such that $\pr(Y_0 \geq y) = y^{-\alpha}$, $y\geq 1$,
 and independent of $(Z_j)_{j\geq 0}$,
\begin{align*}
\lim_{n\toi} \frac{\ex[c_n(1,r_n)]}{r_n \pr(X_0>a_n)} & = \lim_{m\toi} \lim_{n\toi} \frac{\ex[c_n(-m,m)\1{\{M_{-m,-1}\leq a_n\epsilon < X_0\}}]}{\pr(X_0>a_n)} \\
&= \epsilon^{-\alpha} \theta \int_{1}^{\infty} \ex\Big[1- \exp\Big\{-\sum_{j=0}^{\infty}f(\epsilon y \mu_{A}^j, Z_j)\Big\}\Big] \alpha y^{-\alpha-1} \,\dd y \\
& = \theta \int_{\epsilon}^{\infty} \ex\Big[1- \exp\Big\{-\sum_{j=0}^{\infty}f( y \mu_{A}^j, Z_j)\Big\}\Big] \alpha y^{-\alpha-1} \,\dd y \\
&= \theta \int_{0}^{\infty} \ex\Big[1- \exp\Big\{-\sum_{j=0}^{\infty}f( y \mu_{A}^j, Z_j)\Big\}\Big] \alpha y^{-\alpha-1} \,\dd y,
\end{align*}
where the last line follows since $f(x,y)=0$ if $x\leq \epsilon$, $y\in\RR$, and $\mu_A\in(0,1)$, and this concludes the proof of (ii) yielding the statement.
\proofend

\section{From point processes to sums}

The key idea in handling the sums $\sum_{j=1}^n X_j^2$ and $\sum_{j=1}^n X_j M_{j+1}$ (the building blocks of $\hmuA^{(n)} - \mu_A$)
 is to apply usual truncation argument
and then summation to obtain the following result, cf.\ Davis and Hsing \cite[Theorem 3.1]{davis:hsing:1995}.

\begin{theorem}\label{thm:partial_sums}
We have
\begin{align} \label{eq:SumLim}
(V_n^{(1)},V_n^{(2)}):=\left(\frac{1}{a_n^2} \dsum_{j=1}^n X^2_j , \frac{1}{a_n^{3/2}} \dsum_{j=1}^{n} X_j M_{j+1} \right) \dto (V^{(1)},V^{(2)})
\end{align}
as $n \toi$ with
 \begin{align} \label{eq:SumLim2}
   (V^{(1)}, V^{(2)}) \eind \biggl(\frac{1}{1-\mu_A^2} \dsum_{i=1}^\infty P_i^{2}, \frac{1}{(1-\mu_A^3)^{1/2}} \dsum_{i=1}^\infty P_i^{3/2} Z_i\biggr) ,
 \end{align}
 where $\sum_{i=1}^\infty \delta_{P_i}$ is a Poisson point process on $(0,\infty)$ with intensity
 $\theta \dd(-y^{-\alpha})$ such that $P_1 \geq P_2 \geq \ldots$ almost surely, $(Z_i)_{i\geq 0}$ is
 an i.i.d.\ sequence  of $\cN(0,\sigma_A^2)$--distributed random variables independent of $\sum_{i=1}^\infty\delta_{P_i}$ with \ $\theta$ \ given in \eqref{theta},
 and the series on the right hand side of \eqref{eq:SumLim2} are convergent almost surely.

The characteristic function of the vector $(V^{(1)}, V^{(2)})$ has the form
 \[
   \ex\bigl[\ee^{\ii(sV^{(1)}+tV^{(2)})}\bigr]
   = \exp\biggl\{\theta \int_0^\infty \biggl(\exp\biggl\{\frac{\ii s}{1 - \mu_A^2} y^2 - \frac{\sigma_A^2 t^2}{2(1- \mu_A^3)} y^3\biggr\} - 1\biggr) \alpha y^{-\alpha -1} \,\dd y\biggr\}
 \]
 for $s, t \in \R$.
For the marginals, we have
 \[
   \ex\bigl[\ee^{\ii sV^{(1)}}\bigr]
   = \exp\Bigl\{- C_1 |s|^{\alpha/2} \Bigl( 1 - \ii \tan\Bigl(\frac{\pi\alpha}{4}\Bigr) \sgn(s)\Bigr)\Bigr\} , \qquad
   \ex\bigl[\ee^{\ii tV^{(2)}}\bigr] = \ee^{-C_2 |t|^{2\alpha/3}}
 \]
 for $s,t \in \R$ with
 \[
   C_1 := \theta \Gamma\Bigl(1 - \frac{\alpha}{2}\Bigr)  \frac{\cos\Bigl(\frac{\pi\alpha}{4}\Bigr)}{{(1-\mu_A^2)^{\alpha/2}}} , \qquad
   C_2 := \theta \Gamma\Bigl(1 - \frac{\alpha}{3}\Bigr) \biggl(\frac{\sigma_A^2}{2(1-\mu_A^3)}\biggr)^{\alpha/3} ,
 \]
 thus $V^{(1)}$ is an $\alpha/2$-stable positive random variable and $V^{(2)}$ is a symmetric $2\alpha/3$-stable random variable.
\end{theorem}

\noindent{\bf Proof.}
We divide the proof into four steps.

{\bf Step 1.}
Applying continuous mapping theorem to the convergence in (\ref{eq:PPconv}) we show that
\begin{align*}
 &(V_{n,\gamma}^{(1)},V_{n,\gamma}^{(2)})
  := \left( \dsum_{j=1}^{n} (X_{j}/a_n)^{2} \1{\{ X_j/a_n  > \gamma\}}, \dsum_{j=1}^{n} X_j M_{j+1}/a_n^{3/2} \1{\{ X_j/a_n  > \gamma\}} \right) \\
 &\dto
  (V_{\gamma}^{(1)},V_{\gamma}^{(2)}):=\left(\dsum_{i=1}^\infty\dsum_{j=0}^\infty (P_i \mu_A^j )^{2} \1{\{ P_i \mu_A^j > \gamma \}}\, , \,
  \dsum_{i=1}^\infty\dsum_{j=0}^\infty  {(P_i \mu_A^j)^{3/2} Z_{i,j} \1{ \{P_i \mu_A^j > \gamma\}}}\right)
\end{align*}
as $n \toi$ for all $\gamma>0$, where $Z_{i,j}$, $i \in \N$, $j \geq 0$, are given in Theorem \ref{th:PPconv}.
Fix $\gamma > 0$, and consider the mapping $T_\gamma:\mathcal{M}_p(S)\to \RR^2$,
 \begin{align*}
  T_\gamma(\kappa)
  &:= \left( \int_S x^2 \1{ \{x > \gamma\} } \, \kappa( \dd x, \dd y), \int_S x^{3/2} y \1{ \{x > \gamma\} } \, \kappa(\dd x, \dd y) \right) \\
  &= \biggl(\sum_k x_k^2 \1{\{x_k>\gamma\}}, \sum_k x_k^{3/2} y_k \1{\{x_k>\gamma\}}\biggr)
 \end{align*}
 for $\kappa = \sum_k \delta_{(x_k,y_k)} \in \mathcal{M}_p(S)$.
Note that the sums in the definition of $T_\gamma(\kappa)$ are sums with finitely many terms,
 since the set $(\gamma, \infty) \times \RR$ is bounded and $\kappa$ is a locally finite measure on $S$.
By Lemma \ref{lem:vagueConv_pointMeas}, $T_\gamma$ is continuous on
 $C_\gamma := \left\{ \kappa \in \mathcal{M}_p(S):  \kappa( \{ \gamma \} \times \R ) = 0 \right\}$,
 i.e., for any $\kappa \in C_\gamma$ and any sequence $(\kappa_n)_{n\in\NN}$ in $\mathcal{M}_p(S)$ such that $\kappa_n \stackrel{v}{\to} \kappa$ as $n\to\infty$,
 we have $T_\gamma(\kappa_n ) \to T_\gamma(\kappa)$ as $n\to\infty$.
Indeed, by Lemma \ref{lem:vagueConv_pointMeas}, $\kappa_n\stackrel{v}{\to}\kappa$ as $n\to\infty$ yields that
 there exist integers $n_0,M\geq 0$ and a labeling of the points of $\kappa$ and $\kappa_n$, $n \geq n_0$,
 in $(\gamma,\infty)\times \R$ such that
 \[
   \kappa_n|_{(\gamma,\infty)\times \R}=\sum_{k=1}^M \delta_{(x_k^{(n)}, \, y_k^{(n)})},
      \quad n\geq n_0,
     \qquad  \kappa|_{(\gamma,\infty)\times \R}=\sum_{k=1}^M \delta_{(x_k, y_k)},
 \]
 and $(x_k^{(n)},y_k^{(n)})\to (x_k,y_k)$ as $n \toi$ for all $k=1,\dots,M$.
Hence
 \begin{align*}
   &T_\gamma(\kappa_n)=\left(\sum_{k=1}^M (x_k^{(n)})^2, \sum_{k=1}^M  (x_k^{(n)})^{3/2}y_k^{(n)}\right),\qquad n\geq n_0,\\
   &T_\gamma(\kappa)=\left(\sum_{k=1}^M x_k^2, \sum_{k=1}^M  x_k^{3/2}y_k\right),
 \end{align*}
 yielding that $T_\gamma(\kappa_n ) \to T_\gamma(\kappa)$ as $n\to\infty$, as desired.
Finally, $(V_{n,\gamma}^{(1)},V_{n,\gamma}^{(2)})\dto (V_{\gamma}^{(1)},V_{\gamma}^{(2)})$ as $n\to\infty$ follows by an application of the continuous mapping theorem
 (see, e.g., Resnick \cite{resnick:2007}).
Indeed, $\mathcal{M}_p(S)$ is a complete separable metric space with a metric inducing the vague topology (see Kallenberg \cite[Lemma 4.6]{kallenberg:2017}),
 $(V_{n,\gamma}^{(1)},V_{n,\gamma}^{(2)}) = T_\gamma(N_n)$, $n\in\NN$, $T_\gamma(N) = (V_\gamma^{(1)},V_\gamma^{(2)})$, and
 we check that $\pr(N\in C_\gamma)=1$.
At the beginning of the proof of Theorem \ref{th:PPconv}, we already checked that $N$ is a locally finite measure almost surely,
 so it remains to verify that $\pr(N(\{\gamma\} \times \R) = 0) = 1$.
We have
\[
 N(\{\gamma\} \times \R)
     = \sum_{i=1}^\infty \sum_{j=0}^\infty \1{\{P_i\mu_A^j = \gamma\}}
    = \biggl(\sum_{i=1}^\infty \delta_{P_i}\biggr)(\{\gamma, \mu_A^{-1}\gamma,\mu_A^{-2}\gamma, \ldots\})
    = 0  \;\; \text{a.s.} \, ,
\]
 since the intensity measure of the point process $\sum_{i=1}^\infty \delta_{P_{i}}$ is absolutely continuous.

{\bf Step 2.}
We check that $(V_{\gamma}^{(1)}, V_{\gamma}^{(2)}) \Pto (V^{(1)}, V^{(2)})$ as $\gamma\downarrow 0$, where
 \begin{align*}
    V^{(1)}:=  (1-\mu_A^2)^{-1}\dsum_{i=1}^\infty P_i^{2},\qquad
    V^{(2)}:=  (1-\mu_A^3)^{-1/2}\dsum_{i=1}^\infty P_i^{3/2}Z_i
 \end{align*}
 with
 \[
   Z_i:= (1-\mu_A^3)^{1/2} \sum_{j=0}^\infty (\mu_A^j)^{3/2} Z_{i,j},\qquad i\in\NN,
 \]
 where $(P_i)_{i\in\NN}$ and $\{Z_{i,j} : i\in\NN, j\geq 0\}$ are given in Theorem \ref{th:PPconv}.

By monotonicity of $V_{\gamma}^{(1)}$ in $\gamma>0$,
\begin{align*}
V_{\gamma}^{(1)}\to (1-\mu_A^2)^{-1}\dsum_{i=1}^\infty P_i^{2} = V^{(1)}  \qquad \text{as $\gamma\downarrow 0$}
\end{align*}
almost surely.
By Campbell's theorem (see, e.g., Kingman \cite[Section 3.2]{kingman:1993}), we have $\pr(\sum_{i=1}^\infty P_i^\beta<\infty)=1$ for any $\beta \in (\alpha, \infty)$, and hence the series $\sum_{i=1}^\infty P_i^2$ is (absolutely) convergent almost surely, since $\alpha \in (1, 2)$.
Indeed, condition (3.16) in Kingman \cite{kingman:1993} is satisfied, since
 \begin{align*}
  \int_0^\infty (y^\beta\wedge 1)\theta \,\dd (-y^{-\alpha})
     & = \theta\alpha \int_0^\infty  (y^\beta\wedge 1) y^{-\alpha-1}\,\dd y
       = \theta\alpha \int_0^1 y^{\beta-\alpha-1}\,\dd y + \theta\alpha \int_1^\infty y^{-\alpha-1}\,\dd y \\
     & = \theta\alpha \left( \frac{1}{\beta-\alpha} + \frac{1}{\alpha}\right)
       <\infty.
 \end{align*}

Next, we show that the series $\sum_{j=0}^\infty (\mu_A^j)^{3/2} Z_{i,j}$, $i\in\NN$, and $\sum_{i=1}^\infty P_i^{3/2} Z_i$ are convergent almost surely.
Kolmogorov's one series theorem yields that $\sum_{j=0}^\infty  (\mu_A^j)^{3/2} Z_{i,j}$ converges almost surely for each $i\in\NN$, and hence
 $(Z_i)_{i\in\N}$ is an i.i.d. sequence of $\cN(0,\sigma_A^2)$-distributed random variables independent of $(P_i)_{i\in\N}$.
Let $E_j$, $j\in\N$, be i.i.d.\ random variables with an exponential distribution with parameter 1, independent of $Z_i$, $i \in \N$.
Put $\Gamma_i := \sum_{j=1}^i E_j$, $i \in \N$.
By the mapping theorem for Poisson random measures, $\sum_{i=1}^\infty \delta_{\theta^{1/\alpha} \Gamma_i^{-1/\alpha}}$ is a Poisson random measure on
 $(0,\infty)$ with intensity $\theta \dd(-y^{-\alpha})$, hence we have $(P_i)_{i\in\N} \eind (\theta^{1/\alpha} \Gamma_i^{-1/\alpha})_{i\in\N}$, and hence,
 $(P_i, Z_i)_{i\in\N} \eind (\theta^{1/\alpha} \Gamma_i^{-1/\alpha}, Z_i)_{i\in\N}$.
Consequently,
 \[
   \biggl(\sum_{i=1}^n P_i^{3/2} Z_i\biggr)_{n\in\N}
   \eind \biggl(\sum_{i=1}^n \bigl(\theta^{1/\alpha} \Gamma_i^{-1/\alpha}\bigr)^{3/2} Z_i\biggr)_{n\in\N}
   = \biggl(\theta^{3/(2\alpha)} \sum_{i=1}^n \Gamma_i^{-3/(2\alpha)} Z_i\biggr)_{n\in\N} \, ,
 \]
 thus the almost sure convergence of $\sum_{i=1}^\infty P_i^{3/2} Z_i$ will follow from the almost sure convergence of $\sum_{i=1}^\infty \Gamma_i^{-3/(2\alpha)}
Z_i$.
Indeed, if $\sum_{i=1}^\infty \Gamma_i^{-3/(2\alpha)} Z_i$ is convergent almost surely, then for each \ $\vare \in (0, \infty)$,
 by the continuity of probability, we have
 \begin{align*}
  \pr\biggl[\sup_{m\in\NN} \biggl|\sum_{i=n}^{n+m} P_i^{3/2} Z_i\biggr| > \vare\biggr]
  = \pr\biggl[\theta^{3/(2\alpha)} \sup_{m\in\NN} \biggl|\sum_{i=n}^{n+m} \Gamma_i^{-3/(2\alpha)} Z_i\biggr| > \vare\biggr] \to 0 \qquad \text{as $n \toi$,}
 \end{align*}
 implying the almost sure convergence of $\sum_{i=1}^\infty P_i^{3/2} Z_i$, see, e.g., Shiryaev \cite[Chapter II, Section 3, Theorem 1]{shi:1996}.
The almost sure convergence of $\sum_{i=1}^\infty \Gamma_i^{-3/(2\alpha)} Z_i$ follows from Theorem 1.4.5 in Samorodnitsky and Taqqu \cite{samorodnitsky:taqqu:1994},
 since $\frac{2\alpha}{3}\in(0,2)$, $\ex[|Z_1|^{2\alpha/3}] < \infty$, $\ex\bigl[|Z_1 \log(|Z_1|)|\bigr] < \infty$, $\ex[Z_1] = 0$ and $\ex\bigl[Z_1\int_{|Z_1|/i}^{|Z_1|/(i-1)} x^{-2} \sin(x) \, \dd x\bigr] = 0$ for all $i \in \NN$, where $|Z_1|/(i-1) := \infty$ for $i = 1$ (due to the fact that $Z_1$ is symmetric).
Hence the series $\sum_{i=1}^\infty P_i^{3/2} Z_i$ converges almost surely and $V^{(2)}$ is well--defined.

Next, we show that $V_\gamma^{(2)} \Pto V^{(2)}$ as $\gamma\downarrow 0$.
For every $\gamma, \epsilon \in (0, \infty)$ define
\begin{align*}
V_{\gamma,\epsilon}^{(2)} := \sum_{i=1}^\infty \sum_{j=0}^\infty (P_i \mu_A^j)^{3/2} Z_{i,j} \1{\{ P_i>\epsilon, \, P_i \mu_A^j  > \gamma \} } \, .
\end{align*}
Since for every $\epsilon \in (0, \infty)$ there are almost surely only finitely many $P_i$'s greater that $\epsilon$, for every fixed $\epsilon \in (0, \infty)$,
\begin{align*}
V_{\gamma,\epsilon}^{(2)} \to V_{0,\epsilon}^{(2)} :=  (1-\mu_A^3)^{-1/2}
                                              \sum_{i=1}^\infty P_i^{3/2} Z_i \1{\{ P_i>\epsilon \} } \qquad \text{as $\gamma\downarrow 0$}
\end{align*}
 almost surely.
Indeed, by the dominated convergence theorem, for every $i\in\NN$, we have
 $\sum_{j=0}^\infty \mu_A^{3j/2} Z_{i,j} \1{\{ P_i \mu_A^j  > \gamma \}}\to (1-\mu_A^3)^{-1/2} Z_i$ as $\gamma\downarrow 0$ almost surely, since
 $\vert \mu_A^{3j/2} Z_{i,j} \1{ \{P_i \mu_A^j > \gamma \} }\vert
 \leq \mu_A^{3j/2} \vert Z_{i,j}\vert$, $j\geq 0$, yielding that
 \[
   \left\vert \sum_{j=0}^\infty \mu_A^{3j/2} Z_{i,j} \1{\{ P_i \mu_A^j  > \gamma \}} \right\vert
      \leq \sum_{j=0}^\infty \mu_A^{3j/2} \vert Z_{i,j}\vert,
 \]
 where $\ex\Big[\sum_{j=0}^\infty \mu_A^{3j/2} \vert Z_{i,j}\vert\Big] = \ex[\vert Z_{1,1}\vert]/(1-\mu_A^{3/2}) <\infty$
 (especially, $\sum_{j=0}^\infty \mu_A^{3j/2} \vert Z_{i,j}\vert$ converges almost surely).

Now we check that $V_{0,\epsilon}^{(2)} \to V^{(2)}$ as $\vare \downarrow 0$ almost surely.
For each $\vare \in (0, \infty)$, we can write
 \[
   V_{0,\epsilon}^{(2)} = (1-\mu_A^3)^{-1/2} \sum_{i=1}^{K_\vare} P_i^{3/2} Z_i
 \]
 with $K_\vare := \max\{i \in \N : P_i > \vare\}$.
We have $K_\vare \to \infty$ as $\vare \downarrow 0$ almost surely, since $P_i \downarrow 0$ as $i \toi$ almost surely, thus
 the almost sure convergence of $\sum_{i=1}^\infty P_i^{3/2} Z_i$ yields that $V_{0,\epsilon}^{(2)} \to V^{(2)}$ as $\vare \downarrow 0$ almost surely.

For every $\gamma, \epsilon, \eta \in (0, \infty)$, we have
 \begin{align*}
  \pr(|V_\gamma^{(2)} - V^{(2)}|>\eta)
  \leq \pr(|V_\gamma^{(2)} - V_{\gamma,\epsilon}^{(2)}|>\eta/3) + \pr(|V_{\gamma,\epsilon}^{(2)} - V_{0,\epsilon}^{(2)}|>\eta/3)
   + \pr(|V_{0,\epsilon}^{(2)} - V^{(2)}|>\eta/3) .
 \end{align*}
The almost sure convergences $V_{\gamma,\epsilon}^{(2)} \asto V_{0,\epsilon}^{(2)}$ as $\gamma\downarrow 0$ for all $\vare \in (0, \infty)$ and $V_{0,\epsilon}^{(2)} \asto V^{(2)}$ as $\vare\downarrow 0$ imply the corresponding convergences in probability, hence
 \[
   \limsup_{\gamma\downarrow0} \pr(|V_\gamma^{(2)} - V^{(2)}|>\eta) \leq \limsup_{\gamma\downarrow0} \pr(|V_\gamma^{(2)} - V_{\gamma,\epsilon}^{(2)}|>\eta/3) + \pr(|V_{0,\epsilon}^{(2)} - V^{(2)}|>\eta/3)
 \]
 for every $\epsilon, \eta \in (0, \infty)$, and hence
 \[
   \limsup_{\gamma\downarrow0} \pr(|V_\gamma^{(2)} - V^{(2)}|>\eta) \leq \limsup_{\epsilon\downarrow0} \limsup_{\gamma\downarrow0} \pr(|V_{\gamma,\epsilon}^{(2)} - V_\gamma^{(2)}|>\eta/3)
 \]
 for every $\eta \in (0, \infty)$.
Consequently, if we show that for all $\eta \in (0, \infty)$,
\begin{align}\label{limsup_limsup}
\limsup_{\epsilon\downarrow0} \limsup_{\gamma\downarrow0} \pr(|V_{\gamma,\epsilon}^{(2)} - V_\gamma^{(2)}|>\eta)=0 \, ,
\end{align}
 then we obtain $V_\gamma^{(2)} \Pto V^{(2)}$ as $\gamma \downarrow 0$, as desired.
In order to check \eqref{limsup_limsup}, observe
\begin{align*}
  V_\gamma^{(2)} - V_{\gamma,\epsilon}^{(2)}
  = \sum_{i=1}^\infty \sum_{j=0}^\infty (P_i \mu_A^j)^{3/2} Z_{i,j} \1{\{ P_i\leq\epsilon, \, P_i \mu_A^j  > \gamma\} } \, ,
\end{align*}
since the sums defining $V_{\gamma,\epsilon}^{(2)}$ and $V_\gamma^{(2)}$ are sums with finitely many terms almost surely (see Step 1).
Since $\ex[V_{\gamma,\epsilon}^{(2)}- V_\gamma^{(2)}]=0$ for all $\gamma, \epsilon \in (0, \infty)$, we have
 \begin{align*}
  \Var(V_{\gamma,\epsilon}^{(2)}- V_\gamma^{(2)})
  &= \ex[(V_{\gamma,\epsilon}^{(2)}- V_\gamma^{(2)})^2] = \ex[\ex[(V_{\gamma,\epsilon}^{(2)}- V_\gamma^{(2)})^2 \mid (P_i)_{i\in\N}]] \\
  &= \ex\biggl[\sum_{i=1}^\infty \sum_{j=0}^\infty (P_i \mu_A^j)^3 \sigma_A^2 \1{\{ P_i\leq\epsilon, \, P_i \mu_A^j  > \gamma\} }\biggr]
  \leq \ex\biggl[\sum_{i=1}^\infty \sum_{j=0}^\infty (P_i \mu_A^j)^3 \sigma_A^2 \1{\{ P_i\leq\epsilon\} }\biggr] \\
  &\leq \frac{\sigma_A^2}{1-\mu_A^3} \ex\biggl[\sum_{i=1}^\infty P_i^3 \1{\{ P_i\leq\epsilon\} }\biggr]
  = \frac{\theta\sigma_A^2}{1-\mu_A^3} \int_0^\vare x^3 \alpha x^{-\alpha-1} \, \dd x
  = \frac{\theta\sigma_A^2\alpha\vare^{3-\alpha}}{(1-\mu_A^3)(3-\alpha)} ,
 \end{align*}
 where the last but one step follows by Campbell's theorem (see, e.g., Kingman \cite[Section 3.2]{kingman:1993}).
For all $\eta \in (0, \infty)$, Chebyshev's inequality implies that
\begin{align*}
\limsup_{\epsilon\downarrow0} \limsup_{\gamma\downarrow0} \pr(|V_{\gamma,\epsilon}^{(2)}- V_\gamma^{(2)}|>\eta)&\leq \limsup_{\epsilon\downarrow0}
 \frac{\theta \sigma_A^2\alpha\vare^{3-\alpha}}{\eta^2(1-\mu_A^3)(3-\alpha)} =0 ,
\end{align*}
hence we conclude \eqref{limsup_limsup}, as desired.
Altogether,
 \[
   (V^{(1)}_\gamma, V^{(2)}_\gamma) \stoch (V^{(1)}, V^{(2)}) \qquad \text{as $\gamma \downarrow 0$,}
 \]
 and hence in distribution as well.

{\bf Step 3.} By Billingsley \cite[Theorem 4.2]{billingsley:1968} and Steps 1 and 2,
using also that $\Vert (z_1,z_2)\Vert \leq \vert z_1\vert + \vert z_2\vert$, $(z_1,z_2)\in\R^2$,
 to show \eqref{eq:SumLim}, it suffices to prove
 that for all $\epsilon>0$,
\begin{align*}
 \lim_{\gamma\downarrow 0}\limsup_{n\toi} \pr(|V_{n}^{(k)} - V_{n,\gamma}^{(k)}|>\epsilon)= 0, \qquad k = 1, 2 .
\end{align*}
In case of $k=1$, by Markov's inequality, \eqref{a_n_quantile} and Karamata's theorem (see, Lemma \ref{truncated_moments}), we have
 \begin{align*}
  \lim_{\gamma \downarrow 0}\limsup_{n\toi} \pr(|V_{n}^{(1)} - V_{n,\gamma}^{(1)}|>\epsilon)
  &= \lim_{\gamma \downarrow 0}\limsup_{n\toi} \pr\left( \sum_{j=1}^n X_j^2 \1{\{ X_j/a_n \leq \gamma\}}  > a_n^2 \vare \right) \\
  &\leq \lim_{\gamma \downarrow 0}\limsup_{n\toi} \frac{\gamma^2}{\vare} \frac{\ex[X_1^2 \1{\{ X_1\leq a_n\gamma\}}]}{\gamma^2 a_n^2 \pr(X_1 > a_n\gamma)}
                                                    \frac{\pr(X_1 > a_n\gamma)}{\pr(X_1> a_n)} n\pr(X_1>a_n) \\
     &= \lim_{\gamma \downarrow 0} \frac{\alpha}{\vare(2-\alpha)}\gamma^{2-\alpha} =0,
 \end{align*}
 as desired.
In case of $k = 2$, we have to prove that
 \begin{align} \label{eq:VanSmallValues}
 & \lim_{\gamma \downarrow  0} \limsup_{n\toi}
\pr \left( \left| \dsum_{j=1}^{n} X_j M_{j+1} \1{\{ X_j/a_n  \leq \gamma\} }
 \right| > a_n^{3/2} \vare \right) =0
\end{align}
for all $\vare>0$.
Using the definition of $M_{j+1}$, we have
 \begin{align*}
  \pr \left( \left| \dsum_{j=1}^{n} X_j M_{j+1} \1{\{ X_j/a_n  \leq \gamma\} }
       \right| > a_n^{3/2} \vare \right)
  &\leq \pr \left( \left| \dsum_{j=1}^{n} X_j \sum_{i=1}^{X_j} \tilde{A}_i^{(j+1)}  \1{\{ X_j/a_n  \leq \gamma\} }
                   \right| > a_n^{3/2} \vare/2 \right) \\
  &\phantom{\leq}
    + \pr \left( \left| \dsum_{j=1}^{n} X_j \1{\{ X_j/a_n  \leq \gamma\} } \tilde B_{j+1}\right| > a_n^{3/2} \vare/2 \right).
 \end{align*}
Recall that for each $j=1,\dots, n$, $\{\tilde{A}_i^{(j+1)}, i\in \N\}$ are i.i.d.\ random variables
 with $\ex[\tilde{A}_1^{(j+1)}]=0$ and $\var[\tilde{A}_1^{(j+1)}]=\sigma_A^2$, and
independent of $\{\tilde{A}_i^{(k)},  i\in \N, k=2,\dots, j\}$ and $X_1,\dots, X_j$.
This implies that the random variables $X_j \sum_{i=1}^{X_j} \tilde{A}_i^{(j+1)} \1{\{ X_j/a_n  \leq \gamma \} }$, $j \in \{1,\ldots,n\}$, are uncorrelated
 and have zero expectation.
Hence, using Markov's inequality and the law of total variance, we get for all $\vare>0$,
 \begin{align*}
  &\pr \left( \left| \dsum_{j=1}^{n} X_j \sum_{i=1}^{X_j} \tilde{A}_i^{(j+1)}  \1{\{ X_j/a_n  \leq \gamma\} }
                   \right| > a_n^{3/2} \vare \right)
    \leq \frac{1}{a_n^3 \vare^2} \ex\left[  \left( \dsum_{j=1}^{n} X_j \sum_{i=1}^{X_j} \tilde{A}_i^{(j+1)}  \1{\{ X_j/a_n  \leq \gamma\} }
                                               \right)^2 \right] \\
  &= \frac{1}{a_n^3 \vare^2} \ex\left[ \dsum_{j=1}^{n} X_j^2 \left(\sum_{i=1}^{X_j} \tilde{A}_i^{(j+1)} \right)^2 \1{\{ X_j/a_n  \leq \gamma\} } \right]
    = \frac{n}{a_n^3 \vare^2} \ex\left[ X_0^2 \left(\sum_{i=1}^{X_0} \tilde{A}_i^{(1)} \right)^2 \1{\{ X_0/a_n  \leq \gamma\} } \right]\\
   & = \frac{n \Var\left[X_0 \left(\sum_{i=1}^{X_{0}} \tilde{A}_i^{(1)} \right) \1{\{ X_0/a_n  \leq \gamma \} }\right]}{a_n^3 \vare^2}
     = \frac{n \sigma_A^2 \ex[X_0^3 \1{ \{ X_0 \leq a_n \gamma\} }]}{a_n^3 \vare^2} .
\end{align*}
By Karamata's theorem (see, Lemma \ref{truncated_moments}) and \eqref{a_n_quantile}, since $\alpha<3$,
\begin{align*}
\frac{n \sigma_A^2 \ex[X_0^3 \1{ \{ X_0 \leq a_n \gamma \} }]}{a_n^3 \vare^2}
&= \sigma_A^2 \vare^{-2} \frac{\ex[X_0^3 \1{ \{ X_0 \leq a_n \gamma \} }]}{a_n^3 \gamma^3\pr(X_0 > a_n \gamma)} \gamma^3
  \frac{\pr (X_0 > a_n \gamma)}{\pr (X_0 > a_n)} \,n \pr (X_0 > a_n)\\
&\to \sigma_A^2 \vare^{-2} \frac{\alpha}{3-\alpha} \gamma^{3-\alpha}
 \qquad \text{as $n\toi$,}
\end{align*}
 which further goes to $0$ as $\gamma \downarrow 0$.
Hence, (\ref{eq:VanSmallValues}) will follow if we show that
\begin{align} \label{eq:EXB}
\lim_{\gamma \downarrow 0} \limsup_{n\toi}\pr\left(\left|\sum_{j=1}^n X_j  \1{ \{ X_j/a_n  \leq \gamma \} } \tilde{B}_{j+1} \right |
  > a_n^{3/2} \vare \right)=0
\end{align}
for all $\vare>0$.
With the notation $c_{B,n} := \ex [\tilde{B}_{j+1} \1{ \{ |\tilde B_{j+1}|/a_n  \leq 1 \} } ] = \ex [\tilde{B} \1{ \{ |\tilde{B}|/a_n  \leq 1 \} } ]$,
 where $\tilde{B}:=B-\mu_B$, we can write
 \begin{align*}
  \sum_{j=1}^n X_j \1{\{X_j/a_n\leq\gamma\}} \tilde{B}_{j+1}
  &= \sum_{j=1}^n X_j \1{\{X_j/a_n\leq\gamma\}} \Bigl(\tilde{B}_{j+1} \1{\{|\tilde B_{j+1}|/a_n  \leq1\}} - c_{B,n}\Bigr)
  + c_{B,n} \sum_{j=1}^n X_j \1{\{X_j/a_n\leq\gamma\}} \\
  &\quad
  + \sum_{j=1}^n X_j \1{\{X_j/a_n\leq\gamma\}} \tilde{B}_{j+1} \1{\{|\tilde B_{j+1}|/a_n>1\}}
   =: J^{(1)}_{n,\gamma} + J^{(2)}_{n,\gamma} + J^{(3)}_{n,\gamma} .
 \end{align*}
Since
 \begin{align*}
   &\pr(\vert J^{(1)}_{n,\gamma} + J^{(2)}_{n,\gamma} + J^{(3)}_{n,\gamma}\vert>  a_n^{3/2} \vare )\\
   &\leq \pr(\vert J^{(1)}_{n,\gamma} \vert>   a_n^{3/2} \vare /3)
        +  \pr(\vert J^{(2)}_{n,\gamma} \vert>   a_n^{3/2} \vare /3)
        + \pr(\vert J^{(3)}_{n,\gamma} \vert>   a_n^{3/2} \vare /3),
        \qquad \vare>0,
 \end{align*}
 to prove \eqref{eq:EXB}, it is enough to check that
 $\lim_{\gamma \downarrow 0} \limsup_{n\toi} \pr(\vert J^{(i)}_{n,\gamma} \vert>   a_n^{3/2} \vare) = 0$, $i=1,2,3$,
  for all $\vare>0$.

In case of $i=1$, using the independence of $X_j$ and $\tilde B_{j+1}$, Markov's inequality, and the facts that
 the summands in $J_{n,\gamma}^{(1)} $ are uncorrelated and $\Var[\tilde B \1{ \{ |\tilde B|\leq a_n \}} ]
 \leq \ex[\tilde B^2 \1{ \{ |\tilde B|\leq a_n \}}]$, we have
\begin{align*}
 &\pr \left( | J_{n,\gamma}^{(1)} | >  a_n^{3/2} \vare  \right)
   \leq \frac{n}{\vare^2 a_n^3} \ex[X_0^2 \1{ \{ X_0/a_n  \leq \gamma \} }]
        \ex[\tilde B^2 \1{ \{ |\tilde B|\leq a_n \} }] \\
  &\qquad = \frac{1}{\vare^2} \frac{ \ex[X_0^2 \1{ \{ X_0/a_n  \leq \gamma \} }] }{a_n^2 \gamma^2 \pr ( X_0>a_n   \gamma) }
     \frac{ \ex[\tilde B^2 \1{ \{ |\tilde B|\leq a_n \} }]  }{a_n^2  \pr (|\tilde B|> a_n  ) }
      \gamma^2 \frac{\pr ( X_0>a_n\gamma)}{\pr(X_0>a_n)} \, n  \pr( X_0>a_n) a_n \pr(|\tilde B|> a_n) .
\end{align*}
Indeed, the summands in $J_{n,\gamma}^{(1)} $ are uncorrelated, since for all $i<j$, $i,j\in\{1,\ldots,n\}$, we have
 \begin{align*}
  & \ex\left[ X_i  \1{ \{ X_i/a_n  \leq \gamma \} } \left(\tilde{B}_{i+1} \1{ \{ |\tilde B_{i+1}|/a_n  \leq 1 \} } -   c_{B,n}\right)
            X_j  \1{ \{ X_j/a_n  \leq \gamma \} } \left(\tilde{B}_{j+1} \1{ \{ |\tilde B_{j+1}|/a_n  \leq 1 \} } -   c_{B,n}\right)
       \right] \\
  & = \ex\left[ X_i X_j \1{ \{ X_i/a_n  \leq \gamma \} } \1{ \{ X_j/a_n  \leq \gamma \} }
              \left(\tilde{B}_{i+1} \1{ \{ |\tilde B_{i+1}|/a_n  \leq 1 \} } -   c_{B,n}\right)  \right]
      \ex\left[ \left(\tilde{B}_{j+1} \1{ \{ |\tilde B_{j+1}|/a_n  \leq 1 \} } -   c_{B,n}\right) \right] \\
  & = \ex\left[ X_i X_j \1{ \{ X_i/a_n  \leq \gamma \} } \1{ \{ X_j/a_n  \leq \gamma \} }
              \left(\tilde{B}_{i+1} \1{ \{ |\tilde B_{i+1}|/a_n  \leq 1 \} } -   c_{B,n}\right)  \right]\cdot 0
    = 0.
 \end{align*}
Note that, since $\tilde{B}$ is bounded from below and $B$ is regularly varying with tail index $\alpha$, we have
 \begin{align*}
 \pr(|\tilde{B}|>x)\sim \pr(\tilde{B}>x)\sim \pr(B>x)\qquad \text{as $x\toi$,}
 \end{align*}
 where we used that for all $\varepsilon>0$, $\pr(B>x(1+\varepsilon))\leq \pr(\tilde{B}>x)\leq \pr(B>x)$ for large enough $x$.
In particular, $|\tilde{B}|$ is regularly varying with tail index $\alpha$,
 and moreover, by (\ref{eq:tail_of_X}),
 \begin{align}\label{eq:tail_btilde}
  \pr(|\tilde{B}|>x)\sim  (1 - \mu_A^{\alpha})\pr(X_0>x) \qquad \text{as $x\toi$.}
 \end{align}
Consequently, by Karamata's theorem (see, Lemma \ref{truncated_moments}), Bingham et al. \cite[Proposition 1.3.6. (v)]{bingham:goldie:teugels:1987},
 \eqref{eq:tail_btilde} and \eqref{a_n_quantile}, for all $\gamma>0$,
 \begin{align*}
 \limsup_{n\toi}\pr \left( | J_{n,\gamma}^{(1)} | > a_n^{3/2} \vare  \right)
   & \leq \frac{\alpha^2 \gamma^{2-\alpha}}{\vare^2(2-\alpha)^2} \limsup_{n\to\infty} (a_n \pr (|\tilde B|> a_n ) )\\
   & = \frac{\alpha^2 \gamma^{2-\alpha} (1-\mu_A^\alpha)}{\vare^2(2-\alpha)^2} \limsup_{n\to\infty} (a_n\pr(X_0>a_n))
    = \frac{\alpha^2 \gamma^{2-\alpha} (1-\mu_A^\alpha)}{\vare^2(2-\alpha)^2}  \limsup_{n\to\infty} \frac{a_n}{n} \\
   & = \frac{\alpha^2 \gamma^{2-\alpha} (1-\mu_A^\alpha)}{\vare^2(2-\alpha)^2}  \limsup_{n\to\infty} (n^{\frac{1}{\alpha} -1} L(n) )
    \to 0  \qquad \text{as $n\toi$.}
 \end{align*}

In case of $i=2$, since $\ex[\tilde B] = 0$, by Markov's inequality, for all $\vare>0$,
 \begin{align*}
 \pr \left( | J^{(2)}_{n,\gamma} | >  a_n^{3/2} \vare  \right)
  &\leq \frac{n}{\vare a_n^{3/2}} \ex[X_0]\left|   \ex [\tilde{B}  \1{ \{ |\tilde B|/a_n  \leq 1 \} }]  \right|
  = \frac{n}{\vare a_n^{1/2}} \ex[X_0]
    \frac{ \left|  \ex [\tilde{B}  \1{ \{|\tilde B|/a_n > 1\} } ] \right|}
         { a_n  \pr (|\tilde B|> a_n) }
      \pr (|\tilde B|> a_n) \\
   & \leq \frac{\ex[X_0] }{\vare a_n^{1/2}}
    \frac{ \ex [| \tilde{B} | \1{ \{|\tilde B| > a_n\} } ] }
         { a_n  \pr (|\tilde B|> a_n) }
      n\pr (|\tilde B|> a_n)
   \to 0\qquad \text{as $n\toi$,}
 \end{align*}
 since, by Karamata's theorem (see, Lemma \ref{truncated_moments}),
  \[
     \lim_{n\toi} \frac{ \ex [| \tilde{B} | \1{ \{|\tilde B| > a_n\} } ] }
                       { a_n  \pr (|\tilde B|> a_n) } =\frac{\alpha}{\alpha-1},
  \]
 and, by \eqref{eq:tail_btilde} and \eqref{a_n_quantile}, $n \pr (|\tilde B|> a_n) \to 1- \mu_A^\alpha\in(0,1)$ as $n\to\infty$.

In case of $i=3$, similarly as in case of $i=2$, using Markov's inequality and the independence of $X_0$ and $B_1$,
 for all $\vare,\gamma>0$, we have
 \begin{align*}
  \pr \left( | J_{n,\gamma}^{(3)} | > a_n^{3/2} \vare  \right)
  &\leq \frac{n \ex(X_0 \1{ \{ X_0/a_n  \leq \gamma \} } \vert \tilde  B_1 \vert \1{ \{ \vert  \tilde B_1 \vert > a_n \} } ) }{a_n^{3/2} \vare} \\
   &\leq \frac{n}{\vare a_n^{1/2}} \ex[X_0]
    \frac{ \ex [| \tilde{B} | \1{ \{|\tilde B| > a_n\} } ] }
         { a_n  \pr (|\tilde B|> a_n) }
      \pr (|\tilde B|> a_n)
  \to 0\qquad \text{as $n\toi$,}
 \end{align*}
hence we conclude \eqref{eq:SumLim}.

{\bf Step 4.}
Finally, we determine the characteristic function of the random vector $(V^{(1)}, V^{(2)})$.
Using the continuity theorem, conditioning on $(P_i)_{i\in\N}$ and applying the portmanteau lemma, we have for $s, t \in \R$
 \begin{align*}
  \ex\bigl[\ee^{\ii(s V^{(1)} + t V^{(2)})}\bigr]
  &= \ex\biggl[\exp\biggl\{\ii\biggl(\frac{s}{1-\mu_A^2} \sum_{i=1}^\infty P_i^2
 + \frac{t}{(1-\mu_A^3)^{1/2}} \sum_{i=1}^\infty P_i^{3/2} Z_i \biggr)\biggr\}\biggr] \\
  &= \lim_{n\toi} \ex\biggl[\exp\biggl\{\ii\biggl(\frac{s}{1-\mu_A^2} \sum_{i=1}^n P_i^2
 + \frac{t}{(1-\mu_A^3)^{1/2}} \sum_{i=1}^n P_i^{3/2} Z_i \biggr)\biggr\}\biggr] \\
  &= \lim_{n\toi} \ex\biggl[\ex\biggl[\exp\biggl\{\ii\biggl(\frac{s}{1-\mu_A^2} \sum_{i=1}^n P_i^2
 + \frac{t}{(1-\mu_A^3)^{1/2}} \sum_{i=1}^n P_i^{3/2} Z_i \biggr)\biggr\} \,\bigg|\, (P_i)_{i\in\N}\biggr]\biggr] \\
  &= \lim_{n\toi} \ex\biggl[\exp\biggl\{\frac{\ii s}{1- \mu_A^2} \sum_{i=1}^n P_i^2 - \frac{\sigma_A^2 t^2}{2 (1 - \mu_A^3)} \sum_{i=1}^n P_i^3\biggr\}\biggr] \\
  &= \ex\biggl[\exp\biggl\{\frac{\ii s}{1- \mu_A^2} \sum_{i=1}^\infty P_i^2 - \frac{\sigma_A^2 t^2}{2 (1 - \mu_A^3)} \sum_{i=1}^\infty P_i^3\biggr\}\biggr] ,
 \end{align*}
 since the series $\sum_{i=1}^\infty P_i^\beta$ and $\sum_{i=1}^\infty P_i^{3/2}Z_i$ are convergent almost surely
  for any \ $\beta \in (\alpha, \infty)$ (see Step 2).
As in the proof of Campbell's theorem (see, e.g., Kingman \cite[Section 3.2]{kingman:1993}), one can prove that
 \begin{equation}\label{Kingman}
  \ex\biggl[\exp\biggl\{u \sum_{i=1}^\infty P_i^2 + v \sum_{i=1}^\infty P_i^3\biggr\}\biggr]
  = \exp\biggl\{\int_0^\infty \bigl(\ee^{uy^2+vy^3} - 1\bigr)
\theta \alpha y^{-\alpha -1} \, \dd y\biggr\}
 \end{equation}
 for any $u, v \in \C$ with $\Re(u) \leq 0$ and $\Re(v) \leq 0$, where $\Re(z)$ denotes the real part of $z \in \C$.
Indeed, \eqref{Kingman} holds for $u, v \in (-\infty, 0]$ by Campbell's theorem with the function $(0, \infty) \ni y \mapsto u y^2 + v y^3$ satisfying $\int_0^\infty (|u y^2 + v y^3| \land 1) \, \theta \alpha y^{-\alpha -1} \, \dd y < \infty$.
Since for any given  $u\in \C$ with $\Re(u) \leq 0$, both sides of \eqref{Kingman} as functions of $v$ are analytic functions on
 $\{ v \in \C :  \Re(v) < 0\}$, and for any given $v\in \C$ with $\Re(v) \leq 0$, both sides of \eqref{Kingman} as functions of $u$ are
 analytic functions on $\{ u \in \C :  \Re(u) < 0\}$, Hartogs's theorem yields that
 both sides of \eqref{Kingman} are analytic functions on $\{(u, v) \in \C^2 : \Re(u) < 0, \, \Re(v) < 0\}$.
So, by the identity theorem for analytic functions, \eqref{Kingman} holds on $\{(u, v) \in \C^2 : \Re(u) < 0, \, \Re(v) < 0\}$.
Both sides of \eqref{Kingman} are continuous functions on $\{(u, v) \in \C^2 : \Re(u) \leq 0, \, \Re(v) \leq 0\}$, so \eqref{Kingman} holds on $\{(u, v) \in \C^2 : \Re(u) \leq 0, \, \Re(v) \leq 0\}$.
Applying \eqref{Kingman} for $u = \ii s / (1 - \mu_A^2)$ and $v = - \sigma_A^2 t^2 / (2 (1 - \mu_A^3))$, we obtain the formula for $\ex\bigl[\ee^{\ii(s V^{(1)} + t V^{(2)})}\bigr]$, $s, t, \in \R$.

For each $\beta \in (0, 1)$ and $z \in \R$, we have
 \begin{equation}\label{Sato}
  \int_0^\infty (\ee^{\ii zr} - 1) r^{-1-\beta} \, \dd r
  = \Gamma(-\beta) \cos\Bigl(\frac{\pi\beta}{2}\Bigr) |z|^\beta \Bigl(1 - \ii \tan\Bigl(\frac{\pi\beta}{2}\Bigr) \sgn(z)\Bigr) ,
 \end{equation}
 see, e.g., the proof of Theorem 14.10 in Sato \cite{sato:2013}.
Applying \eqref{Sato} for $\beta = \alpha / 2$ and $z = s / (1 - \mu_A^2)$, we obtain for $s \in \R$
 \begin{align*}
  \ex\bigl[\ee^{\ii s V^{(1)}}\bigr]
  &= \exp\biggl\{\theta \int_0^\infty \biggl(\exp\biggl\{\frac{\ii s}{1 - \mu_A^2} y^2\biggr\} - 1\biggr) \alpha y^{-\alpha -1} \, \dd y\biggr\} \\
  &= \exp\biggl\{\theta \int_0^\infty \biggl(\exp\biggl\{\frac{\ii s}{1 - \mu_A^2} r\biggr\} - 1\biggr) \frac{\alpha}{2} r^{-\frac{\alpha}{2}-1} \, \dd r\biggr\} \\
  &= \exp\biggl\{\theta \frac{\alpha}{2} \Gamma\Bigl(-\frac{\alpha}{2}\Bigr) \cos\Bigl(\frac{\pi\alpha}{4}\Bigr) \Bigl|\frac{s}{1 - \mu_A^2}\Bigr|^{\alpha/2} \Bigl(1 - \ii \tan\Bigl(\frac{\pi\alpha}{4}\Bigr) \sgn\Bigl(\frac{s}{1 - \mu_A^2}\Bigr)\Bigr)
  \biggr\} \\
  &= \exp\biggl\{- \frac{\theta}{(1-\mu_A^2)^{\alpha/2}} \Gamma\Bigl(1 - \frac{\alpha}{2}\Bigr) \cos\Bigl(\frac{\pi\alpha}{4}\Bigr) |s|^{\alpha/2} \Bigl(1 - \ii \tan\Bigl(\frac{\pi\alpha}{4}\Bigr) \sgn(s)\Bigr)
  \biggr\} ,
 \end{align*}
 hence we obtain the characteristic function of $V^{(1)}$.

For each $\beta \in (0, 1)$ and $z \in [0, \infty)$, we have
 \begin{equation}\label{Li}
  \int_0^\infty (1 - \ee^{-zr}) \beta r^{-1-\beta} \, \dd r
  = \Gamma(1 - \beta) z^\beta ,
 \end{equation}
 see, e.g., Example 1.4 in Li \cite{li:2011} or the method of the proof of Example 8.11 in Sato \cite{sato:2013}.
Applying \eqref{Li} for $\beta = \alpha / 3$ and $z = \sigma_A^2 t^2 / (2(1 - \mu_A^3))$, we obtain for $t \in \R$
 \begin{align*}
  \ex\bigl[\ee^{\ii t V^{(2)}}\bigr]
  &= \exp\biggl\{\theta \int_0^\infty \biggl(\exp\biggl\{- \frac{\sigma_A^2 t^2}{2(1- \mu_A^3)} y^3\biggr\} - 1\biggr) \alpha y^{-\alpha -1} \,\dd y\biggr\} \\
  &= \exp\biggl\{\theta \int_0^\infty \biggl(\exp\biggl\{- \frac{\sigma_A^2 t^2}{2(1- \mu_A^3)} r\biggr\} - 1\biggr) \frac{\alpha}{3} r^{-\frac{\alpha}{3} -1} \, \dd r\biggr\} \\
  &= \exp\biggl\{-\theta \Gamma\Bigl(1 - \frac{\alpha}{3}\Bigr) \biggl(\frac{\sigma_A^2 t^2}{2(1- \mu_A^3)}\biggr)^{\alpha/3}\biggr\} ,
 \end{align*}
 hence we obtain the characteristic function of $V^{(2)}$, and this finishes the proof.
\proofend

\begin{remark}
 Note that the law of $(V^{(1)},V^{(2)})$ in \eqref{eq:SumLim} does not depend on the choice of the scaling sequence of $(a_n)_{n\in\NN}$ satisfying \eqref{a_n_quantile},
  which can be seen from the form of the characteristic function of $(V^{(1)},V^{(2)})$ given in Theorem \ref{thm:partial_sums}.
 \proofend
 \end{remark}

\begin{remark}
One can check that if $\alpha\in(1,3/2)$, then the series in the definition of $V^{(2)}$ in Theorem \ref{thm:partial_sums} is absolutely convergent
 almost surely.
By the mapping and marking theorems (see, e.g., Kingman \cite[Sections 2.3 and 5.2]{kingman:1993}),
 we have $\sum_{i=1}^\infty \delta_{(P_i^{3/2},|Z_i|)}$ is a Poisson random measure on $(0,\infty)\times(0,\infty)$
 with intensity measure $\theta\dd(-y^{-2\alpha/3})\times f_{|Z_1|}(z)\,\dd z$, where $f_{|Z_1|}$
 denotes the density function of $|Z_1|$.
Using again the mapping theorem, $\sum_{i=1}^\infty \delta_{P_i^{3/2}|Z_i|}$ is a Poisson random measure on $(0,\infty)$
 with intensity measure $\theta \ex[|Z_1|^{2\alpha/3}] \dd(-y^{-2\alpha/3}) $, since for any $t>0$,
 \begin{align*}
 &\int_{\{ (y,z) \in (0,\infty)^2 \,:\, yz>t \}}  \theta f_{|Z_1|}(z)\,\dd(-y^{-2\alpha/3})\,\dd z
    = \int_t^\infty \int_0^\infty \frac{2\alpha\theta}{3} \left(\frac{u}{z}\right)^{-\frac{2\alpha}{3} - 1} f_{|Z_1|}(z) \frac{1}{z}\,\dd u\,\dd z\\
 &= \frac{2\alpha\theta}{3} \ex[|Z_1|^{2\alpha/3}] \int_t^\infty u^{-\frac{2\alpha}{3} - 1}\,\dd u
    = \theta \ex[|Z_1|^{2\alpha/3}] \int_t^\infty 1 \, \dd(-u^{-2\alpha/3}).
 \end{align*}
Hence, by Campbell's theorem, we have $\sum_{i=1}^\infty P_i^{3/2} |Z_i|$ is convergent almost surely, since now
 $0<2\alpha/3<1$ and then
 \begin{align*}
  \int_0^\infty (y\wedge 1) \theta \ex[|Z_1|^{2\alpha/3}] \, \dd(-y^{-2\alpha/3})
   = \frac{2\alpha\theta}{3} \ex[|Z_1|^{2\alpha/3}] \left( \int_0^1 y^{-2\alpha/3}\,\dd y + \int_1^\infty y^{-2\alpha/3-1}\,\dd y \right)
   <\infty.
 \end{align*}
Further, since $\int_0^\infty y\theta \ex[\vert Z_1\vert^{2\alpha/3}]\, \dd(-y^{-2\alpha/3})=\infty$,
 by Kingman \cite[formula (3.18)]{kingman:1993}, we have $\ex[\sum_{i=1}^\infty P_i^{3/2} \vert Z_i\vert]=\infty$.
\proofend
\end{remark}

\section{On the limit behavior of the CLS estimator}
\label{section:properties}

Now we can formulate our main result.

\begin{theorem}\label{th:main}
We have
 \[
   \sqrt{a_n} (\hmuA^{(n)} - \mu_A)  \dto \frac{V^{(2)}}{V^{(1)}}
   \qquad \text{as $n \toi$,}
 \]
 where the sequence $(a_n)_{n\in\NN}$ and the joint characteristic function of $(V^{(1)},V^{(2)})$ is given
 in \eqref{a_n_quantile} and in Theorem \ref{thm:partial_sums}, respectively.
\end{theorem}

\noindent{\bf Proof.}
By Theorem \ref{thm:partial_sums}, $V^{(1)}$ is an $\alpha/2$-stable positive random variable, thus it is absolutely continuous and $\pr(V^{(1)} > 0) = \pr(V^{(1)} \ne 0) = 1$.
For each $n \in \Nset$, by the strong stationarity of $(X_i)_{i\in\Zset}$, we have
 \[
   \sqrt{a_n} (\hmuA^{(n)} - \mu_A)
   \eind \frac{\sum_{j=1}^n X_j M_{j+1}/a_n^{3/2}}{\sum_{j=1}^n (X_j/a_n)^2} .
 \]
Consequently, by Theorem \ref{thm:partial_sums} and the continuous mapping theorem (see, e.g., Billingsley \cite[Theorem 5.1]{billingsley:1968}), we conclude the statement.
\proofend

\begin{remark}
The limit law $V^{(2)}/V^{(1)}$ in Theorem \ref{th:main} can be written in the form
  \[
  \frac{V^{(2)}}{V^{(1)}} \eind
    \frac{(1-\mu_A^2) \sigma_A }{(1-\mu_A^3)^{1/2}(1-\mu_A^\alpha)^{\frac{1}{2\alpha}}}
      \frac{\dsum_{i=1}^\infty (\widetilde P_i)^{\frac{3}{2}}\widetilde Z_i }{\dsum_{i=1}^\infty (\widetilde P_i)^2},
 \]
 where $\sum_{i=1}^\infty \delta_{\widetilde P_i}$ is a Poisson random measure on $(0,\infty)$ with intensity measure $\dd(-y^{-\alpha})$
 being independent of an i.i.d. sequence of $\cN(0,1)$-distributed random variables $(\widetilde Z_i)_{i\in\N}$.
Indeed, using the mapping theorem for Poisson random measures,
 one can check that $\sum_{i=1}^\infty \delta_{\theta^{-\frac{1}{\alpha}} P_i}$ is a Poisson random measure on $(0,\infty)$
 with intensity measure $\dd(-y^{-\alpha})$.
Consequently, we have
 \begin{align*}
   (V^{(1)}, V^{(2)})
     & \eind \left(\frac{1}{1-\mu_A^2} \theta^{\frac{2}{\alpha}} \dsum_{i=1}^\infty (\theta^{-\frac{1}{\alpha}} P_i)^{2},
             \frac{1}{(1-\mu_A^3)^{1/2}} \theta^{\frac{3}{2\alpha}} \dsum_{i=1}^\infty (\theta^{-\frac{1}{\alpha}} P_i)^{\frac{3}{2}} Z_i
       \right)\\
     & \eind \left(\frac{(1-\mu_A^\alpha)^{\frac{2}{\alpha}}}{1-\mu_A^2} \dsum_{i=1}^\infty (\widetilde P_i)^2,
                         \frac{(1-\mu_A^\alpha)^{\frac{3}{2\alpha}}\sigma_A }{(1-\mu_A^3)^{1/2}} \dsum_{i=1}^\infty (\widetilde P_i)^{\frac{3}{2}}  \widetilde Z_i
                   \right),
 \end{align*}
 yielding the statement.

Note that $\sum_{i=1}^\infty (\widetilde P_i)^{\frac{3}{2}} \widetilde Z_i / \sum_{i=1}^\infty (\widetilde P_i)^2$
 does not depend on the parameter $\mu_A$ to be estimated nor on $\sigma_A$.
This gives the possibility for formulating a version of Theorem \ref{th:main} with a random normalization such that the limit law does not depend
 on $\mu_A$ and $\sigma_A$.
\proofend
\end{remark}

In what follows, we collect several interesting properties of $(V^{(1)}, V^{(2)})$ and $V^{(2)}/V^{(1)}$.
The characteristic function of a random vector $\bX$ will be denoted by $\varphi_\bX$.

\begin{proposition}
The distribution of $(V^{(1)}, V^{(2)})$ is operator stable, and the matrix
 $\diag_2\bigl(\frac{2}{\alpha}, \frac{3}{2\alpha}\bigr)$ is an exponent of it.
Particularly, the distribution of $(V^{(1)}, V^{(2)})$ is full and infinitely divisible,
 and has an infinitely differentiable density function, and the partial derivatives of this density function tend to 0 at infinity.
\end{proposition}

\noindent{\bf Proof.}
First, observe that for each $a \in (0, \infty)$ and $s, t \in \R$, by the substitution $a^{1/\alpha} y = x$, we obtain
 \begin{equation}\label{OS}
  \begin{aligned}
   &\varphi_{a^{2/\alpha}V^{(1)},\,a^{3/(2\alpha)}V^{(2)}}(s, t) \\
   &= \exp\biggl\{\theta \int_0^\infty \biggl(\exp\biggl\{\frac{\ii a^{2/\alpha}s}{1 - \mu_A^2} y^2 - \frac{\sigma_A^2 (a^{3/(2\alpha)}t)^2}{2(1- \mu_A^3)} y^3\biggr\} - 1\biggr) \alpha y^{-\alpha -1} \,\dd y\biggr\} \\
   &= \exp\biggl\{a \theta \int_0^\infty \biggl(\exp\biggl\{\frac{\ii s}{1 - \mu_A^2} x^2 - \frac{\sigma_A^2 t^2}{2(1- \mu_A^3)} x^3\biggr\} - 1\biggr) \alpha x^{-\alpha -1} \,\dd x\biggr\} \\
   &= (\varphi_{V^{(1)},\,V^{(2)}}(s, t))^a ,
  \end{aligned}
 \end{equation}
 hence Equation (7.8) in Meerschaert and Scheffler \cite{meerschaert:scheffler:2001} is satisfied with exponent $\diag_2\bigl(\frac{2}{\alpha}, \frac{3}{2\alpha}\bigr)$ and without shifts.

Particularly, with $a = n^{-1}$, $n \in \NN$, we get
 \[
   \varphi_{V^{(1)},\,V^{(2)}}(s, t)
   = \bigl(\varphi_{n^{-2/\alpha}V^{(1)},\,n^{-3/(2\alpha)}V^{(2)}}(s, t)\bigr)^n
 \]
 for all $s, t \in \R$ and $n \in \N$, hence the distribution of $(V^{(1)}, V^{(2)})$ is infinitely divisible.

In order to prove that the distribution of $(V^{(1)}, V^{(2)})$ is full, we have to show that for each $(v_1, v_2) \in \R^2 \setminus \{(0, 0)\}$, the random variable $v_1 V^{(1)} + v_2 V^{(2)}$ is nondegenerate.
If we suppose that, on the contrary, there exist $(v_1, v_2) \in \R^2 \setminus \{(0, 0)\}$ and $x_0 \in \RR$ such that $v_1 V^{(1)} + v_2 V^{(2)} = x_0$ almost surely, then for each $t \in \R$, we would have
 \begin{align*}
  \ee^{\ii tx_0}
  &= \varphi_{v_1V^{(1)}+v_2V^{(2)}}(t)
   = \varphi_{V^{(1)},\,V^{(2)}}(v_1t, v_2t) \\
  &= \exp\biggl\{\theta \int_0^\infty \biggl(\exp\biggl\{\frac{\ii v_1t}{1 - \mu_A^2} y^2 - \frac{\sigma_A^2 (v_2t)^2}{2(1- \mu_A^3)} y^3\biggr\} - 1\biggr) \alpha y^{-\alpha -1} \,\dd y\biggr\} .
 \end{align*}
We have $\varphi_{V^{(1)},\,V^{(2)}}(v_1t, v_2t) \ne 0$ for any $t \in \R$, since the distribution of $(V^{(1)}, V^{(2)})$ is infinitely divisible.
Applying Lemma 7.6 in Sato \cite{sato:2013}, we would obtain
 \[
   \ii tx_0
   = \theta \int_0^\infty \biggl(\exp\biggl\{\frac{\ii v_1t}{1 - \mu_A^2} y^2 - \frac{\sigma_A^2 v_2^2 t^2}{2(1- \mu_A^3)} y^3\biggr\} - 1\biggr) \alpha y^{-\alpha -1} \,\dd y , \qquad t \in \R .
 \]
Taking the real parts of both sides, we would get
 \[
   0 = \theta \int_0^\infty \biggl(\exp\biggl\{- \frac{\sigma_A^2 v_2^2 t^2}{2(1- \mu_A^3)} y^3\biggr\} \cos\biggl(\frac{v_1t}{1 - \mu_A^2} y^2\biggr)  - 1\biggr) \alpha y^{-\alpha -1} \,\dd y , \qquad t \in \R .
 \]
Since the integrand is continuous and nonpositive, we would conclude that it is identically zero, yielding that
 \[
   \cos\biggl(\frac{v_1t}{1 - \mu_A^2} y^2\biggr) = \exp\biggl\{\frac{\sigma_A^2 v_2^2 t^2}{2(1- \mu_A^3)} y^3\biggr\} , \qquad y \in [0, \infty) , \qquad t \in \R,
 \]
 which is a contradiction due to $(v_1,v_2)\ne (0,0)$, hence the distribution of $(V^{(1)}, V^{(2)})$ is full.
By Theorem 7.2.1 in Meerschaert and Scheffler \cite{meerschaert:scheffler:2001}, taking into account \eqref{OS}, we obtain that the distribution of $(V^{(1)}, V^{(2)})$ is operator stable and the matrix $\diag_2\bigl(\frac{2}{\alpha}, \frac{3}{2\alpha}\bigr)$ is an exponent of it.
For the facts that $(V^{(1)}, V^{(2)})$ has an infinitely differentiable density function,
 and the partial derivatives of this density function tend to 0 at infinity,
 see {\L}uczak \cite[Corollary 2.1]{Luc} and Kern and Wedrich \cite[page 387]{KerWed}.
\proofend

\begin{proposition}
The random variables $V^{(1)}$ and $V^{(2)}$ are dependent.
\end{proposition}

\noindent{\bf Proof.}
If we suppose that, on the contrary, $V^{(1)}$ and $V^{(2)}$ are independent, then we would have $\varphi_{V^{(1)},\,V^{(2)}}(s, t) = \varphi_{V^{(1)}}(s) \varphi_{V^{(2)}}(t)$ for all $s, t \in \R$, hence
 \begin{align*}
  &\exp\biggl\{\theta \int_0^\infty \biggl(\exp\biggl\{\frac{\ii s}{1 - \mu_A^2} y^2 - \frac{\sigma_A^2 t^2}{2(1- \mu_A^3)} y^3\biggr\} - 1\biggr) \alpha y^{-\alpha -1} \,\dd y\biggr\}
  \end{align*}
  \begin{align*}
  &= \exp\biggl\{\theta \int_0^\infty \biggl(\exp\biggl\{\frac{\ii s}{1 - \mu_A^2} y^2\biggr\} - 1\biggr) \alpha y^{-\alpha -1} \,\dd y\biggr\} \\
  &\quad\times
  \exp\biggl\{\theta \int_0^\infty \biggl(\exp\biggl\{- \frac{\sigma_A^2 t^2}{2(1- \mu_A^3)} y^3\biggr\} - 1\biggr) \alpha y^{-\alpha -1} \,\dd y\biggr\} \\
  &= \exp\biggl\{\theta \int_0^\infty \biggl(\exp\biggl\{\frac{\ii s}{1 - \mu_A^2} y^2\biggr\} + \exp\biggl\{- \frac{\sigma_A^2 t^2}{2(1- \mu_A^3)} y^3\biggr\} - 2\biggr) \alpha y^{-\alpha -1} \,\dd y\biggr\}
 \end{align*}
 for all $s, t \in \R$.
We have $\varphi_{V^{(1)},\,V^{(2)}}(s, t) \ne 0$, $\varphi_{V^{(1)}}(s) \ne 0$ and $\varphi_{V^{(2)}}(t) \ne 0$ for any $s, t \in \R$, since the distributions of $(V^{(1)}, V^{(2)})$, $V^{(1)}$ and $V^{(2)}$ are infinitely divisible.
Applying Lemma 7.6 in Sato \cite{sato:2013}, we would obtain
 \begin{align*}
  &\theta \int_0^\infty \biggl(\exp\biggl\{\frac{\ii s}{1 - \mu_A^2} y^2 - \frac{\sigma_A^2 t^2}{2(1- \mu_A^3)} y^3\biggr\} - 1\biggr) \alpha y^{-\alpha -1} \,\dd y \\
  &= \theta \int_0^\infty \biggl(\exp\biggl\{\frac{\ii s}{1 - \mu_A^2} y^2\biggr\} + \exp\biggl\{- \frac{\sigma_A^2 t^2}{2(1- \mu_A^3)} y^3\biggr\} - 2\biggr) \alpha y^{-\alpha -1} \,\dd y
 \end{align*}
 for all $s, t \in \R$, hence
 \[
   \int_0^\infty \biggl(\exp\biggl\{\frac{\ii s}{1 - \mu_A^2} y^2\biggr\} - 1\biggr) \biggl(\exp\biggl\{- \frac{\sigma_A^2 t^2}{2(1- \mu_A^3)} y^3\biggr\} - 1\biggr) \alpha y^{-\alpha -1} \,\dd y = 0
 \]
 for all $s, t \in \R$.
Taking the real parts of both sides, we would get
 \[
   \int_0^\infty \biggl(\cos\biggl(\frac{s}{1 - \mu_A^2} y^2\biggr) - 1\biggr) \biggl(\exp\biggl\{- \frac{\sigma_A^2 t^2}{2(1- \mu_A^3)} y^3\biggr\} - 1\biggr) \alpha y^{-\alpha -1} \,\dd y = 0
 \]
 for all $s, t \in \R$.
Since the integrand is continuous and nonnegative, we would conclude that it is identically zero, yielding a contradiction
 unless $s=t=0$, hence the random variables $V^{(1)}$ and $V^{(2)}$ are dependent.
\proofend

Note that, by Theorem \ref{thm:partial_sums}, $\ex[V^{(1)}] = \infty$ for all $\alpha\in(1,2)$,
 and $\ex[V^{(2)}]$ does not exist if $\alpha \in \bigl(1, \frac{3}{2}\bigr]$ and $\ex[V^{(2)}] = 0$ if $\alpha \in \bigl(\frac{3}{2}, 2\bigr)$.
In what follows we show that all the exponential moments of $V^{(2)}/V^{(1)}$ are finite.

\begin{proposition}
For each $t \in \R$, we have
\[
\ex \left[ \exp \left\{ t \frac{V^{(2)}}{V^{(1)}} \right\} \right]< \infty.
\]
\end{proposition}

\noindent{\bf Proof.}
Let $t \in \R$ be fixed.
Using that the series $\sum_{i=1}^\infty P_i^\beta$ and $\sum_{i=1}^\infty P_i^{3/2}Z_i$ are convergent almost surely
 for any \ $\beta \in (\alpha, \infty)$ (see Step 2 of the proof of Theorem \ref{thm:partial_sums}),
 conditioning on $(P_i)_{i\in \N}$, by the continuity theorem,  we have
\begin{equation*} 
\ex \left[ \exp \left\{ t \frac{V^{(2)}}{V^{(1)}} \right\} \right]
= \ex \left[ \exp \left\{ \frac{t^2 \sigma_A^2 (1 - \mu_A^2)^2}{2(1- \mu_A^3)}
\frac{\sum_{i=1}^\infty P_i^3}
{\left(\sum_{i=1}^\infty P_i^2 \right)^2}\right\} \right].
\end{equation*}
So we need to check that all the exponential moments of
\begin{equation} \label{eq:defU}
U := \theta^{1/\alpha}\frac{\sum_{i=1}^\infty P_i^3}
{\left(\sum_{i=1}^\infty P_i^2 \right)^2}
\end{equation}
 are finite.
Then $\pr(U\in(0,\infty))=1$, since the series $\sum_{i=1}^\infty P_i^2$ and $\sum_{i=1}^\infty P_i^3$ are absolutely convergent with positive sums almost surely.
Recall that the Poisson point process $(P_i)_{i\in\N}$ in Theorem \ref{thm:partial_sums} can be represented as
\[
 (P_i)_{i\in\N} \eind \bigl( \theta^{1/\alpha} \Gamma_i^{-1/\alpha} \bigr)_{i\in\N}
\]
with $\Gamma_i = E_1 + \ldots + E_i$, $i\in\N$,
 where $E_j$, $j\in\N$, are i.i.d. random variables with an exponential distribution with parameter 1 independent of $(Z_i)_{i\in\N}$.
Hence
\begin{equation*} 
U \eind \frac{\sum_{i=1}^\infty \Gamma_i^{-3/\alpha}}
{\bigl(\sum_{i=1}^\infty \Gamma_i^{-2/\alpha}\bigr)^2}.
\end{equation*}
Since
\[
\Biggl(\sum_{i=1}^\infty \Gamma_i^{-2/\alpha}\Biggr)^2 \geq
\sum_{i=1}^\infty \Gamma_i^{-4/\alpha} +
\Gamma_1^{-2/\alpha} \sum_{i=2}^\infty \Gamma_i^{-2/\alpha},
\]
we see that $U > x$ with $x > 0$, implies that
\[
\Gamma_1^{-3/\alpha} > \Gamma_1^{-4/\alpha} x \ \qquad \text{or } \qquad
\Gamma_i^{-3/\alpha} > \bigl( \Gamma_i^{-4/\alpha} +
\Gamma_1^{-2/\alpha} \Gamma_i^{-2/\alpha}\bigr) x \quad \text{for some $i \geq 2$.}
\]
In both cases we have $\Gamma_1 > x^\alpha$.
Indeed, if $\Gamma_i^{-3/\alpha} > \bigl(\Gamma_i^{-4/\alpha} + \Gamma_1^{-2/\alpha} \Gamma_i^{-2/\alpha}\bigr) x$
 with some $i\geq 2$, then $\Gamma_i^{-3/\alpha} >  \Gamma_1^{-2/\alpha} \Gamma_i^{-2/\alpha} x$ yielding that $\Gamma_1^{2/\alpha} \Gamma_i^{-1/\alpha} > x$,
 and then, since $\Gamma_1\leq \Gamma_i$, we have $x<\Gamma_1^{1/\alpha}$, as desired.
Summarizing, we have shown that
 \[
 \pr( U > x  ) \leq \pr( \Gamma_1 > x^{\alpha}) = \ee^{-x^\alpha},\qquad x > 0 ,
 \]
 which yields the statement.
Indeed, if $s\leq 0$, then $\ex[\ee^{sU}]\leq 1<\infty$, since $\pr(U\in(0,\infty))=1$, and if $s>0$, then
 \begin{align*}
  \ex[\ee^{sU}] & = \int_0^\infty \pr(\ee^{sU} > x)\,\dd x
                  =  \int_0^\infty \pr(U > \ln(x)/s)\,\dd x
                  \leq \int_0^{\ee^s} 1\,\dd x + \int_{\ee^s}^\infty \ee^{-(\ln(x)/s)^\alpha} \,\dd x \\
                & = \ee^s+ \int_s^\infty \ee^{y-(y/s)^\alpha} \,\dd y
                  <\infty,
 \end{align*}
 since $\alpha\in(1,2)$, and $y - (y/s)^\alpha < -1/(2s^\alpha) y$ for large enough $y>s$.
\proofend

\begin{proposition}\label{Lem_abs_cont}
The random variable $V^{(2)} / V^{(1)}$ has a continuously differentiable density function.
\end{proposition}

\noindent{\bf Proof.}
Using that the series $\sum_{i=1}^\infty P_i^\beta$ and $\sum_{i=1}^\infty P_i^{3/2}Z_i$ are convergent almost surely
 for any \ $\beta \in (\alpha, \infty)$ (see Step 2 of the proof of Theorem \ref{thm:partial_sums}),
 conditioning on $(P_i)_{i\in \N}$, by the continuity theorem, for any $t \in \R$, we have
\[
\ex\left[ \exp \left\{ \ii  t \frac{V^{(2)}}{V^{(1)}} \right\} \right]
= \ex \left[\exp \left\{ - \frac{t^2 \sigma_A^2 (1 - \mu_A^2)^2}{2(1- \mu_A^3)}
\frac{\sum_{i=1}^\infty P_i^3}
{\left(\sum_{i=1}^\infty P_i^2 \right)^2}\right\} \right] .
\]
To prove the existence of a uniformly continuous and continuously differentiable density function,
 it is enough to check that
 \[
 \int_{-\infty}^\infty \vert t\vert \ex [\ee^{- C_U t^2 U/2} ]\, \dd t
     = C_U^{-1} \int_{-\infty}^\infty \vert t\vert \ex [\ee^{-t^2 U/2}] \,\dd t< \infty,
 \]
 where $U$ is given in \eqref{eq:defU} and
\[
 C_U := \frac{\theta^{-1/\alpha} \sigma_A^2 (1 - \mu_A^2)^2}{1- \mu_A^3} ,
\]
 see, e.g., Sato \cite[Proposition 28.1]{sato:2013}.
Here, using Fubini's theorem,
 \begin{align*}
 \int_{-\infty}^\infty \vert t\vert \ex[ \ee^{-t^2 U/2}] \dd t
  &= \sqrt{2 \pi} \ex \left[ U^{-1/2}
         \int_{-\infty}^\infty \vert t\vert
              \frac{1}{\sqrt{2 \pi U^{-1}}}
              \ee^{-t^2/ (2U^{-1})} \dd t \right]  \\
 & =  \sqrt{2 \pi}  \ex[ U^{-1/2} \ex[U^{-1/2}\vert Z\vert \mid U]] \\
 & = \sqrt{2 \pi} \ex[U^{-1} \vert Z\vert ]
    = \sqrt{2 \pi}  \ex [U^{-1}] \ex [|Z|],
\end{align*}
 where $Z$ is a standard normally distributed random variable independent of $U$.
Thus we only have to show that $\ex[U^{-1}] < \infty$, i.e., by \eqref{eq:defU},
\[
 \ex \left[ \frac{ \left(\sum_{i=1}^\infty \Gamma_i^{-2/\alpha} \right)^2}
            {\sum_{i=1}^\infty \Gamma_i^{-3/\alpha} } \right] < \infty.
\]
In what follows we will use the following facts:
 \begin{itemize}
  \item $(\Gamma_i)_{i\geq2} \eind (\Gamma_1 + \Gamma_{i-1}')_{i\geq2}$,
       where $(\Gamma_i')_{i\geq1}$ has the same distribution as
     $(\Gamma_i)_{i\in\N}$ and independent of it,
  \item by Campbell's theorem,
       \[
      \ex \left[\sum_{i=1}^\infty h(\Gamma_i)\right] = \ex\left[ \left(\sum_{i=1}^\infty \delta_{\Gamma_i} \right)(h)\right]  = \int_0^\infty h(y)\, \dd y
    \]
    for any Borel measurable function $h:(0,\infty)\to\RR$ in the sense that the expectations exist on the left hand side
    if and only if the integral on the right hand side converges and then they are equal,
 \item if $\int_0^\infty h(x) \dd x$ converges, then
      \[
       \ex\left[ \left(\sum_{i=1}^\infty h(\Gamma_i)\right)^2\right] = \int_0^\infty h^2(y) \,\dd y + \left(\int_0^\infty h(y) \,\dd y\right)^2,
       \]
    where the right hand side can be finite or infinite as well.
  \end{itemize}
So, by conditioning on $\Gamma_1$ having an exponential distribution with parameter 1, we obtain
 \begin{align*}
 & \ex \left[\frac{\left( \sum_{i=1}^\infty \Gamma_i^{-2/\alpha} \right)^2}
         {\sum_{i=1}^\infty \Gamma_i^{-3/\alpha} }  \right]
   \leq \ex \left[ \left( \sum_{i=1}^\infty \Gamma_i^{-2/\alpha} \right)^2 \Gamma_1^{3/\alpha} \right]\\
 & = \int_0^\infty x^{3/\alpha} \ee^{-x}
  \ex \left[ \left( x^{-2/\alpha} +
              \sum_{i=1}^\infty (x + \Gamma_i')^{-2/\alpha} \right)^2 \right] \dd x \\
 &= \int_0^\infty x^{3/\alpha} \ee^{-x}
     \left( x^{-4/\alpha} + 2x^{-2/\alpha} \ex\left[ \sum_{i=1}^\infty (x + \Gamma_i')^{-2/\alpha} \right]
                      + \ex\left[ \left(\sum_{i=1}^\infty (x + \Gamma_i')^{-2/\alpha} \right)^2  \right] \right)  \dd x \\
 & = \int_0^\infty x^{3/ \alpha} \ee^{-x}
   \left( x^{-4/\alpha} + 2 x^{-2/\alpha} \int_{x}^\infty y^{-2/\alpha} \, \dd y
    + \int_x^\infty y^{-4/\alpha} \, \dd y + \left(\int_x^\infty y^{-2/\alpha}\,\dd y \right)^2 \right) \dd x \\
& = \int_0^\infty \ee^{-x}
    \left( x^{-1/\alpha} + \left( \frac{2 \alpha}{2 - \alpha} + \frac{\alpha}{4 - \alpha} \right) x^{1 - 1/\alpha}
           + \frac{\alpha^2}{(2-\alpha)^2}x^{2-1/\alpha}  \right) \dd x < \infty,
\end{align*}
 as desired, since $\int_0^\infty x^{n-\frac{1}{\alpha}}\ee^{-x}\, \dd x = \Gamma(n+1-\frac{1}{\alpha})<\infty$, $n\geq 0$.
\proofend

\begin{proposition}
For each $x\in\RR$, we have
\[
\pr \left( \frac{V^{(2)}}{V^{(1)}} \leq x \right)
=
\frac{1}{2} - \frac{1}{2 \pi \ii} \int_{-\infty}^\infty \frac{\varphi_{V^{(1)},V^{(2)}}(-ux, u)}{u} \dd u,
\]
 where $\varphi_{V^{(1)},V^{(2)}}$ denotes the joint characteristic function of $(V^{(1)},V^{(2)})$ given in Theorem \ref{thm:partial_sums},
 and $\int_{-\infty}^\infty$ is meant in the sense of Cauchy principal value, i.e.,
 $\int_{-\infty}^\infty:= \lim_{T\to\infty}\lim_{h\to 0}\left( \int_h^T + \int_{-T}^{-h}\right)$.
\end{proposition}

\noindent{\bf Proof.}
By Proposition \ref{Lem_abs_cont}, $V^{(2)}/V^{(1)}$ is absolutely continuous, so the inversion formula for characteristic functions due to Gurland \cite{Gur}
 yields that for each $x\in\RR$,
 \begin{align*}
  \pr \left( \frac{V^{(2)}}{V^{(1)}} \leq x \right)
   = \pr(V^{(2)} - x V^{(1)} \leq 0)
    = \frac{1}{2} - \frac{1}{2\pi\ii} \int_{-\infty}^\infty \ee^{-\ii u 0}\, \frac{\varphi_{V^{(2)} - x V^{(1)}}(u)}{u} \,\dd u,
 \end{align*}
 yielding the statement, where $\varphi_{V^{(2)} - x V^{(1)}}$ denotes the characteristic function of $V^{(2)} - x V^{(1)}$.
\proofend

\appendix

\section{On topological properties of \ $S$}
\label{app:topology}

\begin{lemma}\label{Lem_S_top}
The set $S=(0,\infty)\times \R$ furnished with the metric $d$ given in \eqref{S_metric} is a complete separable metric space,
 and $B \subset S$ is bounded with respect to the metric $d$ if and only if $B$ is separated from
 the vertical line $\{(0,y): y\in\RR\}$, i.e., there exists $\vare>0$ such that $B \subset \{ (x,y)\in S : x>\vare\}$.
Moreover, the topology and the Borel $\sigma$-algebra $\cB(S)$ on $S$ induced by the metric $d$ coincides with
 the topology and the Borel $\sigma$-algebra on $S$ induced by the usual metric $\rho((x,x'), (y,y')) := \sqrt{(x - x')^2 + (y - y')^2}$, $(x,x'), (y,y') \in S$,
 respectively.
\end{lemma}

\noindent{\bf Proof.}
First, we check that $S$ is a complete separable metric space.
If $(x_n,y_n)_{n\in\NN}$ is a Cauchy sequence in $S$, then for all $\vare\in(0,1)$,
 there exists an $N_\vare\in\NN$ such that $d((x_n,y_n), (x_m,y_m))<\vare$ for $n,m\geq N_\vare$.
Hence $\rho((x_n,y_n), (x_m,y_m)) < \vare$ and $\left\vert \frac{1}{x_n} - \frac{1}{x_m}\right\vert<\vare$ for $n,m\geq N_\vare$,
 i.e., $(x_n,y_n)_{n\in\NN}$ and $(1/x_n)_{n\in\NN}$ are Cauchy sequences in $\RR^2$ and in $\RR$, respectively.
 Consequently, there exists an $(x,y)\in[0,\infty)\times \RR$ such that $\lim_{n\to\infty} (x_n,y_n) = (x,y)$
 and $\frac{1}{x_n}$ being convergent as $n\to\infty$, yielding that $x>0$, and so $(x,y)\in(0,\infty)\times \RR$.
By continuity, $\lim_{n\to\infty} d((x_n,y_n),(x,y))=0$, as desired.
The separability of $S$ readily follows, since $S\cap \QQ^2$ is a countable everywhere dense subset of $S$.

\noindent
Next, we check that $B \subset S$ is bounded with respect to the metric $d$ if and only if
 there exists $\vare>0$ such that $B \subset \{ (x,y)\in S : x>\vare\}$.
If $B \subset S$ is bounded, then there exists $r>0$ such that $d((x,y),(1,0))<r$, $(x,y)\in B$,
 yielding that $\vert \frac{1}{x} - 1\vert < r$, $(x,y)\in B$, and then $x>\frac{1}{1+r}$, $(x,y)\in B$,
 so one can choose $\vare=\frac{1}{1+r}$.
If there exists $\vare>0$ such that $B \subset \{ (x,y)\in S : x>\vare\}$, then
 $d((x,y),(1,0)) = \min(\sqrt{(x-1)^2 + y^2},1) + \vert \frac{1}{x} - 1\vert \leq 1+\frac{1}{\vare}+1$, $(x,y)\in B$.
\proofend

Since $S$ is locally compact, second countable and Hausdorff, one could choose a metric such that the relatively compact sets are precisely the bounded ones, see Kallenberg \cite[page 18]{kallenberg:2017}.
The metric $d$ does not have this property, but we do not need it.
For historical fidelity, we note that originally the vague convergence of point measures in ${\cal M}_p(S)$ is defined by the convergence of integrals of some compactly supported functions (see, e.g.,
 Resnick \cite[Section 3.3.5]{resnick:2007}), but recently instead of compactly supported functions one uses functions with bounded support
 (see, e.g., Kallenberg \cite{kallenberg:2017}).
We also follow the latter approach.

\section{Vague convergence of point measures}
\label{app:vague_convergence_of_point_measures}

\noindent{\bf Proof of Lemma \ref{lem:vagueConv_pointMeas}.}
First, let us suppose that $\mu_n\vto \mu$ as $n\to\infty$, and let $\epsilon>0$ be such that $\mu(\{\epsilon\}\times \R)=0$.
Since $\mu$ is locally finite and $(\epsilon,\infty)\times \R$ is bounded, there exist integers $K \geq 0$ and $c_1, \ldots, c_K \in \N$ such that
 \[
   \mu|_{(\epsilon,\infty)\times \R} = \sum_{j=1}^K c_j \delta_{(u_j, v_j)} ,
 \]
 where $(u_1, v_1), \ldots, (u_K, v_K)$ are the atoms of $\mu$ in $(\epsilon,\infty)\times \R$ and $c_1, \ldots, c_K$ are their multiplicities.
Let {$s_0 := 0$, $s_j:=c_1+\cdots+c_j$ for $j \in \{1,\ldots,K\}$ and $M := s_K$, and let us label the points of $\mu|_{(\epsilon,\infty)\times \R}$ such that for each $j \in \{1,\ldots,K\}$ we have $(u_j, v_j) = (x_k, y_k)$ for all $k \in \{s_{j-1}+1,\ldots,s_j\}$, yielding that
 $\mu|_{(\epsilon,\infty)\times \R}=\sum_{i=1}^M \delta_{(x_i, y_i)}$}.
Since $(\epsilon,\infty)\times \R$ is open in $S$ (see Lemma \ref{Lem_S_top}), one can choose pairwise disjoint open sets $G_1,\ldots,G_K\subset (\epsilon,\infty)\times \R$
 such that $(u_j, v_j)\in G_j$, $j=1,\ldots K$.
Especially, we have $\mu(\partial G_j)=0$, $j \in \{1,\ldots,K\}$, where $\partial G_j$ denotes the boundary of $G_j$ in $S$
 (since $\mu|_{(\epsilon,\infty)\times \R}$ puts zero mass outside the points $(u_1, v_1), \ldots, (u_K, v_K)$).
Hence, since $\mu_n\vto \mu$ as $n\to\infty$, using the equivalence of (i) and (iv) in Lemma 4.1 in Kallenberg \cite{kallenberg:2017}, we have
 $\mu_n(G_j)\to \mu(G_j)=c_j$ as $n\to\infty$ for each $j \in \{1,\ldots,K\}$.
Similarly,
 \[
   \mu_n( ((\epsilon,\infty)\times \R)\setminus (G_1\cup \cdots\cup G_K)) \to \mu( ((\epsilon,\infty)\times \R) \setminus (G_1\cup \cdots\cup G_K)) =0
    \qquad \text{as \ $n\to\infty$,}
 \]
 since $((\epsilon,\infty)\times \R)\setminus (G_1\cup \cdots\cup G_K)$ is a bounded (with respect to metric $d$ given in \eqref{S_metric}) Borel subset of $S$
  and $\mu(\partial( ((\epsilon,\infty)\times \R)\setminus (G_1\cup \cdots\cup G_K)))=0$ (using also the assumption $\mu(\{\epsilon\}\times \R)=0$).
Consequently, since $\mu_n$, $n\in\NN$, and $\mu$ are integer-valued measures, there exists an integer $n_0\geq 0$ such that
 $\mu_n(G_j)=c_j$, $j \in \{1,\ldots,K\}$, and $\mu_n(((\epsilon,\infty)\times \R)\setminus (G_1\cup \cdots\cup G_K)) =0$ for all $n\geq n_0$,
 yielding that $\mu_n((\epsilon,\infty)\times \R) = \mu((\epsilon,\infty)\times \R)=c_1+\cdots+c_K=M$ for all $n\geq n_0$.
So for each $n\geq n_0$, one can label the points of $\mu_n|_{(\epsilon,\infty)\times \R}$ such that
 for each $j \in \{1,\ldots,K\}$ we have $(x_k^{(n)}, y_k^{(n)}) \in G_j$ for all $k \in \{s_{j-1}+1,\ldots,s_j\}$, yielding that
 $\mu_n|_{(\epsilon,\infty)\times \R}=\sum_{i=1}^M \delta_{(x_i^{(n)}, y_i^{(n)})}$.
Shrinking the open sets $G_1,\ldots,G_K$ onto $(u_1,v_1),\ldots, (u_K,v_K)$, respectively, we have
 $x_i^{(n)}\to x_i$ and $y_i^{(n)}\to y_i$ as $n \toi$ for all $i=1,\dots,M$.

Now, let us prove the reverse direction.
Let us suppose that $f : S \to [0, \infty)$ is a bounded, continuous function with bounded support.
Then, using Lemma \ref{Lem_S_top}, there exists $\vare_0 \in(0,\infty)$ such that $f(x,y)=0$ for all $(x,y) \in (0,\vare_0]\times\R$.
Since the function $(0,\vare_0)\ni\vare\mapsto \mu((\vare,\infty)\times\R)\in[0,\infty)$ is decreasing,
 there exists $\vare\in(0,\vare_0)$ such that $\mu(\{\epsilon\}\times \R)=0$.
Due to our assumptions, there exist integers $n_0,M\geq 0$ and a labeling of the points of $\mu$ and $\mu_n$, $n \geq n_0$, in $(\epsilon,\infty)\times \R$ such that
\begin{align*}
\mu_n|_{(\epsilon,\infty)\times \R}=\sum_{i=1}^M \delta_{(x_i^{(n)}, \, y_i^{(n)})} \, , \qquad  \mu|_{(\epsilon,\infty)\times \R}=\sum_{i=1}^M \delta_{(x_i, y_i)} \, ,
\end{align*}
and $x_i^{(n)}\to x_i$ and $y_i^{(n)}\to y_i$ as $n \toi$ for all $i=1,\dots,M$.
Consequently,
 \begin{align*}
  \mu_n(f)
  &= \iint_{(\vare_0,\infty)\times\RR} f(x,y) \, \mu_n(\dd x, \dd y)
   = \iint_{(\vare,\infty)\times\RR} f(x,y) \, \mu_n(\dd x, \dd y)
   = \sum_{i=1}^M f(x_i^{(n)}, \, y_i^{(n)}) \\
  &\to \sum_{i=1}^M f(x_i, y_i)
   = \iint_{(\vare,\infty)\times\RR} f(x,y) \, \mu(\dd x, \dd y)
   = \iint_{(\vare_0,\infty)\times\RR} f(x,y) \, \mu(\dd x, \dd y)= \mu(f)
 \end{align*}
 as $n\to\infty$, hence we have $\mu_n\vto \mu$ as $n\to\infty$, as desired.
\proofend

\section{Approximation of Laplace functionals}
\label{app:laplace}

First, we recall an auxiliary lemma stating that $(X_i)_{i\geq 0}$ is strongly mixing with a geometric rate
 from Basrak et al.\ \cite[Remark 3.1]{basrak:kulik:palmowski:2013} and Barczy et al.\ \cite[Lemma F.1]{barczy:nedenyi:pap:2019}.

A strongly stationary sequence $(Y_k)_{k\geq 0}$ is called strongly mixing with a rate function $(\alpha_h)_{h\in\NN}$
 if its strongly stationary extensions $(Y_k)_{k\in\ZZ}$ admit this property, namely,
 \begin{align}\label{alpha}
   \alpha_h := \sup_{A\in{\mathcal F}^Y_{-\infty,0},\; B\in{\mathcal F}^Y_{h,\infty}} |\pr(A \cap B) - \pr(A) \pr(B)|
   \to 0 \qquad \text{as \ $h \to \infty$,}
 \end{align}
 where ${\mathcal F}^Y_{-\infty,0} := \sigma(\ldots, Y_{-1}, Y_0)$, ${\mathcal F}^{Y}_{h,\infty} := \sigma(Y_h, Y_{h+1}, \ldots)$, $h \in \NN$.

\begin{lemma}\label{lemma:strong_mixing}
The strongly stationary Markov chain $(X_i)_{i\geq 0}$ is strongly mixing with a geometric rate, i.e., there exists a constant
 $q\in(0,1)$ such that $\alpha_\ell = O(q^\ell)$ as $\ell\toi$.
\end{lemma}

Note that in this paper we need only that $(X_i)_{i\geq 0}$ is strongly mixing, and we will not use that the mixing rate is geometric,
 however, for completeness, we decided to recall it in Lemma \ref{lemma:strong_mixing} as well.

Next, we show that the process
 $(X_i \1{\{X_i>0\}}, \frac{M_{i+1}}{\sqrt{X_i}} \1{\{X_i>0\}} )_{i\geq0}$ satisfies a certain mixing condition
  (for the definition of $M_{i+1}$, $i\geq 0$, see the Introduction).

\begin{lemma}\label{lemma:mixing}
There exists a sequence of positive integers $(r_n)_{n\in\NN}$ with $r_n \to \infty$ and $r_n/n \to 0$ as $n \toi$ such that for each bounded, continuous function $f:S\to [0,\infty)$ having the property $f(x,y)=0$ for all $(x,y)\in (0,\vare]\times\RR$ for some $\epsilon>0$, we have
 \begin{equation*}
   \ex\biggl[\exp\biggl\{- \sideset{}{^{*}} \sum_{i=1}^n f\biggl(\frac{X_i}{a_n},  \frac{M_{i+1}}{\sqrt{X_i}} \biggr)\biggr\}\biggr]
   - \biggl(
     \ex\biggl[\exp\biggl\{- \sideset{}{^{*}} \sum_{i=1}^{r_n} f\biggl(\frac{X_i}{a_n}, \frac{M_{i+1}}{\sqrt{X_i}}  \biggr)\biggr\}\biggr]\biggr)^{k_n}
    \to 0
 \end{equation*}
 as $n \to \infty$ with $k_n := \lfloor n/r_n\rfloor$, where we recall that $S=(0,\infty)\times \R$ and
 \[
  \sideset{}{^{*}} \sum_{i=1}^m = \sum_{\{ j\in\{1,\ldots,m\} : X_j>0\}}\;\;, \qquad m\in\NN.
 \]
\end{lemma}

\noindent{\bf Proof.}
We follow the proof of Proposition 1.34 in Krizmani\'c \cite{Kri} (see also Basrak \cite[Lemma 2.3.9]{basrakPhD}).
Let $(\ell_n)_{n\in\N}$ be a sequence of positive integers with $\ell_n \toi$ and $\ell_n/n^{1/8} \to 0$ as $n \toi$.
We will show that the sequence
 \[
   r_n := \lfloor \max\{n \sqrt{\alpha_{\ell_n}}, n^{2/3}\}\rfloor + 1 , \qquad n \in \N ,
 \]
 is a good choice with $\alpha_\ell$, $\ell\in\NN$, given in \eqref{alpha}.
Clearly, $r_n \toi$ as $n \toi$.
By Lemma \ref{lemma:strong_mixing}, the strongly stationary Markov chain $(X_i)_{i\geq 0}$ is strongly mixing,
 i.e., $\alpha_{\ell_n} \to 0$ as $n \toi$, yielding $r_n/n \to 0$ as $n\toi$ and
 \begin{equation}\label{limits}
  k_n \toi , \qquad k_n \alpha_{\ell_n} \to 0 , \qquad \frac{k_n\ell_n}{n} \to 0
 \end{equation}
 as $n \toi$.

Fix a bounded, continuous function $f:S\to [0,\infty)$ having the property $f(x,y)=0$ for all $(x,y)\in (0,\vare]\times\RR$ for some $\epsilon>0$.
Put $M := \sup_{(x,y)\in S} f(x,y) < \infty$.
We have to show that $I(n) \to 0$ as $n \toi$ with
 \[
   I(n) :=
   \biggl|\ex\biggl[\exp\biggl\{- \sideset{}{^{*}} \sum_{i=1}^n f\biggl(\frac{X_i}{a_n},  \frac{M_{i+1}}{\sqrt{X_i}} \biggr)\biggr\}\biggr]
   - \biggl(
     \ex\biggl[\exp\biggl\{- \sideset{}{^{*}} \sum_{i=1}^{r_n} f\biggl(\frac{X_i}{a_n}, \frac{M_{i+1}}{\sqrt{X_i}} \biggr)\biggr\}\biggr]\biggr)^{k_n}\biggr| .
 \]
We have
 \[
   I(n) \leq I_1(n) + I_2(n) + I_3(n) + I_4(n) , \qquad n \in \N ,
 \]
 with
 \begin{align*}
  I_1(n) &:= \biggl|\ex\biggl[\exp\biggl\{- \sideset{}{^{*}} \sum_{i=1}^n f\biggl(\frac{X_i}{a_n}, \frac{M_{i+1}}{\sqrt{X_i}} \biggr)\biggr\}\biggr]
   - \ex\biggl[\exp\biggl\{-  \sideset{}{^{*}} \sum_{i=1}^{k_nr_n} f\biggl(\frac{X_i}{a_n}, \frac{M_{i+1}}{\sqrt{X_i}} \biggr)\biggr\}\biggr]\biggr| , \\
  I_2(n) &:= \biggl|\ex\biggl[\exp\biggl\{- \sideset{}{^{*}} \sum_{i=1}^{k_nr_n} f\biggl(\frac{X_i}{a_n}, \frac{M_{i+1}}{\sqrt{X_i}} \biggr)\biggr\}\biggr]
  - \ex\biggl[\exp\biggl\{-\sum_{k=1}^{k_n} \, \sideset{}{^{*}} \sum_{i=(k-1)r_n+1}^{kr_n-\ell_n} f\biggl(\frac{X_i}{a_n}, \frac{M_{i+1}}{\sqrt{X_i}} \biggr)\biggr\}\biggr]\biggr| , \\
  I_3(n) &:= \biggl|\ex\biggl[\exp\biggl\{-\sum_{k=1}^{k_n} \, \sideset{}{^{*}} \sum_{i=(k-1)r_n+1}^{kr_n-\ell_n} f\biggl(\frac{X_i}{a_n}, \frac{M_{i+1}}{\sqrt{X_i}} \biggr)\biggr\}\biggr]
  - \biggl(
     \ex\biggl[\exp\biggl\{- \sideset{}{^{*}} \sum_{i=1}^{r_n-\ell_n} f\biggl(\frac{X_i}{a_n}, \frac{M_{i+1}}{\sqrt{X_i}} \biggr)\biggr\}\biggr]\biggr)^{k_n}\biggr| , \\
  I_4(n) &:= \biggl|\biggl(
     \ex\biggl[\exp\biggl\{- \sideset{}{^{*}} \sum_{i=1}^{r_n-\ell_n} f\biggl(\frac{X_i}{a_n}, \frac{M_{i+1}}{\sqrt{X_i}} \biggr)\biggr\}\biggr]\biggr)^{k_n}
     - \biggl(
     \ex\biggl[\exp\biggl\{- \sideset{}{^{*}} \sum_{i=1}^{r_n} f\biggl(\frac{X_i}{a_n}, \frac{M_{i+1}}{\sqrt{X_i}} \biggr)\biggr\}\biggr]\biggr)^{k_n}\biggr| ,
 \end{align*}
 where, by \eqref{limits}, $kr_n - \ell_n\to\infty$ as $n\to\infty$ for each $k\in\N$.
By the strong stationarity of $(X_k)_{k\geq0}$ and using the inequality
$1 - \ee^{-x} \leq x$ for any $x \in (0, \infty)$, we obtain
 \begin{align*}
  I_1(n) &\leq \ex\biggl[\exp\biggl\{- \sideset{}{^{*}} \sum_{i=1}^{k_nr_n} f\biggl(\frac{X_i}{a_n}, \frac{M_{i+1}}{\sqrt{X_i}} \biggr)\biggr\}
  \biggl|1 - \exp\biggl\{-\sideset{}{^{*}} \sum_{i=k_nr_n+1}^n f\biggl(\frac{X_i}{a_n}, \frac{M_{i+1}}{\sqrt{X_i}} \biggr)\biggr\}\biggr|\biggr] \\
  &\leq  \ex\biggl[\;\sideset{}{^{*}} \sum_{i=k_nr_n+1}^n f\biggl(\frac{X_i}{a_n}, \frac{M_{i+1}}{\sqrt{X_i}}\biggr)\biggr]
   = \sideset{}{^{*}} \sum_{i=k_nr_n+1}^n \ex\biggl[f\biggl(\frac{X_i}{a_n}, \frac{M_{i+1}}{\sqrt{X_i}}\biggr) \1{ \bigl\{\frac{X_i}{a_n}>\vare\bigr\} }\biggr]\\
  &\leq \sum_{i=k_nr_n+1}^n M \pr(X_i > \vare a_n)
   = (n - k_n r_n) M \pr(X_0 > \vare a_n) .
 \end{align*}
In a similar manner we obtain
 \begin{align*}
  I_2(n) &\leq \ex\biggl[\exp\biggl\{-\sum_{k=1}^{k_n} \, \sideset{}{^{*}} \sum_{i=(k-1)r_n+1}^{kr_n-\ell_n} f\biggl(\frac{X_i}{a_n}, \frac{M_{i+1}}{\sqrt{X_i}}\biggr)\biggr\}
  \biggl|1 - \exp\biggl\{-\sum_{k=1}^{k_n} \, \sideset{}{^{*}} \sum_{i=kr_n-\ell_n+1}^{kr_n}
          f\biggl(\frac{X_i}{a_n}, \frac{M_{i+1}}{\sqrt{X_i}} \biggr)\biggr\}\biggr|\biggr] \\
  &\leq \ex\biggl[\sum_{k=1}^{k_n} \, \sideset{}{^{*}} \sum_{i=kr_n-\ell_n+1}^{kr_n} f\biggl(\frac{X_i}{a_n}, \frac{M_{i+1}}{\sqrt{X_i}} \biggr)\biggr]
   = \sum_{k=1}^{k_n} \, \sideset{}{^{*}} \sum_{i=kr_n-\ell_n+1}^{kr_n} \ex\biggl[f\biggl(\frac{X_i}{a_n},  \frac{M_{i+1}}{\sqrt{X_i}} \biggr) \1{ \{ \frac{X_i}{a_n}>\vare \} }\biggr]\\
  &\leq \sum_{k=1}^{k_n} \sum_{i=kr_n-\ell_n+1}^{kr_n} M \pr(X_i > \vare a_n)
   = k_n \ell_n M \pr(X_0 > \vare a_n) .
 \end{align*}
We have
 \[
   I_3(n) \leq I_5(n) + I_6(n) , \qquad n \in \N ,
 \]
 with
 \begin{align*}
  I_5(n) &:= \biggl|\ex\biggl[\exp\biggl\{-\sum_{k=1}^{k_n} \, \sideset{}{^{*}} \sum_{i=(k-1)r_n+1}^{kr_n-\ell_n} f\biggl(\frac{X_i}{a_n}, \frac{M_{i+1}}{\sqrt{X_i}} \biggr)\biggr\}\biggr] \\
  &\phantom{:= \biggl|}
  - \ex\biggl[\exp\biggl\{- \sideset{}{^{*}} \sum_{i=1}^{r_n-\ell_n} f\biggl(\frac{X_i}{a_n}, \frac{M_{i+1}}{\sqrt{X_i}} \biggr)\biggr\}\biggr]
    \ex\biggl[\exp\biggl\{-\sum_{k=2}^{k_n} \, \sideset{}{^{*}} \sum_{i=(k-1)r_n+1}^{kr_n-\ell_n} f\biggl(\frac{X_i}{a_n},  \frac{M_{i+1}}{\sqrt{X_i}} \biggr)\biggr\}\biggr]\biggr| , \\
  I_6(n) &:= \biggl|\ex\biggl[\exp\biggl\{- \sideset{}{^{*}} \sum_{i=1}^{r_n-\ell_n} f\biggl(\frac{X_i}{a_n}, \frac{M_{i+1}}{\sqrt{X_i}} \biggr)\biggr\}\biggr]
    \ex\biggl[\exp\biggl\{-\sum_{k=2}^{k_n} \, \sideset{}{^{*}} \sum_{i=(k-1)r_n+1}^{kr_n-\ell_n} f\biggl(\frac{X_i}{a_n},  \frac{M_{i+1}}{\sqrt{X_i}} \biggr)\biggr\}\biggr] \\
  &\phantom{:= \biggl|}
    - \biggl(
     \ex\biggl[\exp\biggl\{- \sideset{}{^{*}} \sum_{i=1}^{r_n-\ell_n} f\biggl(\frac{X_i}{a_n}, \frac{M_{i+1}}{\sqrt{X_i}} \biggr)\biggr\}\biggr]\biggr)^{k_n}\biggr| .
 \end{align*}
Since $(X_i)_{i\geq0}$ is strongly mixing, we have
 \[
   |\ex[\xi \eta] - \ex[\xi] \ex[\eta]| \leq 4 C_1 C_2 \alpha_m
 \]
 for any $\cF^X_{0,j}$-measurable random variable $\xi$ and any $\cF_{j+m,\infty}^X$-measurable random variable $\eta$ with $j, m \in \N$, $|\xi| \leq C_1$ and $|\eta| \leq C_2$ (see, e.g., Lemma 1.2.1 in Lin and Lu \cite{lin:lu:1996}).
Hence, using that the random variables
 \[
   \sideset{}{^{*}} \sum_{i=1}^{r_n-\ell_n} f\biggl(\frac{X_i}{a_n}, \frac{M_{i+1}}{\sqrt{X_i}} \biggr) \qquad \text{and} \qquad
   \sum_{k=2}^{k_n} \, \sideset{}{^{*}} \sum_{i=(k-1)r_n+1}^{kr_n-\ell_n} f\biggl(\frac{X_i}{a_n},  \frac{M_{i+1}}{\sqrt{X_i}} \biggr)
 \]
 are $\cF^X_{0,r_n-\ell_n+1}$-measurable and $\cF^X_{r_n+1,\infty}$-measurable, respectively, we have
 \[
   I_5(n) \leq 4 \alpha_{\ell_n} , \qquad n\in\NN.
 \]
It is easy to obtain that
 \begin{align*}
  &I_6(n) = \ex\biggl[\exp\biggl\{-\sideset{}{^{*}} \sum_{i=1}^{r_n-\ell_n} f\biggl(\frac{X_i}{a_n}, \frac{M_{i+1}}{\sqrt{X_i}} \biggr)\biggr\}\biggr] \\
  &\hspace*{10mm}
   \times \biggl|\ex\biggl[\exp\biggl\{-\sum_{k=2}^{k_n} \, \sideset{}{^{*}} \sum_{i=(k-1)r_n+1}^{kr_n-\ell_n} f\biggl(\frac{X_i}{a_n},  \frac{M_{i+1}}{\sqrt{X_i}} \biggr)\biggr\}\biggr]
  - \biggl(
     \ex\biggl[\exp\biggl\{- \sideset{}{^{*}} \sum_{i=1}^{r_n-\ell_n} f\biggl(\frac{X_i}{a_n}, \frac{M_{i+1}}{\sqrt{X_i}} \biggr)\biggr\}\biggr]\biggr)^{k_n-1}\biggr| \\
  &\leq \biggl|\ex\biggl[\exp\biggl\{-\sum_{k=2}^{k_n} \, \sideset{}{^{*}} \sum_{i=(k-1)r_n+1}^{kr_n-\ell_n} f\biggl(\frac{X_i}{a_n}, \frac{M_{i+1}}{\sqrt{X_i}} \biggr)\biggr\}\biggr]
  - \biggl(
     \ex\biggl[\exp\biggl\{-  \sideset{}{^{*}} \sum_{i=1}^{r_n-\ell_n} f\biggl(\frac{X_i}{a_n}, \frac{M_{i+1}}{\sqrt{X_i}} \biggr)\biggr\}\biggr]\biggr)^{k_n-1}\biggr| ,
 \end{align*}
 hence we recursively obtain (we repeat the same procedure for the above estimation of $I_6(n)$ as we did for $I_3(n)$ and so on)
 \[
   I_3(n) \leq 4 k_n \alpha_{\ell_n} .
 \]
Strong stationarity of $(X_i)_{i\geq0}$ and Lemma 4.3 in Chapter 2 in Durrett \cite{durrett:1996} imply
 \begin{align*}
  I_4(n) &\leq k_n \biggl|\ex\biggl[\exp\biggl\{-\sideset{}{^{*}} \sum_{i=1}^{r_n-\ell_n} f\biggl(\frac{X_i}{a_n}, \frac{M_{i+1}}{\sqrt{X_i}} \biggr)\biggr\}\biggr]
     - \ex\biggl[\exp\biggl\{-\sideset{}{^{*}} \sum_{i=1}^{r_n} f\biggl(\frac{X_i}{a_n},  \frac{M_{i+1}}{\sqrt{X_i}} \biggr)\biggr\}\biggr]\biggr| \\
  &\leq k_n \ex\biggl[\exp\biggl\{-\sideset{}{^{*}} \sum_{i=1}^{r_n-\ell_n} f\biggl(\frac{X_i}{a_n}, \frac{M_{i+1}}{\sqrt{X_i}} \biggr)\biggr\}
        \biggl|1 - \exp\biggl\{-\sideset{}{^{*}} \sum_{i=r_n-\ell_n+1}^{r_n} f\biggl(\frac{X_i}{a_n},  \frac{M_{i+1}}{\sqrt{X_i}} \biggr)\biggr\}
        \biggr|\biggr] \\
  &\leq k_n \ex\biggl[  \; \sideset{}{^{*}} \sum_{i=r_n-\ell_n+1}^{r_n}  f\biggl(\frac{X_i}{a_n}, \frac{M_{i+1}}{\sqrt{X_i}} \biggr)\biggr]
   = k_n \sideset{}{^{*}} \sum_{i=r_n-\ell_n+1}^{r_n} \ex\biggl[f\biggl(\frac{X_i}{a_n}, \frac{M_{i+1}}{\sqrt{X_i}} \biggr) \1{\{ \frac{X_i}{a_n}>\vare \} }\biggr] \\
  &\leq k_n \sum_{i=r_n-\ell_n+1}^{r_n} M \pr(X_i > \vare a_n)
   = k_n \ell_n M \pr(X_0 > \vare a_n) .
 \end{align*}
Since $X_0$ is regularly varying with tail index $\alpha$, by \eqref{a_n_quantile} and \eqref{limits}, we obtain
 \begin{align*}
  I(n) &\leq (n - k_n r_n + 2 k_n \ell_n) M \pr(X_0 > \vare a_n) + 4 k_n \alpha_{\ell_n} \\
  &= \frac{(n-k_nr_n+2k_n\ell_n)M}{n} \cdot n \pr(X_0 > a_n) \cdot
     \frac{\pr(X_0>\vare a_n)}{\pr(X_0>a_n)} + 4 k_n \alpha_{\ell_n}
   \to 0
 \end{align*}
 as $n \toi$, since
 \[
   \frac{n-k_nr_n}{n}\leq \frac{n-\left(\frac{n}{r_n}-1\right)r_n}{n}
                     = \frac{r_n}{n}\to0
                     \qquad \text{as $n\to\infty$,}
 \]
 as desired.
\proofend

\section{Conditional Slutsky's lemma, conditional continuous mapping theorem}
\label{app:slutsky+CMT}

First, we prove the analogues of parts (iv) and (v) of Theorem 2.7 in van der Vaart \cite{vaart:1998} for probability measures instead of random vectors.

\begin{lemma}\label{lemma:Slutsky1}
Let $\mu_n$, $n \in \Nset$, be probability measures on $(\Rset^{2k}, \cB(\Rset^{2k}))$ with some $k \in \Nset$.
For each $n \in \Nset$, consider the marginal probability measures $\mu_n^{(1)}$ and $\mu_n^{(2)}$ on $(\Rset^k, \cB(\Rset^k))$ defined by $\mu_n^{(1)}(B) := \mu_n(B \times \Rset^k)$ and $\mu_n^{(2)}(B) := \mu_n(\Rset^k \times B)$ for $B \in \cB(\Rset^k)$.
If $\mu_n^{(1)} \wto \mu^{(1)}$ as $n \to \infty$ with some probability measure $\mu^{(1)}$ on $(\Rset^k, \cB(\Rset^k))$ and $\mu_n(\{(\bx, \by) \in \Rset^k \times \Rset^k : \|\bx - \by\| > \vare\}) \to 0$ as $n \to \infty$ for all $\vare \in (0, \infty)$, then $\mu_n^{(2)} \wto \mu^{(1)}$ as $n \to \infty$.
\end{lemma}

\noindent{\bf Proof.}
For each bounded Lipschitz function $g : \Rset^k \to \Rset$ and for each $n \in \Nset$, we have
 \[
   \Delta_n^{(g)} := \biggl|\int_{\Rset^k} g(\by) \, \mu_n^{(2)}(\dd\by) - \int_{\Rset^k} g(\bx) \, \mu^{(1)}(\dd\bx)\biggr| \leq I_n^{(g)} + J_n^{(g)} ,
 \]
 where
 \begin{gather*}
  I_n^{(g)} := \biggl|\int_{\Rset^k} g(\by) \, \mu_n^{(2)}(\dd\by) - \int_{\Rset^k} g(\bx) \, \mu_n^{(1)}(\dd\bx)\biggr| , \\
  J_n^{(g)} := \biggl|\int_{\Rset^k} g(\bx) \, \mu_n^{(1)}(\dd\bx) - \int_{\Rset^k} g(\bx) \, \mu^{(1)}(\dd\bx)\biggr| .
 \end{gather*}
By the portmanteau lemma (see, e.g., van der Vaart \cite[Lemma 2.2]{vaart:1998}), the convergence $\mu_n^{(1)} \wto \mu^{(1)}$ as $n \to \infty$
 implies $J_n^{(g)} \to 0$ as $n \to \infty$.
Moreover, for each $\vare \in (0, \infty)$, by Fubini's theorem, we have
 \begin{align*}
  I_n^{(g)} &= \biggl|\int_{\Rset^{2k}} g(\by) \, \mu_n(\dd\bx, \dd\by) - \int_{\Rset^{2k}} g(\bx) \, \mu_n(\dd\bx, \dd\by)\biggr| \\
  &\leq \int_{\Rset^{2k}} |g(\by) - g(\bx)| \, \mu_n(\dd\bx, \dd\by) \\
  &= \int_{\|\bx-\by\|\leq\vare} |g(\bx) - g(\by)| \, \mu_n(\dd\bx, \dd\by) + \int_{\|\bx-\by\|>\vare} |g(\bx) - g(\by)| \, \mu_n(\dd\bx, \dd\by) \\
  &\leq \vare \sup_{\bx,\by\in\Rset^k,\,\bx\ne\by} \frac{|g(\bx)-g(\by)|}{\|\bx-\by\|}
   + 2 \mu_n(\{(\bx, \by) \in \Rset^k \times \Rset^k : \|\bx - \by\| > \vare\}) \sup_{\bx\in\Rset^k} |g(\bx)| .
 \end{align*}
By the assumptions, for each $\vare \in (0, \infty)$, we get
 \[
   \limsup_{n\to\infty} I_n^{(g)} \leq \vare \sup_{\bx,\by\in\Rset^k,\,\bx\ne\by} \frac{|g(\bx)-g(\by)|}{\|\bx-\by\|} ,
 \]
 thus $\limsup_{n\to\infty} I_n^{(g)} = 0$, and hence $\lim_{n\to\infty} I_n^{(g)} = 0$.
Consequently, for each bounded Lipschitz function $g : \Rset^k \to \Rset$, we obtain $\Delta_n^{(g)} \to 0$ as $n \to \infty$.
By the portmanteau lemma (see, e.g., van der Vaart \cite[Lemma 2.2]{vaart:1998}), we conclude $\mu_n^{(2)} \wto \mu^{(1)}$ as $n \to \infty$.
\proofend

\begin{lemma}\label{lemma:Slutsky2}
Let $\mu_n$, $n \in \Nset$, be probability measures on $(\Rset^{k+\ell}, \cB(\Rset^{k+\ell}))$ with some $k, \ell \in \Nset$.
For each $n \in \Nset$, consider the marginal probability measures \ $\mu_n^{(1)}$ and $\mu_n^{(2)}$ on $(\Rset^k, \cB(\Rset^k))$ and $(\Rset^\ell, \cB(\Rset^\ell))$, respectively, defined by $\mu_n^{(1)}(B_1) := \mu_n(B_1 \times \Rset^\ell)$ for $B_1 \in \cB(\Rset^k)$ and $\mu_n^{(2)}(B_2) := \mu_n(\Rset^k \times B_2)$ for $B_2 \in \cB(\Rset^\ell)$.
If $\mu_n^{(1)} \wto \mu^{(1)}$ as $n \to \infty$ with some probability measure $\mu^{(1)}$ on $(\Rset^k, \cB(\Rset^k))$ and $\mu_n^{(2)} \wto \delta_\bc$ as $n \to \infty$ with some $\bc \in \Rset^\ell$, then $\mu_n \wto \mu^{(1)} \times \delta_\bc$ as $n \to \infty$.
\end{lemma}

\noindent{\bf Proof.}
For each $n \in \Nset$, consider the probability measure $\tmu_n$ on $(\Rset^{k+\ell}\times \Rset^{k+\ell}, \cB(\Rset^{k+\ell}\times \Rset^{k+\ell}))$ defined by
 \[
   \tmu_n(H) := \mu_n(\{(\bx, \by) \in \Rset^k \times \Rset^\ell : (\bx, \by, \bx, \bc) \in H\}) , \qquad H \in \cB(\Rset^{k+\ell}\times \Rset^{k+\ell}) .
 \]
For each $n \in \Nset$, consider the marginal probability measures $\tmu_n^{(1)}$ and $\tmu_n^{(2)}$ on $(\Rset^{k+\ell}, \cB(\Rset^{k+\ell}))$ defined by $\tmu_n^{(1)}(A) := \tmu_n(A \times \Rset^{k+\ell})$ and $\tmu_n^{(2)}(A) := \tmu_n(\Rset^{k+\ell} \times A)$ for $A \in \cB(\Rset^{k+\ell})$.
Note that for each $n \in \Nset$ and $A \in \cB(\Rset^{k+\ell})$, we have
 \[
   \tmu_n^{(1)}(A) = \mu_n(\{(\bx, \by) \in \Rset^k \times \Rset^\ell : (\bx, \by) \in A\}) = \mu_n(A) ,
 \]
 hence $\tmu_n^{(1)} = \mu_n$.
Moreover, for each $n \in \Nset$ and $A \in \cB(\Rset^{k+\ell})$, we have
 \begin{align*}
  \tmu_n^{(2)}(A)
  &= \mu_n(\{(\bx, \by) \in \Rset^k \times \Rset^\ell : (\bx, \bc) \in A\}) \\
  &= \mu_n(\{\bx \in \Rset^k : (\bx, \bc) \in A\} \times \Rset^\ell) \\
  &= \int_{\Rset^\ell} \mu_n(\{\bx \in \Rset^k : (\bx, \by) \in A\} \times \Rset^\ell) \, \delta_\bc(\dd\by) \\
  &= \int_{\Rset^\ell} \mu_n^{(1)}(\{\bx \in \Rset^k : (\bx, \by) \in A\}) \, \delta_\bc(\dd\by) \\
  &= \int_{\Rset^\ell} \biggl(\int_{\Rset^k} \1{A}(\bx, \by) \, \mu_n^{(1)}(\dd\bx)\biggr) \delta_\bc(\dd\by)
   = (\mu_n^{(1)} \times \delta_\bc)(A) ,
 \end{align*}
 hence $\tmu_n^{(2)} = \mu_n^{(1)} \times \delta_\bc$.
Further, for each $\vare \in (0, \infty)$, we have
 \begin{align*}
  &\tmu_n(\{((\bx, \by), (\bu, \bv)) \in \Rset^{k+\ell} \times \Rset^{k+\ell} : \|(\bx, \by) - (\bu, \bv)\| > \vare\}) \\
  &= \mu_n(\{(\bx, \by) \in \Rset^k \times \Rset^\ell : \|(\bx, \by) - (\bx, \bc)\| > \vare\}) \\
  &= \mu_n(\{(\bx, \by) \in \Rset^k \times \Rset^\ell : \|\by - \bc\| > \vare\}) \\
  &= \mu_n^{(2)}(\{\by \in \Rset^\ell : \|\by - \bc\| > \vare\})
   \to 0 \qquad \text{as \ $n \to \infty$,}
 \end{align*}
 since $\mu_n^{(2)} \wto \delta_\bc$ as $n \to \infty$.
Thus, according to Lemma \ref{lemma:Slutsky1}, to prove the statement
 it suffices to show that $\mu_n^{(1)} \times \delta_\bc \wto \mu^{(1)} \times \delta_\bc$ as $n \to \infty$.
For every continuous, bounded function $g : \Rset^{k+\ell} \times \Rset^{k+\ell} \to \Rset$, by the portmanteau lemma (see, e.g., van der Vaart \cite[Lemma 2.2]{vaart:1998}), we have
 \begin{align*}
  &\int_{\Rset^{k+\ell}} g(\bx, \bz) \, (\mu_n^{(1)} \times \delta_\bc)(\dd\bx, \dd\bz)
   = \int_{\Rset^k} g(\bx, \bc) \, \mu_n^{(1)}(\dd\bx) \\
  &\to \int_{\Rset^k} g(\bx, \bc) \, \mu^{(1)}(\dd\bx)
   = \int_{\Rset^{k+\ell}} g(\bx, \bz) \, (\mu^{(1)} \times \delta_\bc)(\dd\bx, \dd\bz)
 \end{align*}
 as $n \to \infty$, since $\mu_n^{(1)} \wto \mu^{(1)}$ as $n \to \infty$, and the function \ $\Rset^k \ni \bx \mapsto g(\bx, \bc) \in \Rset$ is a continuous, bounded function.
Again by the portmanteau lemma, we conclude $\mu_n^{(1)} \times \delta_\bc \wto \mu^{(1)} \times \delta_\bc$ as $n \to \infty$, as desired.
\proofend

We will use the following conditional continuous mapping theorem and a consequence of it.
Recall that for a random vector $\bX$ and an event $A \in \cA$ such that $\pr(A) > 0$, the conditional
 law of $\bX$ with respect to $A$ is denoted by $\cL(\bX | A)$.

\begin{lemma}\label{lemma:cmt}
For each $n \in \Nset$, let $A_n \in \cA$ such that $\pr(A_n) > 0$.
Let $\bX$ and $\bX_n$, $n \in \Nset$, be $\R^k$-valued random vectors and let $h: \Rset^k \to \Rset^\ell$
 be a Borel measurable function with some $k, \ell \in \Nset$.
Suppose that $\cL(\bX_n | A_n) \wto \cL(\bX)$ as $n \to \infty$ and $\pr(\bX \in D_h) = 0$, where $D_h$ denotes the set of discontinuities of $h$.
\begin{enumerate}
 \item[\textup{(i)}]
  Then $\cL(h(\bX_n) | A_n) \wto \cL(h(\bX))$ as $n \to \infty$.
 \item[\textup{(ii)}]
  If, in addition, $h$ is bounded, then $\ex[h(\bX_n) | A_n] \to \ex[h(\bX)]$ as $n \to \infty$.
\end{enumerate}
\end{lemma}

\noindent{\bf Proof.}
(i).
For each $B \in \cB(\RR^k)$, let $\mu(B) := \pr(\bX \in B)$ and $\mu_n(B) := \pr(\bX_n \in B | A_n)$, $n \in \N$.
By the assumption, $\mu_n \wto \mu$ as $n \to \infty$.
By Billingsley \cite[Lemma 5.1]{billingsley:1968}, we obtain $\nu_n \wto \nu$ as $n \to \infty$, where the probability measures $\nu_n$, $n \in \Nset$, and $\nu$ on $(\Rset^\ell, \cB(\Rset^\ell))$ are defined by $\nu_n(B) := \mu_n(h^{-1}(B)) = \pr(\bX_n \in h^{-1}(B) | A_n) = \pr(h(\bX_n) \in B | A_n)$, $n \in \Nset$, and $\nu(B) := \mu(h^{-1}(B)) = \pr(\bX \in h^{-1}(B)) = \pr(h(\bX) \in B)$ for $B \in \cB(\Rset^\ell)$.
Consequently, we obtain $\cL(h(\bX_n) | A_n) \wto \cL(h(\bX))$ as $n \to \infty$, as desired.

(ii).
By Billingsley \cite[part (iii) of Lemma 5.2]{billingsley:1968}, we obtain $\int_{\RR^k} h(\bx) \, \mu_n(\dd\bx) \to \int_{\RR^k} h(\bx) \, \mu(\dd\bx)$ as $n \to \infty$.
Consequently, we obtain $\ex[h(\bX_n) | A_n] \to \ex[h(\bX)]$ as $n \to \infty$, as desired.
\proofend

Next, we prove a conditional analogue of part (v) of Theorem 2.7 in van der Vaart \cite{vaart:1998} together with one of its useful consequences.

\begin{lemma}\label{lemma:Slutsky4}
For each $n \in \Nset$, let $A_n \in \cA$ such that $\pr(A_n) > 0$.
Let $\bX$ and $\bX_n$, $n \in \Nset$, be $\Rset^k$-valued random vectors and let $\bY_n$, $n \in \Nset$, be $\Rset^\ell$-valued random vectors with some $k, \ell \in \N$.
Suppose that $\cL(\bX_n | A_n) \wto \cL(\bX)$ and $\cL(\bY_n | A_n) \wto \delta_\bc$ as $n \to\infty$ with some $\bc \in \Rset^\ell$.
\begin{enumerate}
 \item[\textup{(i)}]
  Then $\cL((\bX_n, \bY_n) | A_n) \wto \cL(\bX) \times \delta_\bc = \cL((\bX, \bc))$ as $n \to \infty$.
 \item[\textup{(ii)}]
  If, in addition, $h: \Rset^{k+\ell} \to \Rset^m$ is a Borel measurable function with some $m \in \Nset$ such that $h$ is continuous at every $(\bx_0, \bc)$, $\bx_0 \in \R^k$, then $\cL(h(\bX_n, \bY_n) | A_n) \wto \cL(h(\bX, \bc))$ as $n \to \infty$.
\end{enumerate}
\end{lemma}

\noindent{\bf Proof.}
(i).
We apply Lemma \ref{lemma:Slutsky2} for the probability measures $\mu_n := \cL((\bX_n, \bY_n) | A_n)$, $n \in \Nset$, on $(\Rset^{k+\ell}, \cB(\Rset^{k+\ell}))$.
Then we have
$\mu_n^{(1)} = \cL(\bX_n | A_n) \wto \cL(\bX)$ and $\mu_n^{(2)} = \cL(\bY_n | A_n) \wto \delta_\bc$ as $n \to \infty$, hence we obtain $\cL((\bX_n, \bY_n) | A_n) = \mu_n \wto \cL(\bX) \times \delta_\bc = \cL((\bX, \bc))$ as $n \to \infty$.

(ii).
By the assumption, $D_h \subset \Rset^k \times (\R^\ell \setminus \{\bc\})$, hence $\pr((\bX, \bc) \in D_h) = 0$.
Consequently, part (i) of Lemma \ref{lemma:cmt} implies $\cL(h(\bX_n, \bY_n) | A_n) \wto \cL(h(\bX, \bc))$ as $n \to \infty$.
\proofend

Finally, we provide a conditional Slutsky's lemma.

\begin{lemma}\label{lemma:Slutsky3}
For each $n \in \Nset$, let $A_n \in \cA$ such that $\pr(A_n) > 0$.
Let $\bX$ and $\bX_n$, $n \in \Nset$, be $\Rset^{k\times\ell}$-valued random matrices such that
 $\cL(\bX_n | A_n) \wto \cL(\bX)$ as $n \toi$ with some $k, \ell \in \N$.
\begin{enumerate}
 \item[\textup{(i)}]
  If $\bY_n$, $n \in \Nset$, are $\Rset^{k\times\ell}$-valued random matrices such that $\cL(\bY_n | A_n) \wto \delta_\bC$ as $n \toi$ with some $\bC \in \Rset^{k\times\ell}$, then $\cL(\bX_n + \bY_n | A_n) \wto \cL(\bX + \bC)$ as $n \to \infty$.
 \item[\textup{(ii)}]
  If $\bY_n$, $n \in \Nset$, are $\Rset^{m\times k}$-valued random matrices such that $\cL(\bY_n | A_n) \wto \delta_\bC$ as $n \toi$ with some $\bC \in \Rset^{m\times k}$
   and $m\in\NN$, then $\cL(\bY_n \bX_n | A_n) \wto \cL(\bC \bX)$ as $n \to \infty$.
 \item[\textup{(iii)}]
  If $\bY_n$, $n \in \Nset$, are $\Rset^{k\times k}$-valued random matrices such that $\cL(\bY_n | A_n) \wto \delta_\bC$ as $n \toi$ with some invertible $\bC \in \Rset^{k\times k}$, then $\cL(\bY_n^{-1} \bX_n | A_n) \wto \cL(\bC^{-1} \bX)$ as $n \to \infty$.
\end{enumerate}
\end{lemma}

\noindent{\bf Proof.}
Identifying $\Rset^{k\times\ell}$, $\Rset^{m\times k}$ and $\Rset^{m\times \ell}$ with $\Rset^{k\ell}$, $\Rset^{mk}$ and $\Rset^{m\ell}$, respectively, in a natural way, we can apply part (ii) of Lemma \ref{lemma:Slutsky4} for the Borel measurable functions
 \begin{gather*}
  \Rset^{k\times\ell} \times \Rset^{k\times\ell} \ni (\bU, \bV) \mapsto \bU + \bV \in \Rset^{k\times\ell} , \qquad
   \Rset^{k\times\ell} \times \Rset^{m\times k} \ni (\bU, \bV) \mapsto \bV\bU \in \Rset^{m\times\ell} , \\
   \Rset^{k\times\ell} \times \Rset^{k\times k} \ni (\bU, \bV) \mapsto
  \begin{cases}
   \bV^{-1} \bU \in \Rset^{k\times\ell} , &\text{if $\bU$ is invertible,} \\
   \bzero \in \Rset^{k\times\ell} , &\text{otherwise,}
  \end{cases}
 \end{gather*}
 and we obtain the statements.
\proofend

\section{Regular variation of a related process}
\label{app:regvarrelproc}

First, we recall Karamata's theorem for truncated moments, see, e.g., Bingham et al.\ \cite[pages 26-27]{bingham:goldie:teugels:1987}
 or Buraczewski et al.\ \cite[Appendix B.4]{BurDamMik}.

\begin{lemma}[Karamata's theorem for truncated moments]\label{truncated_moments}
Consider a non-negative regularly varying random variable \ $X$ \ with tail index \ $\alpha>0$.
\ Then
 \begin{align*}
  \lim_{x\to\infty}
   \frac{x^\beta\pr(X>x)}{\ex(X^\beta\1{\{X\leq x\}})}
  &= \frac{\beta-\alpha}{\alpha}, \qquad \text{if \ $\beta \in [\alpha, \infty)$,} \\
  \lim_{x\to\infty}
   \frac{x^\beta\pr(X>x)}{\ex(X^\beta\1{\{X>x\}})}
  &= \frac{\alpha-\beta}{\alpha}, \qquad \text{if \ $\beta \in (-\infty, \alpha)$.}
 \end{align*}
\end{lemma}

Now, we give a representation of the strongly stationary Markov chain $(X_i)_{i\in\Zset}$.

\begin{lemma}\label{lem:reprX}
We have
 \begin{equation}\label{repr}
  (X_k)_{k\in\Zset}
  \eind \biggl(B_k + \sum_{i=1}^\infty \theta_k^{(k-i)} \circ \cdots \circ \theta_{k-i+1}^{(k-i)} \circ B_{k-i}\biggr)_{k\in\Zset} ,
 \end{equation}
 where $\{B_k : k \in \Zset\}$ are independent random variables with the same distribution as $B$, and $\theta_k^{(\ell)}$, $k, \ell \in \Zset$, are given by
 \[
   \theta_k^{(\ell)} \circ i
    := \begin{cases}
       \sum_{j=1}^i A_{k,j}^{(\ell)} , & \text{if \ $i \in \Nset$,} \\
       0 , & \text{if \ $i = 0$,}
      \end{cases}
 \]
 where $A_{k,j}^{(\ell)}$, $j \in \Nset$, $k, \ell \in \Zset$, have the same distribution as $A$, and $\{B_k : k \in \Zset\}$ and $\theta_k^{(\ell)}$, $k, \ell \in \Zset$,
 are independent in the sense that the families $\{B_k : k \in \Zset\}$ and $\{A_{k,j}^{(\ell)} : j \in \Nset\}$, $k, \ell \in \Zset$,
 occurring in $\theta_k^{(\ell)}$, $k, \ell \in \Zset$,
 are independent families of independent random variables, and the series in the representation \eqref{repr} converge with probability one.
\end{lemma}

\noindent{\bf Proof.}
Due to Basrak  et al.~\cite[Lemma 2.2.1]{basrak:kulik:palmowski:2013}, the series in the representation \eqref{repr} converge with probability one.
Clearly, for each $k \in \Zset$ and $\ell  \in \Nset$, we have
 \begin{align}\label{help15}
 \begin{split}
  &\biggl(B_k + \sum_{i=1}^\infty \theta_k^{(k-i)} \circ \cdots \circ \theta_{k-i+1}^{(k-i)} \circ B_{k-i}, \ldots,
          B_{k+\ell} + \sum_{i=1}^\infty \theta_{k+\ell}^{(k+\ell-i)} \circ \cdots \circ \theta_{k+\ell-i+1}^{(k+\ell-i)} \circ B_{k+\ell-i}\biggr) \\
  &\eind
   \biggl(B_0 + \sum_{i=1}^\infty \theta_0^{(-i)} \circ \cdots \circ \theta_{-i+1}^{(-i)} \circ B_{-i}, \ldots, B_\ell + \sum_{i=1}^\infty \theta_\ell^{(\ell-i)} \circ \cdots \circ \theta_{\ell-i+1}^{(\ell-i)} \circ B_{\ell-i}\biggr) .
 \end{split}
 \end{align}
Indeed,
 since $\{B_i : i \in \Zset\}$ are identically distributed, independent of $\{\theta_i^{(j)} : i, j \in \Zset\}$,
 and for each $N \in \Nset$ and $i_1, \ldots, i_N \geq 0$,
 the distribution of the random vector $\bigl(\theta_k^{(k-1)} \circ i_1, \ldots, \theta_k^{(k-N)} \circ \cdots \circ \theta_{k-N+1}^{(k-N)} \circ i_N\bigr)$
 is invariant with respect to a shift of $k \in \Zset$,
 we have
 \begin{align*}
  &\biggl(B_k + \sum_{i=1}^n \theta_k^{(k-i)} \circ \cdots \circ \theta_{k-i+1}^{(k-i)} \circ B_{k-i}, \ldots,
          B_{k+\ell} + \sum_{i=1}^n \theta_{k+\ell}^{(k+\ell-i)} \circ \cdots \circ \theta_{k+\ell-i+1}^{(k+\ell-i)} \circ B_{k+\ell-i}\biggr)_{n\in \NN} \\
  &\eind
   \biggl(B_0 + \sum_{i=1}^n \theta_0^{(-i)} \circ \cdots \circ \theta_{-i+1}^{(-i)} \circ B_{-i}, \ldots, B_\ell
           + \sum_{i=1}^n \theta_\ell^{(\ell-i)} \circ \cdots \circ \theta_{\ell-i+1}^{(\ell-i)} \circ B_{\ell-i}\biggr)_{n\in \NN}.
 \end{align*}
Thus using that for a sequence of random variables $\xi_i$, $i\in\NN$, the series $\sum_{i=1}^\infty\xi_i$ is convergent almost surely
 if and only if $\pr(\sup_{m\in\NN} \Big\vert \sum_{i=n}^{n+m} \xi_i\Big\vert >\epsilon)\to 0$ as $n\to\infty$ for each $\epsilon>0$,
 the almost sure convergence of the series on the left and right hand sides of \eqref{help15} yields \eqref{help15}
 (for a similar argument, see Step 2 of the proof of Theorem \ref{thm:partial_sums}).
Hence the right hand side of \eqref{repr} defines a strongly stationary process.
Moreover, for each $k \in \Zset$, we have
 \begin{align*}
  &B_k + \sum_{i=1}^\infty \theta_k^{(k-i)} \circ \cdots \circ \theta_{k-i+1}^{(k-i)} \circ B_{k-i}\\
  &= B_k + \theta_k^{(k-1)} \circ B_{k-1} + \theta_k^{(k-2)} \circ \theta_{k-1}^{(k-2)} \circ B_{k-2}
         + \theta_k^{(k-3)} \circ \theta_{k-1}^{(k-3)} \circ \theta_{k-2}^{(k-3)} \circ B_{k-2} + \cdots \\
  &\eind B_k + \theta_k^{(k-1)} \circ \bigl(B_{k-1} + \theta_{k-1}^{(k-2)} \circ B_{k-2} + \theta_{k-1}^{(k-3)} \circ \theta_{k-2}^{(k-3)} \circ B_{k-3} + \cdots\bigr) \\
  &= B_k + \theta_k^{(k-1)} \circ \biggl(B_{k-1} + \sum_{i=1}^\infty\theta_{k-1}^{(k-i-1)} \circ \cdots \circ \theta_{k-i}^{(k-i-1)} \circ B_{k-1-i}\biggr) ,
 \end{align*}
 since $\theta_k^{(k-i)}$, $i \in \Nset$, are independent, and independent of
$\{B_{k-i} : i \geq 0\} \cup \{\theta_{k-\ell}^{(k-i-1)} : \ i, \ell \in \Nset\}$.
Consequently, the stochastic process given on the right hand side of \eqref{repr} is a time homogeneous Markov process with the same transition probabilities
 as the Galton--Watson process $(X_k)_{k\in\Z}$ \ with immigration satisfying \eqref{GWI} such that the distribution of $X_0$ is the unique stationary distribution
 of $(X_k)_{k\in \ZZ}$.
\proofend

It turns out that the process $(X_i^{3/2}, X_i M_{i+1})_{i\in\Zset}$ is regularly varying with
 tail index $\frac{2\alpha}{3}$ with an explicitly given forward tail process.

\begin{theorem}\label{th:joint_reg2}
As $x \to \infty$,
 \[
   \cL\biggl(\biggl(\frac{1}{x}
   (X_k^{3/2}, X_k M_{k+1})\biggr)_{k\geq0} \,\bigg|\, (X_0^{3/2} \lor X_0 |M_1|) > x\biggr)
   \fidi
   \cL\bigl((\mu_A^{3k/2} \tY, \mu_A^{3k/2} \tY \tZ_k)_{k\geq0}\bigr) ,
 \]
 where $\tZ_k$, $k\in\N$, is an i.i.d. sequence of $\cN(0, \sigma_A^2)$-distributed random variables,
 the distribution of $(\tY, \tZ_0)$ is given by
 \begin{align}\label{help6}
  \pr(\tY > y, \tZ_0 > v_0)
  = \frac{\ex\bigl(\bigl(y\lor(1\lor|Z_0|)^{-1}\bigr)^{-2\alpha/3}\1{(v_0,\infty)}(Z_0)\bigr)}
         {\ex\bigl((1\lor|Z_0|)^{2\alpha/3}\bigr)}
 \end{align}
 for $y, v_0 \in \Rset$, where $Z_0 \eind \cN(0, \sigma_A^2)$, and the random vector $(\tY, \tZ_0)$ \ is independent from the variables $\tZ_k$, $k \in \Nset$.
Consequently, the strongly stationary process $(X_k^{3/2}, X_k M_{k+1})_{k\in\Z}$ is jointly regularly varying with tail index $\frac{2\alpha}{3}$, i.e., all its finite dimensional distributions are regularly varying with tail index $\frac{2\alpha}{3}$.
The process
 \[
   \bigl(\mu_A^{3k/2} \tY, \mu_A^{3k/2} \tY \tZ_k\bigr)_{k\geq0}
 \]
 is the forward tail process of $(X_k^{3/2}, X_k M_{k+1})_{k\in\Z}$.
Moreover, there exists a (whole) tail process of $(X_k^{3/2}, X_k M_{k+1})_{k\in\ZZ}$ as well.
\end{theorem}

\noindent{\bf Proof.}
By Lemma \ref{lem:reprX}, we may and do suppose that $(X_k)_{k\geq0}$ is the right hand side of \eqref{repr}.
First, we give a useful representation of the random vectors $(X_0, X_1, \ldots, X_n, M_1, \ldots, M_{n+1})$, $n \in \N$.
For each $k \in \Nset$, we obtain
 \begin{equation*}
  \begin{aligned}
   X_k
   &= B_k + \sum_{i=1}^\infty \theta_k^{(k-i)} \circ \cdots \circ \theta_{k-i+1}^{(k-i)} \circ B_{k-i} \\
   &= B_k + \sum_{i=1}^{k-1} \theta_k^{(k-i)} \circ \cdots \circ \theta_{k-i+1}^{(k-i)} \circ B_{k-i} \\
   &\quad
    + \theta_k^{(0)} \circ \cdots \circ \theta_1^{(0)} \circ B_0 + \sum_{i=k+1}^{\infty} \theta_k^{(k-i)} \circ \cdots \circ \theta_{k-i+1}^{(k-i)} \circ B_{k-i} \\
   &\eind \kappa_k + \theta_k^{(0)} \circ \cdots \circ \theta_1^{(0)} \circ \biggl(B_0 + \sum_{j=1}^\infty \theta_0^{(-j)} \circ \cdots \circ
            \theta_{-j+1}^{(-j)} \circ B_{-j}\biggr) \\
   &= \kappa_k + \theta_k^{(0)} \circ \cdots \circ \theta_1^{(0)} \circ X_0 ,
  \end{aligned}
 \end{equation*}
 where $\kappa_k := B_k + \sum_{i=1}^{k-1} \theta_k^{(k-i)} \circ \cdots \circ \theta_{k-i+1}^{(k-i)} \circ B_{k-i}$,
 since for each $k\in\NN$, $\theta^{(0)}_k\circ\cdots\circ\theta^{(0)}_1\circ\theta^{(-j)}_0\circ\cdots\circ\theta^{(-j)}_{-j+1}\circ B_{-j}$, $j\geq 1$,
 has the same distribution as $\theta^{(k-i)}_k\circ\cdots \circ \theta^{(k-i)}_{k-i+1}\circ B_{k-i}$, $i\geq k+1$,
 and $\theta_k^{(0)} \circ \cdots \circ \theta_1^{(0)} \circ B_0$, $\theta^{(0)}_k\circ\cdots\circ\theta^{(0)}_1\circ\theta^{(-j)}_0\circ\cdots\circ\theta^{(-j)}_{-j+1}\circ B_{-j}$, $j\geq 1$, and $\kappa_k$ are independent.
Note also that $\kappa_k$, $\theta_k^{(0)} \circ \cdots \circ \theta_1^{(0)}$ and $X_0$ are independent for any $k\in\N$ (in the sense given in Lemma \ref{lem:reprX}).
In the same way, we get
 \begin{equation}\label{X0Xk}
  (X_0, X_1, \ldots, X_k)
  \eind
  (X_0, \kappa_1 + \theta_1^{(0)} \circ X_0, \ldots,
   \kappa_k + \theta_k^{(0)} \circ \cdots \circ \theta_1^{(0)} \circ X_0) .
 \end{equation}
Moreover, for each $k \in\NN$, using \eqref{X0Xk} for $(X_0, X_1, \ldots, X_{k+1})$, we obtain
 \begin{align*}
  &M_{k+1}
   = X_{k+1} - \mu_A X_k - \mu_B \\
  &\eind
   \kappa_{k+1} + \theta_{k+1}^{(0)} \circ \cdots \circ \theta_1^{(0)} \circ X_0
   - \mu_A (\kappa_k + \theta_k^{(0)} \circ \cdots \circ \theta_1^{(0)} \circ X_0)
   - \mu_B \\
  &= B_{k+1}
     + \sum_{i=1}^k \theta_{k+1}^{(k+1-i)} \circ \cdots \circ \theta_{k+2-i}^{(k+1-i)} \circ B_{k+1-i}
     + \theta_{k+1}^{(0)} \circ \cdots \circ \theta_1^{(0)} \circ X_0 \\
  &\quad
     - \mu_A \biggl(B_k
     + \sum_{i=1}^{k-1} \theta_k^{(k-i)} \circ \cdots \circ \theta_{k+1-i}^{(k-i)} \circ B_{k-i} \biggr)
     - \mu_A \bigl(\theta_k^{(0)} \circ \cdots \circ \theta_1^{(0)} \circ X_0\bigr)
     - \mu_B \\
  &= B_{k+1} - \mu_B
     + \theta_{k+1}^{(k)} \circ B_k - \mu_A B_k \\
  &\quad
     + \sum_{i=2}^k \theta_{k+1}^{(k+1-i)} \circ \cdots \circ \theta_{k+2-i}^{(k+1-i)} \circ B_{k+1-i}
     - \mu_A \sum_{i=1}^{k-1} \theta_k^{(k-i)} \circ \cdots \circ \theta_{k+1-i}^{(k-i)} \circ B_{k-i} \\
  &\quad
     + \theta_{k+1}^{(0)} \circ \cdots \circ \theta_1^{(0)} \circ X_0
     - \mu_A \bigl(\theta_k^{(0)} \circ \cdots \circ \theta_1^{(0)} \circ X_0\bigr) \\
  &= B_{k+1} - \mu_B + \sum_{i=0}^{k-1} \ttheta_{k+1}^{(k-i)} \circ \theta_k^{(k-i)} \circ \cdots \circ \theta_{k+1-i}^{(k-i)} \circ B_{k-i} +
  \ttheta_{k+1}^{(0)} \circ \theta_k^{(0)} \circ \cdots \circ \theta_1^{(0)} \circ X_0 \\
  &= \tkappa_{k+1} + \ttheta_{k+1}^{(0)} \circ \theta_k^{(0)} \circ \cdots \circ \theta_1^{(0)} \circ X_0 ,
 \end{align*}
 where $\ttheta_k^{(\ell)}$, $k, \ell \in \Zset$, are given by
 \[
   \ttheta_k^{(\ell)} \circ i
   := \begin{cases}
       \theta^{(\ell)}_k\circ i - i\mu_A =\sum_{j=1}^i (A_{k,j}^{(\ell)} - \mu_A) , & \text{for \ $i \in \Nset$,} \\
       0 , & \text{for \ $i = 0$,}
      \end{cases}
 \]
 and $\tkappa_{k+1} := B_{k+1} - \mu_B + \sum_{i=0}^{k-1} \ttheta_{k+1}^{(k-i)} \circ \theta_k^{(k-i)} \circ \cdots \circ \theta_{k+1-i}^{(k-i)} \circ B_{k-i}$.
Note that
 \begin{align*}
  \tkappa_{k+1}
     & = B_{k+1} - \mu_B +   \ttheta_{k+1}^{(k)} \circ B_k
           + \sum_{i=1}^{k-1} \ttheta_{k+1}^{(k-i)} \circ \theta_k^{(k-i)} \circ \cdots \circ \theta_{k+1-i}^{(k-i)} \circ B_{k-i}\\
     & \eind B_{k+1} - \mu_B + \ttheta_{k+1}^{(k)} \circ \kappa_k,
 \end{align*}
 since $\ttheta_{k+1}^{(k-i)}\circ j$, $i \in \{0,\ldots,k-1\}$, $j\geq 0$, are independent having the same distribution as $\ttheta_{k+1}^{(k)}\circ j$,
 $j\geq 0$, and $\ttheta_{k+1}^{(k)}$ is independent of $B_{k+1}$.
Further, $\tkappa_{k+1}$, \ $\ttheta_{k+1}^{(0)} \circ \theta_k^{(0)} \circ \cdots \circ \theta_1^{(0)}$ and $X_0$ are independent
 (in the sense given in Lemma \ref{lem:reprX} ).
Moreover, we have $M_1 = \ttheta_1^{(0)} \circ X_0 + B_1 - \mu_B = \tkappa_1 + \ttheta_1^{(0)} \circ X_0$ with $\tkappa_1 := B_1 - \mu_B$.
In the same way, we get
 \begin{equation}\label{X0XnM0Mn+1}
  \begin{aligned}
   &(X_0, X_1, \ldots, X_n, M_1, \ldots, M_{n+1}) \\
   &\eind
   (X_0, \kappa_1 + \theta_1^{(0)} \circ X_0, \ldots,
   \kappa_n + \theta_n^{(0)} \circ \cdots \circ \theta_1^{(0)} \circ X_0, \\
   &\phantom{\eind\,(}
    \tkappa_1 + \ttheta_1^{(0)} \circ X_0, \ldots, \tkappa_{n+1} + \ttheta_{n+1}^{(0)} \circ \theta_n^{(0)} \circ \cdots \circ \theta_1^{(0)} \circ X_0) ,
    \qquad n\geq 0.
  \end{aligned}
 \end{equation}

{\bf Step 1.}
First, we check that
 \begin{align}\label{help4_2}
  \cL\biggl(\frac{X_0^{3/2}}{x} \,\bigg|\, (X_0^{3/2} \lor X_0 |M_1|) > x\biggr)
  \wto \cL(\tY)
  \qquad \text{as \ $x \to \infty$.}
 \end{align}
For each $x, y \in (0, \infty)$, we have
 \[
   \pr\biggl(\frac{X_0^{3/2}}{x} > y \,\bigg|\, (X_0^{3/2} \lor X_0 |M_1|) > x\biggr)
   = \frac{P_1(x,y)}{Q_1(x)}
 \]
 with
 \[
   P_1(x,y)
   := \frac{\pr(X_0^{3/2}>xy,\,(X_0^{3/2}\lor X_0|M_1|)>x)}
           {\pr(X_0^{3/2}>x)} , \qquad
   Q_1(x)
   := \frac{\pr((X_0^{3/2}\lor X_0|M_1|)>x)}
           {\pr(X_0^{3/2}>x)} .
 \]
For each $x \in (0, \infty)$ and $c \in (0, 1)$, we can write $Q_1(x) = Q_{1,1}(x,c) + Q_{1,2}(x,c)$ with
 \begin{gather*}
  Q_{1,1}(x,c)
  := \frac{\pr((X_0^{3/2}\lor X_0|M_1|)>x,\,X_0^{3/2}>cx)}
          {\pr(X_0^{3/2}>x)} , \\
  Q_{1,2}(x,c)
  := \frac{\pr((X_0^{3/2}\lor X_0|M_1|)>x,\,X_0^{3/2}\leq cx)}
          {\pr(X_0^{3/2}>x)} ,
 \end{gather*}
 where $Q_{1,1}(x,c) = Q_{1,1,1}(x,c) Q_{1,1,2}(x,c)$ with
 \[
   Q_{1,1,1}(x,c)
   := \frac{\pr((X_0^{3/2}\lor X_0|M_1|)>x,\,X_0^{3/2}>cx)}
           {\pr(X_0^{3/2}>cx)} , \qquad
   Q_{1,1,2}(x,c)
   := \frac{\pr(X_0^{3/2}>cx)}
           {\pr(X_0^{3/2}>x)} .
 \]
For each $c \in (0, 1)$, we have
 \begin{equation}\label{Q111}
  \lim_{x\to\infty} Q_{1,1,2}(x,c)
  = \lim_{x\to\infty}
     \frac{\pr(X_0>c^{2/3}x^{2/3})}
          {\pr(X_0>x^{2/3})}
  = c^{-2\alpha/3} ,
 \end{equation}
 since $X_0$ is regularly varying with tail index $\alpha$.
For each $c \in (0, 1)$, using $X_0^{3/2} \lor X_0 |M_1| = X_0^{3/2} \bigl(1 \lor \frac{|M_1|}{\sqrt{X_0}}\bigr)$, if $X_0>0$,
 and that \eqref{eq:pseudoTail} yields $\cL(x^{-1}X_0, \frac{M_1}{\sqrt{X_0}} \mid X_0>x) \wto \cL(Y_0,Z_0)$ as $x\to\infty$
 (since $X_0\vee 1 = X_0$ if $X_0>1$), by Lemma \ref{lemma:cmt}, we obtain
 \begin{align*}
  Q_{1,1,1}(x,c)
  &= \pr((X_0^{3/2} \lor X_0 |M_1|) > x \mid X_0^{3/2} > c x) \\
  &= \pr\biggl(\biggl(\frac{X_0}{(cx)^{2/3}}\biggr)^{3/2}
               \biggl(1 \lor \frac{|M_1|}{\sqrt{X_0}}\biggr)
               > \frac{1}{c}
               \,\bigg|\,
               X_0 > (cx)^{2/3}\biggr) \\
  &\to \pr(Y_0^{3/2} (1 \lor |Z_0|) > c^{-1})
   \qquad \text{as $x \to \infty$,}
 \end{align*}
 where $Y_0$ is a Pareto distributed random variable such that $\pr(Y_0\geq y) = y^{-\alpha}$, $y\geq 1$.
Consequently, for each $c \in (0, 1)$, we have
 \begin{equation*}
  Q_{1,1}(x,c)
  \to c^{-2\alpha/3}
      \pr(Y_0^{3/2} (1 \lor |Z_0|) > c^{-1})
  \qquad \text{as $x \to \infty$.}
 \end{equation*}
For each $c \in (0, 1)$, by the tower rule and using $\pr(Y_0 > y) =
y^{-\alpha} \land 1$, $y \in (0, \infty)$, we have
 \begin{align}\label{help7}
 \begin{split}
  &c^{-2\alpha/3}
  \pr(Y_0^{3/2} (1 \lor |Z_0|) > c^{-1}) \\
  &= c^{-2\alpha/3}
     \ex(\pr(Y_0^{3/2} (1 \lor |Z_0|) > c^{-1} \mid Z_0)) \\
  &= c^{-2\alpha/3}
     \ex(\pr(Y_0 > c^{-2/3} (1 \lor |Z_0|)^{-2/3} \mid Z_0)) \\
  &= c^{-2\alpha/3}
     \ex\bigl(\bigl(c^{-2/3} (1 \lor |Z_0|)^{-2/3}\bigr)^{-\alpha} \land 1\bigr) \\
  &= \ex\bigl((1 \lor |Z_0|)^{2\alpha/3} \land c^{-2\alpha/3}\bigr) ,
  \end{split}
 \end{align}
 hence we have $\lim_{c\downarrow 0} \lim_{x \to \infty} Q_{1,1}(x,c) = \ex\bigl((1 \lor |Z_0|)^{2\alpha/3}\bigr)$.
Moreover, for each $x \in (0, \infty)$ and $c \in (0, 1)$, we have $Q_{1,2}(x,c) = Q_{1,2,1}(x,c) Q_{1,1,2}(x,c)$ with
 \begin{align*}
   Q_{1,2,1}(x,c)
   &:= \frac{\pr( (X_0^{3/2}\vee X_0|M_1|)>x,\,X_0^{3/2}\leq cx)}
           {\pr(X_0^{3/2}>cx)}
    = \frac{\pr(X_0|M_1|>x,\,X_0^{3/2}\leq cx)}
           {\pr(X_0^{3/2}>cx)} \\
   & = \frac{\pr(|M_1|X_0\1{\{X_0^{3/2}\leq cx\}}>x)}
          {\pr(X_0^{3/2}>cx)} .
 \end{align*}
By \eqref{X0XnM0Mn+1}, we have $(X_0, M_1) \eind (X_0, \tkappa_1 + \ttheta_1^{(0)} \circ X_0)$.
Hence, for each $x \in (0, \infty)$, $c \in (0, 1)$ and $\delta \in \bigl(0, \frac{\alpha}{3}\bigr)$,
 by Markov's inequality and the independence of $\tkappa_1$ and $X_0$,
 \begin{align*}
  &\pr\bigl(|M_1| X_0 \1{\{X_0^{3/2} \leq cx\}} > x\bigr) \\
  &\leq
   \pr\Bigl(|\tkappa_1| X_0 \1{\{X_0^{3/2} \leq cx\}} > \frac{x}{2}\Bigr)
   + \pr\Bigl(|\ttheta_1^{(0)} \circ X_0| X_0 \1{\{X_0^{3/2} \leq cx\}} > \frac{x}{2}\Bigr) \\
  &\leq \biggl(\frac{2}{x}\biggr)^{\alpha-\delta}
        \ex(|\tkappa_1|^{\alpha-\delta})
        \ex\bigl(X_0^{\alpha-\delta} \1{\{X_0^{3/2}\leq cx\}}\bigr)
        + \biggl(\frac{2}{x}\biggr)^2
          \ex\bigl((\ttheta_1^{(0)} \circ X_0)^2 X_0^2 \1{\{X_0^{3/2} \leq cx\}}\bigr) .
 \end{align*}
We have $\ex(|\tkappa_1|^{\alpha-\delta}) = \ex(|B - \mu_B|^{\alpha-\delta}) < \infty$, since $\vert B-\mu_B\vert$ is regularly varying
 with tail index $\alpha$, see \eqref{eq:tail_btilde}.
Moreover, $\ex\bigl(X_0^{\alpha-\delta} \1{\{X_0^{3/2}\leq cx\}}\bigr) \leq \ex(X_0^{\alpha-\delta}) < \infty$, since $X_0$ is regularly varying
 with tail index $\alpha$.
Further,
 \begin{align*}
  \ex\bigl({(\ttheta_1^{(0)} \circ X_0)^2} X_0^2 \1{\{X_0^{3/2} \leq cx\}}\bigr)
  = \ex\bigl(\ex\bigl({(\ttheta_1^{(0)} \circ X_0)^2} X_0^2 \1{\{X_0^{3/2} \leq cx\}} \, \big|\, X_0\bigr)\bigr)
  = \sigma_A^2 \ex\bigl(X_0^3 \1{\{X_0^{3/2} \leq cx\}}\bigr) .
 \end{align*}
Consequently, for
each $x \in (0, \infty)$, $c \in (0, 1)$ and $\delta \in \bigl(0, \frac{\alpha}{3}\bigr)$,
 \[
   Q_{1,2,1}(x,c)
   \leq \frac{2^{\alpha-\delta}\ex(|B - \mu_B|^{\alpha-\delta})\ex(X_0^{\alpha-\delta})}
             {x^{\alpha-\delta}\pr(X_0^{3/2}>cx)}
        + 4 \sigma_A^2
          \frac{\ex(X_0^3\1{\{X_0^{3/2}\leq cx\}})}
               {x^2\pr(X_0^{3/2}>cx)} .
 \]
The random variable $X_0^{3/2}$ is regularly varying with tail index $\frac{2\alpha}{3}$, since $X_0$ is regularly varying with tail index $\alpha$.
Hence $\alpha - \delta > \frac{2\alpha}{3}$ yields $x^{\alpha-\delta} \pr(X_0^{3/2} > c x) \to \infty$ as $x \to \infty$
 (see, e.g., Bingham et al.~\cite[Proposition 1.3.6. (v)]{bingham:goldie:teugels:1987}).
Applying Karamata's theorem (see, Lemma \ref{truncated_moments}), we obtain
 \[
   \frac{\ex(X_0^3\1{\{X_0^{3/2}\leq cx\}})}
        {x^2\pr(X_0^{3/2}>cx)}
   = c^2
     \frac{\ex((X_0^{3/2})^2\1{\{X_0^{3/2}\leq cx\}})}
          {(cx)^2\pr(X_0^{3/2}>cx)}
   \to c^2 \frac{\frac{2\alpha}{3}}
                {2-\frac{2\alpha}{3}}
   = \frac{c^2\alpha}{3-\alpha}
 \]
 as \ $x \to \infty$.
Consequently, by \eqref{Q111}, for each $c \in (0, 1)$, we obtain
 \begin{align*}
  \limsup_{x\to\infty}
   Q_{1,2}(x,c)
  \leq 4 c^{-2\alpha/3} \sigma_A^2
       \frac{c^2\alpha}{3-\alpha}
  = \frac{4\alpha\sigma_A^2}{3-\alpha} c^{2(3-\alpha)/3} ,
 \end{align*}
 hence
 \begin{align}\label{help8}
  \lim_{c\downarrow 0} \liminf_{x \to \infty} Q_{1,2}(x,c) =\lim_{c\downarrow 0} \limsup_{x \to \infty} Q_{1,2}(x,c) = 0.
 \end{align}
Summarizing, we get
 \begin{equation}\label{Q1}
  Q_1(x) \to \ex\bigl((1 \lor |Z_0|)^{2\alpha/3}\bigr) \qquad \text{as $x \to \infty$.}
 \end{equation}

Now we consider the term $P_1(x,y)$, $x,y \in (0, \infty)$.
For each $x,y \in (0, \infty)$, we have $P_1(x,y)= P_{1,1}(x,y)P_{1,2}(x,y)$, where
 \begin{align*}
  &P_{1,1}(x,y):= \frac{\pr(X_0^{3/2}>xy)}{\pr(X_0^{3/2}>x)},\\
  &P_{1,2}(x,y):=
     \frac{\pr(X_0^{3/2}>xy,\,(X_0^{3/2}\lor X_0|M_1|)>x)}
          {\pr(X_0^{3/2}>xy)}.
 \end{align*}
Since $X_0$ is regularly varying with tail index $\alpha$, for each $y\in(0,\infty)$, we have $\lim_{x\to\infty} P_{1,1}(x,y) = y^{-2\alpha/3}$.
Further, for each $y\in(0,\infty)$, using that \eqref{eq:pseudoTail} yields $\cL(x^{-1}X_0, \frac{M_1}{\sqrt{X_0}} \mid X_0>x) \wto \cL(Y_0,Z_0)$ as
 $x\to\infty$ (since $X_0\vee 1 = X_0$ if $X_0>1$), by Lemma \ref{lemma:cmt}, we have
 \begin{align*}
  P_{1,2}(x,y)
  &= \pr((X_0^{3/2} \lor X_0 |M_1|) > x \mid X_0^{3/2} > xy) \\
  &= \pr\biggl( \biggl(\frac{X_0}{(xy)^{2/3}}\biggr)^{3/2}
               \biggl(1 \lor \frac{|M_1|}{\sqrt{X_0}}\biggr)
               > \frac{1}{y}
               \,\bigg|\,
               X_0 > (xy)^{2/3}\biggr) \\
  &\to \pr(Y_0^{3/2} (1 \lor |Z_0|) > y^{-1})
  \qquad \text{as \ $x \to \infty$.}
 \end{align*}
Consequently, for each $y \in (0, \infty)$, by the tower rule and using $\pr(Y_0>y) = y^{-\alpha}\wedge 1$, $y\in(0,\infty)$,
 we have
 \begin{align*}
  P_1(x,y)
  & \to y^{-2\alpha/3}\pr(Y_0^{3/2} (1 \lor |Z_0|) > y^{-1})
    = y^{-2\alpha/3} \pr(Y_0 > (y(1 \lor |Z_0|))^{-2/3}) \\
  & = y^{-2\alpha/3} \ex[\pr( Y_0 > (y (1 \lor |Z_0|))^{-2/3} \mid Z_0) ]
    = y^{-2\alpha/3} \ex[ ( y^{2\alpha/3} (1 \lor |Z_0|)^{2\alpha/3} ) \wedge 1 ] \\
  & = \ex[ (1 \lor |Z_0|)^{2\alpha/3} \wedge y^{-2\alpha/3} ]
    = \ex[  (y \vee (1 \lor |Z_0|)^{-1})^{-2\alpha/3} ]
  \qquad \text{as $x \to \infty$.}
 \end{align*}
By \eqref{Q1}, for each $y \in (0, \infty)$, we obtain
 \[
   \pr\biggl(\frac{X_0^{3/2}}{x} > y \,\bigg|\, (X_0^{3/2} \lor X_0 |M_1|) > x\biggr)
   \to \frac{\ex\bigl((y \lor (1 \lor |Z_0|)^{-1})^{-2\alpha/3}\bigr)}
            {\ex\bigl((1 \lor |Z_0|)^{2\alpha/3}\bigr)}
 \]
 as $x \to \infty$, thus we conclude \eqref{help4_2}.

{\bf Step 2.}
We check that for each $k \in \Nset$,
 \begin{align}\label{help5_2}
  \cL\biggl(\biggl(\frac{X_k}{X_0}\biggr)^{3/2} \,\bigg|\, (X_0^{3/2} \lor X_0 |M_1|) > x\biggr)
  \wto \delta_{\mu_A^{3k/2}}
  \qquad \text{as $x \to \infty$.}
 \end{align}
For each $k \in \Nset$, $x, y \in (0, \infty)$, we have
 \[
   \pr\biggl(\biggl(\frac{X_k}{X_0}\biggr)^{3/2} > y \,\bigg|\, (X_0^{3/2} \lor X_0 |M_1|) > x\biggr)
   = \frac{P_2(x,y)}{Q_1(x)}
 \]
 with
 \[
   P_2(x,y)
   := \frac{\pr\bigl(\bigl(\frac{X_k}{X_0}\bigr)^{3/2}>y,\,(X_0^{3/2}\lor X_0|M_1|)>x\bigr)}
          {\pr(X_0^{3/2}>x)} .
 \]
For each $k \in \Nset$, $x, y \in (0, \infty)$, and $c \in (0, 1)$, we can write $P_2(x,y) = P_{2,1}(x,y,c) + P_{2,2}(x,y,c)$ with
 \begin{gather*}
  P_{2,1}(x,y,c)
  := \frac{\pr\bigl(\bigl(\frac{X_k}{X_0}\bigr)^{3/2}>y,\,(X_0^{3/2}\lor X_0|M_1|)>x,\,X_0^{3/2}>cx\bigr)}
          {\pr(X_0^{3/2}>x)} , \\
  P_{2,2}(x,y,c)
  :=  \frac{\pr\bigl(\bigl(\frac{X_k}{X_0}\bigr)^{3/2}>y,\,(X_0^{3/2}\lor X_0|M_1|)>x,\,X_0^{3/2}\leq cx\bigr)}
          {\pr(X_0^{3/2}>x)} ,
 \end{gather*}
 where $P_{2,1}(x,y,c) = P_{2,1,1}(x,y,c) Q_{1,1,2}(x,c)$ with
 \[
   P_{2,1,1}(x,y,c)
   := \frac{\pr\bigl(\bigl(\frac{X_k}{X_0}\bigr)^{3/2}>y,\,(X_0^{3/2}\lor X_0|M_1|)>x,\,X_0^{3/2}>cx\bigr)}
          {\pr(X_0^{3/2}>cx)} .
 \]
For each $k \in \Nset$, $y \in (0,\infty) \setminus \{\mu_A^{3k/2}\}$ and $c \in (0, 1)$,
 using \eqref{tail_process}, \eqref{eq:pseudoTail} and Lemma \ref{lemma:cmt}, we obtain
 \begin{align*}
  &P_{2,1,1}(x,y,c)
   = \pr\biggl(\biggl(\frac{X_k}{X_0}\biggr)^{3/2} > y, \, (X_0^{3/2} \lor X_0 |M_1|) > x \,\bigg|\, X_0^{3/2} > c x\biggr) \\
  &= \pr\biggl(\biggl(\frac{(cx)^{-2/3}X_k}{(cx)^{-2/3}X_0}\biggr)^{3/2} > y, \ \biggl(\frac{X_0}{(cx)^{2/3}}\biggr)^{3/2}
               \biggl(1 \lor \frac{|M_1|}{\sqrt{X_0}}\biggr)
               > \frac{1}{c}
               \,\bigg|\,
               X_0 > (cx)^{2/3}\biggr) \\
  &\to \pr(\mu_A^{3k/2} > y, \,
  Y_0^{3/2} (1 \lor |Z_0|) > c^{-1})
  \qquad \text{as $x \to \infty$.}
 \end{align*}
Consequently, by \eqref{Q111} and \eqref{help7}, for each $k \in \Nset$, $y \in (0,\infty) \setminus \{\mu_A^{3k/2}\}$ and $c \in (0, 1)$, we have
 \begin{equation*}
  \begin{aligned}
   P_{2,1}(x,y,c)
   &\to c^{-2\alpha/3}
        \1{(y,\infty)}(\mu_A^{3k/2})
        \pr(Y_0^{3/2} (1 \lor |Z_0|) > c^{-1}) \\
   &= \1{(y,\infty)}(\mu_A^{3k/2})
      \ex\bigl((1 \lor |Z_0|)^{2\alpha/3} \wedge c^{-2\alpha/3}\bigr)
   \qquad \text{as $x \to \infty$,}
  \end{aligned}
 \end{equation*}
 hence, for each $k\in\NN$ and $y \in (0,\infty) \setminus \{\mu_A^{3k/2}\}$, we have
 \[
   \lim_{c\downarrow 0} \lim_{x\to\infty} P_{2,1}(x,y,c)
   = \1{(y,\infty)}(\mu_A^{3k/2}) \ex\bigl((1 \lor |Z_0|)^{2\alpha/3}\bigr) .
 \]
Moreover, for each $x, y \in (0, \infty)$ and $c \in (0, 1)$, we have $P_{2,2}(x,y,c) \leq Q_{1,2}(x,c)$,
 hence, by \eqref{help8}, $\lim_{c\downarrow 0} \liminf_{x\to\infty} P_{2,2}(x,y,c)=\lim_{c\downarrow 0} \limsup_{x\to\infty} P_{2,2}(x,y,c) = 0$.
Summarizing, for each $k \in \Nset$ and $y \in (0,\infty) \setminus \{\mu_A^{3k/2}\}$, we get
 \[
   P_2(x,y) \to \1{(y,\infty)}(\mu_A^{3k/2}) \ex\bigl((1 \lor |Z_0|)^{2\alpha/3}\bigr) \qquad \text{as $x \to \infty$.}
 \]
By \eqref{Q1}, for each $k \in \Nset$ and $y \in (0,\infty) \setminus \{\mu_A^{3k/2}\}$, we obtain
 \[
   \pr\biggl(\biggl(\frac{X_k}{X_0}\biggr)^{3/2} > y \,\bigg|\, (X_0^{3/2} \lor X_0 |M_1|) > x\biggr)
   \to \1{(y,\infty)}(\mu_A^{3k/2})
   \qquad \text{as $x \to \infty$,}
 \]
 thus we conclude \eqref{help5_2}.

{\bf Step 3.}
We check that for each $n \geq 0$,
 \begin{align}\label{help5+_2}
  \cL\biggl(\frac{X_0^{3/2}}{x}, \, \frac{M_1}{\sqrt{X_0}}, \ldots, \frac{M_{n+1}}{\sqrt{X_0}} \,\bigg|\,  (X_0^{3/2} \lor X_0 |M_1|) > x\biggr)
  \dto \cL(\tY, \tZ_0, \mu_A^{1/2} \tZ_1,\ldots, \mu_A^{n/2} \tZ_n)
 \end{align}
 as $x \to \infty$.
For each $n \geq 0$, $x, y \in (0, \infty)$ and $v_0, \ldots, v_n \in \Rset$, we have
 \[
   \pr\biggl(\frac{X_0^{3/2}}{x} > y, \frac{M_1}{\sqrt{X_0}} > v_0, \ldots, \frac{M_{n+1}}{\sqrt{X_0}} > v_n \,\bigg|\, (X_0^{3/2} \lor X_0 |M_1|) > x\biggr)
   = \frac{P_3(x,y,v_0,\ldots,v_n)}{Q_1(x)}
 \]
 with
 \[
   P_3(x,y,v_0,\ldots,v_n)
   := \frac{\pr\bigl(\frac{X_0^{3/2}}{x}>y,\frac{M_1}{\sqrt{X_0}}>v_0,\ldots,\frac{M_{n+1}}{\sqrt{X_0}}>v_n,\,(X_0^{3/2}\lor X_0|M_1|)>x\bigr)}
          {\pr(X_0^{3/2}>x)} .
 \]
For each $n \geq 0$, $x,y \in (0, \infty)$, and $v_0, \ldots, v_n \in \Rset$, we can write $P_3(x,y,v_0,\ldots,v_n)
 = P_{3,1}(x,y,v_0,\ldots,v_n) P_{3,2}(x,y)$ with
 \begin{align*}
  P_{3,1}(x,y,v_0,\ldots,v_n)
  &:= \frac{\pr\Bigl(\frac{X_0^{3/2}}{x}>y,\frac{M_1}{\sqrt{X_0}}>v_0,\ldots,\frac{M_{n+1}}{\sqrt{X_0}}>v_n,\,(X_0^{3/2}\lor X_0|M_1|)>x\Bigr)}
          {\pr(X_0^{3/2}>yx)} \\
  &= \pr\Bigl( (X_0^{3/2}\lor X_0|M_1|)>x,  \frac{M_1}{\sqrt{X_0}}>v_0,\ldots,\frac{M_{n+1}}{\sqrt{X_0}}>v_n \mid X_0^{3/2}>yx \Bigr) ,\\
  P_{3,2}(x,y)
    &:=  \frac{\pr(X_0^{3/2}>yx)}{\pr(X_0^{3/2}>x)}.
 \end{align*}
Since $X_0^{3/2}$ is regularly varying with tail index $2\alpha/3$, we have for each $y \in (0, \infty)$,
 \[
   \lim_{x\to\infty} P_{3,2}(x,y)  = y^{-2\alpha/3}.
 \]
Further, using that \eqref{eq:pseudoTail} holds with $W_i'$ replaced by $W_i$ (see the proof of Proposition \ref{prop:pseudoTail}),
 by Lemma \ref{lemma:cmt}, we have
 \begin{align*}
  &P_{3,1}(x,y,v_0,\ldots,v_n) \\
  & = \pr\biggl( \biggl(\frac{X_0}{(yx)^{2/3}}\biggr)^{3/2} \biggl(1 \lor \frac{|M_1|}{\sqrt{X_0}}\biggr) > \frac{1}{y}, \\
  &\phantom{= \pr\biggl(}
                 \frac{M_1}{\sqrt{X_0}} > v_0, \frac{M_2}{\sqrt{X_0\mu_A}} > \frac{v_1}{\sqrt{\mu_A}}, \ldots, \frac{M_{n+1}}{\sqrt{X_0 \mu_A^n}}
                  > \frac{v_n}{\sqrt{\mu_A^n}}
               \,\bigg|\,
               X_0 > (yx)^{2/3}\biggr) \\
  &\to \pr(Y_0^{3/2} (1 \lor |Z_0|) > y^{-1}, Z_0 > v_0, Z_1 > \mu_A^{-1/2} v_1, \ldots, Z_n > \mu_A^{-n/2}v_n) \\
  &=\pr(Y_0^{3/2} (1 \lor |Z_0|) > y^{-1}, Z_0 > v_0) \pr(Z_1 > \mu_A^{-1/2} v_1, \ldots, Z_n > \mu_A^{-n/2}v_n)
 \end{align*}
 as $x \to \infty$.
Consequently, for each $n\geq 0$, $y\in(0,\infty)$ and $v_0, \ldots, v_n \in \Rset$, we get
 \begin{align*}
  &\lim_{x\to\infty}P_3(x,y,v_0,\ldots,v_n) \\
  &= y^{-2\alpha/3}
     \pr(Y_0^{3/2} (1 \lor |Z_0|) > y^{-1}, Z_0 > v_0) \pr(Z_1 > \mu_A^{-1/2} v_1, \ldots, Z_n > \mu_A^{-n/2}v_n) \\
  &= y^{-2\alpha/3}
     \ex\big[ \pr(Y_0^{3/2} (1 \lor |Z_0|) > y^{-1}, Z_0 > v_0 \mid Z_0 ) \big]
      \pr(Z_1 > \mu_A^{-1/2} v_1, \ldots, Z_n > \mu_A^{-n/2}v_n) \\
  &= y^{-2\alpha/3}
     \ex\big[ \1{\{Z_0 > v_0\}}\pr(Y_0 > (y(1 \lor |Z_0|))^{-2/3} \mid Z_0 ) \big]
      \pr(Z_1 > \mu_A^{-1/2} v_1, \ldots, Z_n > \mu_A^{-n/2}v_n) \\
  & = y^{-2\alpha/3}
      \ex\big[ \1{\{Z_0 > v_0\}} ( (y(1 \lor |Z_0|))^{2\alpha/3}\wedge 1 ) \big]
      \pr(Z_1 > \mu_A^{-1/2} v_1, \ldots, Z_n > \mu_A^{-n/2}v_n) \\
  & = \ex\big[ \1{\{Z_0 > v_0\}} ( (1 \lor |Z_0|)^{2\alpha/3}\wedge y^{-2\alpha/3} ) \big]
     \pr(Z_1 > \mu_A^{-1/2} v_1, \ldots, Z_n > \mu_A^{-n/2}v_n) \\
  & = \ex\bigl(\bigl(y \lor (1 \lor |Z_0|)^{-1}\bigr)^{-2\alpha/3}\1{(v_0,\infty)}(Z_0)\bigr)
        \pr(\mu_A^{1/2} Z_1 > v_1, \ldots, \mu_A^{n/2} Z_n > v_n).
 \end{align*}
By \eqref{Q1}, for each $n\geq 0$, $y\in(0,\infty)$ and $v_0, \ldots, v_n \in \Rset$, we obtain
 \begin{align*}
  &\pr\biggl(\frac{X_0^{3/2}}{x} > y, \frac{M_1}{\sqrt{X_0}} > v_0, \ldots, \frac{M_{n+1}}{\sqrt{X_0}} > v_n \,\bigg|\, (X_0^{3/2} \lor X_0 |M_1|) > x\biggr) \\
  &\to \frac{\ex\bigl(\bigl(y \lor (1 \lor |Z_0|)^{-1}\bigr)^{-2\alpha/3}\1{(v_0,\infty)}(Z_0)\bigr)}
            {\ex\bigl((1 \lor |Z_0|)^{2\alpha/3}\bigr)}
       \pr(\mu_A^{1/2} Z_1 > v_1, \ldots, \mu_A^{n/2} Z_n > v_n) \\
  &= \pr(\tY > y, \, \tZ_0 > v_0)
     \pr\bigl(\mu_A^{1/2} \tZ_1 > v_1, \ldots, \mu_A^{n/2} \tZ_n > v_n\bigr) \\
  &= \pr\bigl(\tY > y, \, \tZ_0 > v_0, \, \mu_A^{1/2} \tZ_1 > v_1, \ldots, \mu_A^{n/2} \tZ_n > v_n\bigr)
 \end{align*}
 as $x \to \infty$, thus we conclude \eqref{help5+_2}.

{\bf Step 4.}
For all $n \geq 0$, we have
 \begin{align*}
  &\cL\biggl(\biggl(\frac{1}{x}
   (X_k^{3/2}, X_k M_{k+1})\biggr)_{k\in\{0,1,\ldots,n\}} \,\bigg|\, (X_0^{3/2} \lor X_0 |M_1|) > x\biggr) \\
  &= \cL\biggl(\biggl(\frac{X_0^{3/2}}{x} \biggl(\frac{X_k}{X_0}\biggr)^{3/2}, \frac{X_0^{3/2}}{x} \frac{X_k}{X_0} \frac{M_{k+1}}{\sqrt{X_0}}\biggr)\biggr)_{k\in\{0,1,\ldots,n\}} \,\bigg|\, (X_0^{3/2} \lor X_0 |M_1|) > x\biggr)
 \end{align*}
 \begin{align*}
  &\dto
   \cL\bigl((\tY \mu_A^{3k/2}, \tY \mu_A^k \mu_A^{k/2} \tZ_k)_{k\in\{0,1,\ldots,n\}}\bigr)
  \qquad \text{as \ $x \to \infty$.}
 \end{align*}
Indeed, \eqref{help5_2}, Lemmas \ref{lemma:cmt} and \ref{lemma:Slutsky4} yield
 \[
   \cL\biggl(\biggl(\frac{X_k}{X_0}\biggr)_{k\in\{0,1,\ldots,n\}} \,\bigg|\, (X_0^{3/2} \lor X_0 |M_1|) > x\biggr)
   \wto \delta_{(1,\mu_A,\ldots,\mu_A^n)} \qquad \text{as $x \to \infty$,}
 \]
 and then, identifying $\R^{(2n+2)\times(2n+2)}$ with $\R^{(2n+2)^2}$ in a natural way, we can use again Lemma \ref{lemma:cmt} to obtain
 \begin{align*}
  &\cL\biggl(\diag_{2n+2}
    \biggl( 1, 1, \biggl(\frac{X_1}{X_0}\biggr)^{3/2}, \frac{X_1}{X_0},\ldots, \biggl(\frac{X_n}{X_0}\biggr)^{3/2}, \frac{X_n}{X_0} \biggr)
        \,\bigg|\, (X_0^{3/2} \lor X_0 |M_1|) > x\biggr) \\
  &\qquad \wto \delta_{  \diag_{2n+2}(1, 1,\mu_A^{3/2},\mu_A, \ldots,\mu_A^{3n/2},\mu_A^n) }\qquad \text{as $x \to \infty$.}
 \end{align*}
Next, \eqref{help5+_2} and the conditional version of the continuous mapping theorem (see Lemma \ref{lemma:cmt}) imply
 \begin{align*}
  &\cL\biggl(\biggl(\biggl(\frac{X_0}{x^{2/3}}\biggr)^{3/2}, \biggl(\frac{X_0}{x^{2/3}}\biggr)^{3/2} \frac{M_{k+1}}{\sqrt{X_0}}\biggr)\biggr)_{k\in\{0,1,\ldots,n\}} \,\bigg|\, (X_0^{3/2} \lor X_0 |M_1|) > x\biggr) \\
  &\wto \cL\bigl((\tY^{3/2}, \tY^{3/2} \mu_A^{k/2} \tZ_k)_{k\in\{0,1,\ldots,n\}}\bigr)
  \qquad \text{as $x \to \infty$.}
 \end{align*}
Finally, identifying $\Rset^{n+1} \times \Rset^{n+1}$ and $\Rset^{2(n+1)}$ in a natural way and applying Lemma \ref{lemma:Slutsky3}, we obtain the convergence statement of the theorem.

The jointly regularly varying property of $(X_k^{3/2}, X_k M_{k+1})_{k\in\Z}$ follows by Theorem 2.1 in Basrak and Segers \cite{basrak:segers:2009}.
The existence of a (whole) tail process of $(X_k^{3/2}, X_kM_{k+1})_{k\in\ZZ}$ follows by Basrak and Segers \cite [Theorem 2.1]{basrak:segers:2009}.
\proofend

In the next remark we point out that $\tY$  given in Theorem \ref{th:joint_reg2} is not a Pareto-distributed random variable.

\begin{remark}\label{tY}
Note that \ $\eqref{help6}$ \  readily yields that
 \[
   \pr(\tY > y)
   = \frac{\ex\bigl(\bigl(y \lor (1 \lor |Z_0|)^{-1}\bigr)^{-2\alpha/3}\bigr)}
            {\ex\bigl((1 \lor |Z_0|)^{2\alpha/3}\bigr)} , \qquad
   y \in \RR .
 \]
Consequently, $\pr(\tY \in(0,\infty)) = 1$ and
 \[
   \pr(\tY > y)
   \leq \frac{1}{\ex\bigl((1 \lor |Z_0|)^{2\alpha/3}\bigr)} y^{-2\alpha/3} \qquad
   \text{for \ $y \in (0,\infty)$,}
 \]
 and equality holds for $y \in [1, \infty)$.
Indeed, for each $y \in (0,\infty)$, we have $y \lor (1 \lor |Z_0|)^{-1} \geq y$ almost surely,
 thus $\bigl(y \lor (1 \lor |Z_0|)^{-1}\bigr)^{-2\alpha/3} \leq y^{-2\alpha/3}$
 almost surely, hence \ $\ex\bigl(\bigl(y \lor (1 \lor |Z_0|)^{-1}\bigr)^{-2\alpha/3}\bigr) \leq y^{-2\alpha/3}$,
 and for $y \in [1, \infty)$, we have $y \lor (1 \lor |Z_0|)^{-1} = y$ almost surely.
\proofend
\end{remark}

\section*{Acknowledgements}
 We are grateful to M\'arton Isp\'any for raising the problem in 2002 and to Peter Kern for calling attention to Corollary 2.1 in {\L}uczak \cite{Luc} on the existence of
 an infinitely differentiable density function for operator stable laws.
We would like to thank the referees for their comments that helped us to improve the paper.
This paper has been revised after the sudden death of Gyula Pap, the fourth author and the main instigator of this project, in October 2019.

\bibliographystyle{plain}
\bibliography{CLS}

\end{document}